\documentclass[letterpaper,11pt]{article}
\usepackage[
  left=0.56in,
  right=0.56in,
  top=0.56in,
  bottom=0.59in
]{geometry}

\usepackage[english]{babel}
\usepackage{chngcntr}
\usepackage[utf8]{inputenc}
\usepackage{amsmath,mathrsfs}
\usepackage{indentfirst}
\usepackage{fancyhdr}
\usepackage{amssymb}
\usepackage{amsthm}
\usepackage{graphicx}
\usepackage{color}
\usepackage{xcolor}
\usepackage{cancel}
\usepackage{subcaption}
\usepackage[dvipsnames]{xcolor}
\usepackage{eucal}
\usepackage{latexsym}
\usepackage{chngcntr}
\usepackage[expansion=false]{microtype}
\usepackage{newlfont}
\usepackage{enumitem}
\usepackage{geometry}
\usepackage{tikz}
\usepackage{todonotes}

\usepackage{amsmath, amsfonts, amsthm,amssymb}

\usepackage{tikz-cd}
\usepackage{multicol}
\usepackage{cancel}
\usepackage{mathtools}
\usepackage{graphicx}

\usepackage[colorlinks=true,    
            linkcolor=blue,    
            citecolor=blue,     
            urlcolor=blue]{hyperref}

\newcommand{\e}{\epsilon}
\newcommand{\cO}{\mathcal{O}}

\newcommand{\cw}{\mathtt{c}}

\newcommand{\pa}{\partial}

\newcommand{\N}{{\mathbb{N}}}
\newcommand{\R}{{\mathbb{R}}}

\newcommand{\T}{{\mathbb{T}}}

\newcommand{\vet}[2]{\begin{bmatrix}#1 \\ #2 \end{bmatrix}}
\newcommand{\pvet}[2]{\begin{pmatrix}#1 \\ #2 \end{pmatrix}}
\newcommand{\cB}{\mathcal{B}}

\newcommand{\la}{\langle}
\newcommand{\ra}{\rangle}

\newcommand{\im}{\textup{i}}

\newcommand{\vertiii}[1]{{\left\vert\kern-0.25ex\left\vert\kern-0.25ex\left\vert #1 
    \right\vert\kern-0.25ex\right\vert\kern-0.25ex\right\vert}}
\newcommand{\opnorm}[1]{{\vert\kern-0.25ex\vert\kern-0.25ex\vert #1 
    \vert\kern-0.25ex\vert\kern-0.25ex\vert}}

\newtheorem{theorem}{Theorem}[section]
\newtheorem{proposition}[theorem]{Proposition}

\newtheorem{lemma}[theorem]{Lemma}

\theoremstyle{remark}
\newtheorem{remark}[theorem]{Remark}
\newtheorem*{remark*}{Remark}
\newtheorem*{theorem*}{Theorem}
\newtheorem*{definition*}{Definition}
\newtheorem*{notation*}{Notation}

\numberwithin{equation}{section}

\title{Unstable Stokes waves \protect\\
in constant vorticity flows of infinite depth}

\begin{document}

\author{
Ting-Yang Hsiao%
\thanks{International School for Advanced Studies (SISSA),
Via Bonomea 265, 34136 Trieste, Italy
\newline
\textit{Email:} \texttt{thsiao@sissa.it}}
,\quad
Vera Mikyoung Hur%
\thanks{Department of Mathematics,
University of Illinois Urbana--Champaign,
Urbana, IL 61801, USA
\newline
\textit{Email:} \texttt{verahur@illinois.edu}}
,\quad
Zhao Yang%
\thanks{State Key Laboratory of Mathematical Sciences, Academy of Mathematics and Systems Science, Chinese Academy of Sciences,
Beijing 100190, China
\newline
\textit{Email:} \texttt{yangzhao@amss.ac.cn}}
}

\date{}
\maketitle

\begin{abstract}
We study the spectral stability of small-amplitude Stokes waves with constant vorticity in infinite depth. We describe the spectrum near the origin and near the first nonzero double eigenvalue. For nonzero vorticity $\omega<2\sqrt{g/3}$, where $g$ is the gravitational constant, the eigenvalues near the origin form a figure eight, whereas for $\omega>2\sqrt{g/3}$ the local spectrum remains purely imaginary. Away from the origin, the interaction of two Fourier modes whose indices differ by two produces an unstable isola when $\omega>-\sqrt{g/6}$, $\omega\ne0$, and the resonant wave numbers are noninteger. This curve is a real analytic oval with an explicit leading ellipse and maximal growth of quadratic order in the amplitude. We also prove instability at the integer resonances needed to cover the remaining vorticities. Together with the classical irrotational result, these two mechanisms show that sufficiently small nonzero Stokes waves are spectrally unstable for every fixed constant vorticity. The proof uses finite-dimensional spectral reductions of the full Euler linearization. Near the origin, 
the finite-dimensional linearized operator has a one-sided limit as the Floquet exponent vanishes, which differs from the periodic linearization. We use the Kato transformation to construct a three-dimensional reduction with biorthogonal bases and use all three eigenvectors of the reference operator to calculate the Benjamin--Feir coefficient. At the nonzero collision, a real two-dimensional reduction and an explicit second-order calculation determine the unstable curve.
\end{abstract}

\tableofcontents

\section{Introduction}\label{sec1}

We consider two-dimensional gravity water waves of infinite depth with constant vorticity. The waves are periodic in the horizontal direction and travel at a constant speed. Our purpose is to determine how the background shear affects the spectral stability of a small-amplitude Stokes wave. We study longitudinal perturbations, which depend on the direction of propagation and preserve the prescribed constant vorticity.

There are two relevant parts of the spectrum. Eigenvalues close to the origin describe the Benjamin--Feir, or modulational, instability. They form a figure eight when the wave is unstable to long-wave perturbations. Other instabilities arise from collisions at nonzero imaginary eigenvalues. The unstable branches produced by these collisions form small closed curves, usually called \emph{isolas}. We describe the low-frequency spectrum and the first such isola. Their parameter ranges overlap, and together they imply spectral instability for every fixed constant vorticity and sufficiently small nonzero amplitude.

\begin{figure}[ht]
\centering
\includegraphics[width=0.92\textwidth]{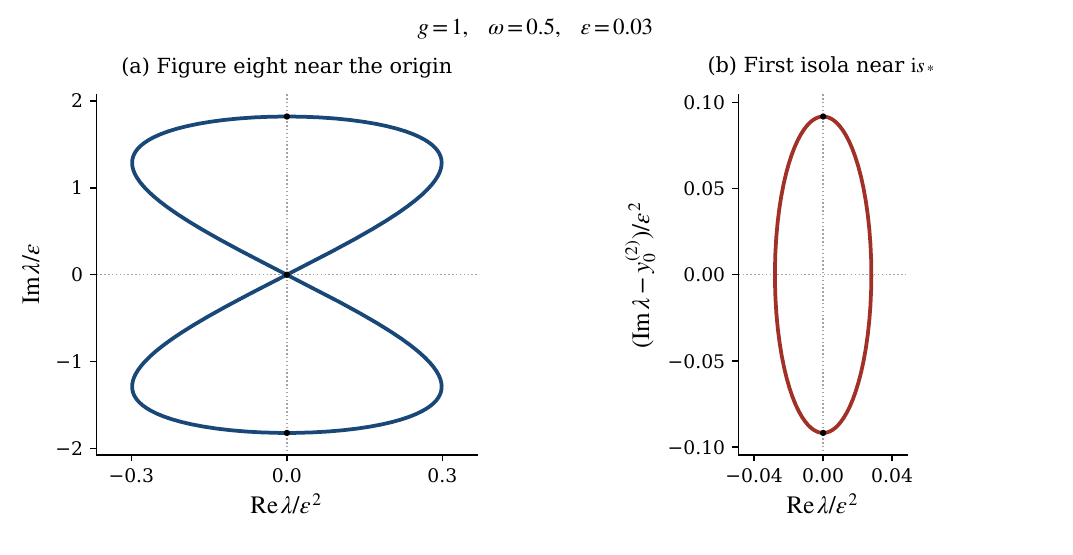}
\caption{Leading spectral curves for $g=1$, $\omega=0.5$, and $\e=0.03$. (a) The Benjamin--Feir figure eight near the origin. (b) The first isola near $\im s_*$. The axis labels show the different scales. Here $y_0^{(2)}$ is the second-order approximation to the center in \eqref{center y0}. The panels illustrate the leading analytic formulas in \eqref{intro:BF-pair}--\eqref{intro:DeltaBF} and \eqref{isola:quadratic-eigenvalues}, with the remainders omitted. The black dots mark the crossing at the origin and the collisions on the imaginary axis.}
\label{fig:figure-eight-first-isola}
\end{figure}

Figure~\ref{fig:figure-eight-first-isola} shows the leading curves for one vorticity in the overlap of the two instability ranges. The different axis scalings display their different sizes.

\paragraph{Related work}
The modulational instability of irrotational Stokes waves was discovered by Benjamin and Feir \cite{BF}. Related modulation criteria were obtained by Lighthill \cite{Li} and Zakharov \cite{Zak68}. The full Euler problem requires more than an envelope approximation: one must locate the eigenvalues of the linearized water-wave equations and control the error terms on the scale where the instability occurs. Bridges and Mielke \cite{BrM} proved the Benjamin--Feir instability in finite depth by spatial dynamics and a center-manifold reduction. Nguyen and Strauss \cite{NS} proved spectral modulational instability in deep water by a Lyapunov--Schmidt reduction and obtained the unstable branches near the origin. Hur and Yang \cite{HY} developed a periodic Evans-function method and obtained instability results both near and away from the origin in finite depth. 
The shape of the unstable spectrum was determined by Berti, Maspero and Ventura \cite{BMV1}. They reduced the four-dimensional spectral space at the origin to two $2\times2$ blocks and proved that the unstable pair traces a complete figure eight. Their argument uses a symplectic Kato transformation and the Hamiltonian and reversible structure of the linearized equations. They subsequently treated finite depth in \cite{BMV3,BMV_ed}. These results describe the local spectral geometry of irrotational gravity waves.

Instability away from the origin has a different source. At zero amplitude, two imaginary eigenvalues belonging to distinct Fourier modes may coincide. Their interaction at nonzero amplitude can produce eigenvalues with nonzero real parts. Deconinck and Oliveras \cite{DO} computed the spectrum of periodic gravity waves and described such instability bands numerically. Creedon, Deconinck and Trichtchenko \cite{CDT} developed a perturbation calculation for the first few high-frequency instabilities and compared the asymptotic curves with numerical spectra. The rigorous analysis in \cite{HY} also detects instability at nonzero resonances. Berti, Maspero and Ventura \cite{BMV_isola} proved the existence of the first unstable isola in irrotational deep water and determined its local geometry. In that problem the second- and third-order interactions vanish, and a fourth-order calculation is needed to open the isola. Berti, Corsi, Maspero and Ventura \cite{BCMV_many} studied higher resonances at finite depth. For each integer $p\ge2$, they obtained an explicit leading coupling coefficient, analytic in the depth and not identically zero. The corresponding isola is present whenever this coefficient is nonzero. Their further analysis in \cite{BCMV_isolas} describes the deep-water asymptotics of the leading coefficients for $p=2,3,4$. These results also explain why the order needed in a finite-depth expansion can change in the deep-water limit.

For constant vorticity, the velocity is the sum of a linear shear and the gradient of a harmonic potential. The equations can therefore still be written in variables on the free surface. Constantin, Ivanov and Prodanov \cite{CIP08} identified a nearly-Hamiltonian structure, and Wahl\'en \cite{Wahlen07} gave a Hamiltonian formulation in canonical variables.
Related existence results include the quasi-periodic waves of Berti, Franzoi and Maspero \cite{BFM1,BFM2} and the Stokes-wave bifurcation results of Barbieri, Berti, Maspero and Mazzucchelli \cite{BBMM26}.

The effect of shear on modulation has also been studied through asymptotic models and numerical calculations. Thomas, Kharif and Manna \cite{TKM12} derived a nonlinear Schr\"odinger equation in finite depth. Their calculation shows that the vorticity contributes to the modulational coefficient. Guyenne, Kairzhan and Sulem \cite{GKS22} derived a Hamiltonian Dysthe equation in deep water and compared its predictions with the full Euler equations. Murashige and Choi \cite{MC20} carried out full-Euler numerical stability calculations for deep-water waves on a linear shear current.

The finite-depth problem with nonzero constant vorticity has also been rigorously investigated. The first author established a modulational instability criterion in his PhD thesis \cite{HsiaoThesis} while Fan and Sun \cite{FanSun2026} obtained the result by a different approach using Kato's similarity transform. Our companion work \cite{HHYfinite} extends the analysis of \cite{HsiaoThesis}, improving the modulational instability result and establishing spectral instability away from the origin. These works provide rigorous finite-depth counterparts of the present analysis.

The present paper gives a direct spectral analysis of the deep-water Euler equations with constant vorticity. Near the origin, we determine the three eigenvalue branches of $\mathcal L_{\mu,\e}$, $\mu\ne0$, including the contribution from the eigenvector $f_0$ in \eqref{def:right-reference-basis}. Away from the origin, we compute the coefficients of $\e^2$ at the $p=2$ resonance and prove that they open an isola in the stated parameter range. In particular, a shear which removes the Benjamin--Feir instability does not remove spectral instability: the nonzero resonance supplies unstable eigenvalues in that range.

\paragraph{Main results}
Fix the spatial period to be $2\pi$, and let $g>0$ be the gravitational acceleration. We use the vorticity convention of Section~\ref{sec 2} and write
\begin{equation}\label{intro:speed}
\cw=\frac{\omega+\sqrt{\omega^2+4g}}2,
\qquad \cw^2-\omega\cw-g=0.
\end{equation}
The parameter $\e$ is the Stokes amplitude used in Theorem~\ref{Thm: Stokes expansion}. In a frame moving with the wave, a Bloch mode has the form
\[
e^{\lambda t}e^{\im\mu x}v(x),\qquad v(x+2\pi)=v(x).
\]
Here $\lambda$ is the temporal eigenvalue and $\mu$ is the Bloch exponent, represented in $[-\tfrac12,\tfrac12)$ and understood modulo one. We study the spectrum of the real-line linearization $\mathcal A_\e$ on $L^2(\mathbb R;\mathbb C^2)$, defined in \eqref{eq:flattened-real-line-operator}. Spectral instability means
\[
\sigma_{L^2(\mathbb R;\mathbb C^2)}(\mathcal A_\e)
\cap\{\lambda\in\mathbb C:\operatorname{Re}\lambda>0\}\ne\varnothing.
\]
The continuous unstable Bloch branches constructed below imply this property by Lemma~\ref{lem:Bloch-to-real-line-spectrum}. All smallness thresholds below may depend on the fixed pair $(g,\omega)$.

We first state the complete spectral description near the origin. The Benjamin--Feir coefficient is
\begin{equation}\label{intro:ThetaBF}
\mathrm{ind}_{\mathrm{BF}}
:=
\frac{g(3g-\cw^2)(\cw^4+g^2)}{(\cw^2+g)^3}.
\end{equation}
It is the coefficient of $\e^2$ in the discriminant below.

\begin{theorem}[The figure-eight spectrum]
\label{thm:main1}
Let $\omega\ne0$ and $\omega\ne2\sqrt{g/3}$. There exist $\e_0,\mu_0,r_*>0$ such that, for $0<\e<\e_0$ and $0<|\mu|<\mu_0$, the Bloch operator has exactly three eigenvalues in $\{|\lambda|<r_*\}$, counted with algebraic multiplicity. In terms of the analytic blocks $A^{(3)}$ and $d^{(3)}$ of the exactly triangularized matrix in \eqref{def:M3-section5}, the purely imaginary branch is
\begin{equation}\label{intro:lambda0}
\lambda_0(\mu,\e)
=
\im\frac{\cw^2}{\omega}\mu
+\im\mu\cO(|\mu|+|\e|)
\in\im\mathbb R.
\end{equation}
and the remaining pair is
\begin{equation}\label{intro:BF-pair}
\lambda_\pm(\mu,\e)
=
\im\cw\mu
-\frac{\im\cw g}{\cw^2+g}\mu
\pm\frac{|\mu|}{2}\sqrt{\Delta_{\mathrm{BF}}(|\mu|,\e)}+\im \mu\cO\big((|\mu|+|\e|)^2\big).
\end{equation}
The blocks $A^{(3)}$ and $d^{(3)}$ are purely imaginary for real parameters, and
\begin{align}
\Delta_{\mathrm{BF}}(\mu,\e):=
\mathrm{ind}_{\mathrm{BF}}\e^2
-\frac{4\cw^6g^4}{(\cw^2+g)^6}\mu^2
+\cO\big((|\mu|+|\e|)^3\big).
\label{intro:DeltaBF}
\end{align}
The discriminant is real analytic. For a negative discriminant, the square root in \eqref{intro:BF-pair} is understood as $\im\sqrt{-\Delta_{\mathrm{BF}}}$. For negative $\mu$, the formulas describe the complex-conjugate eigenvalue set at $|\mu|$.

If $\omega<2\sqrt{g/3}$, there is a unique positive collision point in $(0,\mu_0)$,
\begin{equation}\label{intro:BF-endpoint}
\mu_*(\e)
=
\frac{(\cw^2+g)^3}{2\cw^3g^2}
\sqrt{\mathrm{ind}_{\mathrm{BF}}}\,\e
+\cO(\e^2).
\end{equation}
The pair has opposite nonzero real parts for $0<|\mu|<\mu_*(\e)$, collides on the imaginary axis at $|\mu|=\mu_*(\e)$, and is purely imaginary for $\mu_*(\e)<|\mu|<\mu_0$. After defining $\lambda_\pm(0,\e)=0$ by continuity, the two branches for $|\mu|\le\mu_*(\e)$ form a figure eight.

If $\omega>2\sqrt{g/3}$, all three eigenvalues in this neighborhood of the origin are purely imaginary.
\end{theorem}

The exact eigenvalue formulas and the collision analysis are proved in Theorems~\ref{thm:exact-eigenvalues-section5}.

We next describe the first nonzero resonance. Suppose $\omega>-\sqrt{g/6}$. Define
\begin{equation}\label{intro:isola-anchor}
\xi_*:=\frac{3\cw^2-2g}{4g},\qquad
s_*:=\frac{3\cw^3}{4g},\qquad
\bar\mu:=\xi_*-\left\lfloor\xi_*+\frac12\right\rfloor
\in\left[-\frac12,\frac12\right).
\end{equation}
At zero amplitude, the two positive wave numbers $\xi_*$ and $\xi_*+2$ have the same temporal eigenvalue $\im s_*$. This collision exists for $\omega>-\sqrt{g/6}$. The positive coefficient governing the size of the resulting isola is
\begin{equation}\label{intro:isola-size}
\sqrt{\mathrm{ind}_{\mathrm{high}}}=\frac{3\sqrt3\,|\omega|^3(\cw^2+g)}{64g^2}
\sqrt{\frac{(\cw^2+2g)(3\cw^2-2g)}
{(2\cw^2-g)(2\cw^2+g)}}.
\end{equation}

Local intervals of Bloch exponents are understood on $\mathbb R/\mathbb Z$, using a real coordinate near $\bar\mu$. In particular, the endpoints $-\tfrac12$ and $\tfrac12$ are identified.

\begin{theorem}[The first elliptical isola]\label{thm:main2}
Suppose $\omega>-\sqrt{g/6}$, $\omega\ne0$, and $\xi_*\notin\mathbb Z$. For every sufficiently small $\e>0$, the two eigenvalues near $\im s_*$ have opposite nonzero real parts precisely on a local interval $\mu_-(\e)<\mu<\mu_+(\e)$. On the closed interval $\mu_-(\e)\leq\mu\leq\mu_+(\e)$, the two eigenvalues satisfy
\begin{equation}
\begin{gathered}
\lambda^{(2)}_\pm(\mu,\e)
=\im\left(s_*+\frac{4\cw^5}{4\cw^4-g^2}(\mu-\bar\mu)
+\frac{C_{11}+C_{22}}2\e^2+\cO(|\e|^3)\right)\\
\qquad\pm\sqrt{
\mathrm{ind}_{\mathrm{high}}\e^4-
\left(\dfrac{2g\cw^3}{4\cw^4-g^2}(\mu-\bar\mu)+\frac{C_{11}-C_{22}}{2}\e^2\right)^2
+\cO(|\e|^5)},
\end{gathered}
\end{equation}
where $C_{11}$ and $C_{22}$ are defined in \eqref{isola:C11} and \eqref{isola:C22} respectively. Complex conjugation gives a second curve near $-\im s_*$.
\end{theorem}

Let us make some comments:

\begin{enumerate}
\item {\sc Benjamin--Feir figure eight.} For $\omega\ne0$ and $\mathrm{ind}_{\mathrm{BF}}>0$ (cf. $\omega<2\sqrt{g/3}$) for values of the Floquet parameter $ 0<\mu <  \mu_*(\e) $, the eigenvalues $\lambda_\pm (\mu, \epsilon) $ have opposite non-zero real part, according to \eqref{intro:BF-pair}. As $ \mu $ tends to $  \mu_*(\e)$, the two eigenvalues $\lambda_\pm (\mu, \epsilon)$ collide on the imaginary axis {\it far} from $ 0 $, along which they  keep moving  for $ \mu > \mu_*(\e) $. In particular for $ \mu \in [0, \mu_*(\e)]$ we obtain the upper part of the figure  ``8'', see Figure \ref{fig:figure-eight-first-isola},  which is well approximated by the curves
\begin{equation}\label{intro:figure-eight-parametrization}
\begin{gathered}
\mu\longmapsto
\left(
\pm\frac{|\mu|}{2}
\sqrt{\mathrm{ind}_{\mathrm{BF}}\e^2
-\frac{4\cw^6g^4}{(\cw^2+g)^6}\mu^2},
\ \frac{\cw^3}{\cw^2+g}\mu
\right).
\end{gathered}
\end{equation}
The upper and lower lobes correspond to $\mu>0$ and $\mu<0$, respectively.

\item {\sc First elliptical isola.} For $\omega\neq 0$, $\omega>-\sqrt{g/6}$, and $\xi_*\notin\mathbb Z$ (cf. \eqref{isola:theorem-assumptions}) for values of the Floquet parameter $\mu_-(\e)<\mu<\mu_+(\e)$, the eigenvalues $\lambda^{(2)}_\pm (\mu, \epsilon) $ have opposite non-zero real part, according to \eqref{isola:quadratic-eigenvalues}. As $ \mu $ tends to $  \mu_\pm(\e)$, the two eigenvalues $\lambda^{(2)}_\pm (\mu, \epsilon)=x_{\pm}(\mu,\e)+\im y(\mu,\e)$ collide on the imaginary axis, along which they  keep moving  for $ \mu > \mu_+(\e)$ and $\mu<\mu_-(\e)$. In particular, for $\mu_-(\e)<\mu<\mu_+(\e)$ we obtain the elliptical curve, see Figure \ref{fig:figure-eight-first-isola},  which is well approximated by the equation
\begin{equation}
x^2_{\pm}+\frac{g^2}{4\cw^4}\left(y-y^{(2)}_0\right)^2=\mathrm{ind}_{\mathrm{high}}\e^4+\cO(|\e|^5).
\end{equation}
where
\begin{align} \label{center y0}
y_0^{(2)}(\e)
:=s_*+
\left(
\frac{C_{11}+C_{22}}2
-\frac{\cw^2}{g}(C_{11}-C_{22})
\right)\e^2.
\end{align}
Complex conjugation gives the corresponding isola near $-\im s_*$.
\end{enumerate}

\begin{theorem}[Instability for every constant vorticity]
\label{thm:main3}
For every $g>0$ and every fixed $\omega\in\mathbb R$, there exists $\e_0=\e_0(g,\omega)>0$ such that each Stokes wave of amplitude $0<\e<\e_0$ is spectrally unstable under longitudinal perturbations in $L^2(\mathbb R;\mathbb C^2)$. More precisely, $\mathcal L_{\mu,\e}$ has an eigenvalue with positive real part for some $\mu\ne0$.
\end{theorem}

Indeed, if $\omega<2\sqrt{g/3}$ and $\omega\ne0$, the figure-eight theorem supplies an unstable eigenvalue. At $\omega=0$ the classical result \cite{NS,HY,BMV1} supplies one. If $\omega\ge2\sqrt{g/3}$, the first isola or Proposition~\ref{prop:isola-integer} gives instability. At the boundary $\omega=2\sqrt{g/3}$, one has $\xi_*=7/4$ and $b=\sqrt g/2$, so Theorem~\ref{thm:main2} applies directly. This proves the theorem; see also Theorem~\ref{thm:all-vorticities-unstable}. The assertion is for each fixed vorticity, with an amplitude bound that may depend on it.

\begin{remark*}\rm
While finalizing this manuscript, we became aware of the recent arXiv posting by Chen and Nguyen \cite{CN2026}, which independently establishes the Benjamin--Feir instability. Our Theorem \ref{thm:main1}, obtained independently of \cite{CN2026}, establishes the same instability criterion. Although the two proofs share some underlying ideas, they were developed independently and differ in their technical implementation.

Our work goes further by investigating the spectral instability away from the origin. In Theorem \ref{thm:main2}, we establish the existence and determine to the leading order in the amplitude the first unstable isola. Particularly, although sufficiently large positive constant vorticity suppresses the Benjamin--Feir instability, it does not suppress spectral instability altogether. Combining the two instability mechanisms in Theorems \ref{thm:main1} and \ref{thm:main2}, we prove in Theorem \ref{thm:main3} that every sufficiently small amplitude Stokes waves in infinite depth is spectrally unstable for any constant vorticity.
\end{remark*}

\paragraph{Outline of the proof}
The spectral analysis near the origin is based on the one-sided limits of $\mathcal L_{\mu,\e}$ as $\mu\to0$ with $\mu\ne0$. In the presence of constant vorticity, the term $\omega\partial_x^{-1}G(\eta)\psi$ distinguishes the zero Fourier mode. After the good-unknown and conformal changes of variables \cite{Alinhac,LDbook}, the real-line linearization contains the Hilbert transform. Its Bloch multiplier $-\im\operatorname{sign}(k+\mu)$ has distinct one-sided limits at $k=0$. The positive Bloch family therefore extends to
\[
\mathcal L_{0^+,\e}:=\lim_{\mu\downarrow0}\mathcal L_{\mu,\e},
\]
which differs from the auxiliary periodic linearization $\mathcal L^{\mathrm{per}}_{\e}$ in \eqref{auxiliary periodic operator}. Their difference is the finite-rank operator in \eqref{finite rank jump}, already present at zero amplitude. For positive Bloch exponents, the reference operator is therefore $\mathscr L_0=\mathcal L_{0^+,0}$ in \eqref{def:scriptL-star}.

For fixed $\omega\ne0$, zero is an isolated semisimple eigenvalue of $\mathscr L_0$ with multiplicity three, with eigenvectors $f_+,f_-,f_0$ in \eqref{def:right-reference-basis}. The other eigenvalue of the constant Fourier block is $-\im\omega$; at zero vorticity it joins the spectral cluster at the origin, whose algebraic dimension becomes four \cite{BMV1}. We transport the three right and left reference vectors using $U_{\mu,\e}f_j$ and $(U_{\mu,\e}^{-1})^*g_j$, where $U_{\mu,\e}$ is the Kato transformation in \eqref{def:Kato-U} \cite{Kato1966}. This gives an analytic matrix $\mathbf M(\mu,\e)$ representing the shifted operator $\mathscr L_{\mu,\e}$. Its eigenvalues give those of $\mathcal L_{\mu,\e}$ near the origin after adding $\im\cw\mu$, as in \eqref{reduced-spectrum-original-L}. The reduction is carried out for fixed nonzero vorticity; the comparison with the irrotational coefficients is given below. 

At nonzero amplitude, the generalized kernel of $\mathscr L_{0,\e}=\mathcal L_{0^+,\e}$ is obtained from translation, differentiation along the Stokes branch, and vertical translation. The corrections due to \eqref{finite rank jump} give the exact relations in Lemma~\ref{lem:one-sided-Jordan-section5}. For $\e\ne0$, zero has algebraic multiplicity three and geometric multiplicity two, and the restriction to its generalized kernel is nilpotent of index two. In particular, $\mathbf M(0,\e)^2=0$. Conservation of the mean of the first component also gives a zero last row at $\mu=0$. The analytic change of coordinates $R_\e$ then puts $\mathbf M^{(1)}(0,\e)$ in the form \eqref{eq:M1-at-zero-section5}. The parity of the corrected translation vector and the adjoint kernel vectors gives $\partial_\mu\mathbf M^{(1)}_{12}(0,\e) =\partial_\mu\mathbf M^{(1)}_{32}(0,\e)=0$, as proved in Lemma~\ref{lem:M1-properties-section5}.

These identities make the singular scaling
\[
Q(\mu)=\operatorname{diag}(\mu,1,\mu),
\qquad
\mathbf M^{(2)}(\mu,\e)
=\frac1\mu Q(\mu)^{-1}\mathbf M^{(1)}(\mu,\e)Q(\mu)
\]
analytic through $\mu=0$. Since $\sigma(\mathbf M)=\mu\sigma(\mathbf M^{(2)})$ for $\mu\ne0$, all three eigenvalues near the origin are $\cO(\mu)$. At $(\mu,\e)=(0,0)$, the scalar block has eigenvalue $\im g/\omega$, separated from the double eigenvalue $-\im\cw g/(\cw^2+g)$ of the upper-left block. After this scaling, the analytic implicit function theorem gives a real analytic change of coordinates. More precisely, a real analytic row vector $X$ defines the transformation
\[
T=\begin{pmatrix}I_2&0\\X&1\end{pmatrix}
\]
and gives the exact block upper triangular matrix
\[
\mathbf M^{(3)}
=T^{-1}\mathbf M^{(2)}T
=\begin{pmatrix}
A^{(2)}+B^{(2)}X&B^{(2)}\\
0&d^{(2)}-XB^{(2)}
\end{pmatrix},
\]
as stated in Lemma~\ref{lem:exact-block-triangularization}.

The blocks $A^{(3)}=A^{(2)}+B^{(2)}X$ and $d^{(3)}=d^{(2)}-XB^{(2)}$ are purely imaginary for real parameters. Thus the scalar block gives a purely imaginary eigenvalue, while the remaining pair is determined directly by $\operatorname{tr}A^{(3)}$ and $\det A^{(3)}$. The first-order term of $x_1$ determines the contribution of $B^{(2)}X$ to the quadratic term of $A^{(3)}_{21}$. Consequently, the Benjamin--Feir coefficient is
\[
\mathrm{ind}_{\mathrm{BF}}
=
\frac{8\cw^3g^2}{(\cw^2+g)^3}\Theta
-\frac{\omega^2g}{\cw^2+g},
\]
in agreement with \eqref{intro:ThetaBF}. Only the bound $x_2=\cO(|\mu|+|\e|)$ is needed here. The upper-right entry of $A^{(3)}$ has a nonzero constant term, the lower-left entry starts at quadratic order, and the diagonal difference is $\cO((|\mu|+|\e|)^2)$. Therefore the identity
\[
\Delta_{\mathrm{BF}}
=
\big(A^{(3)}_{11}-A^{(3)}_{22}\big)^2
+4A^{(3)}_{12}A^{(3)}_{21}
\]
shows that the quadratic expansion of $\mathbf M$ suffices to obtain
\[
\Delta_{\mathrm{BF}}(\mu,\e)
=
\mathrm{ind}_{\mathrm{BF}}\e^2
-\frac{4\cw^6g^4}{(\cw^2+g)^6}\mu^2
+\cO\big((|\mu|+|\e|)^3\big).
\]
Multiplication by $\mu$ and addition of $\im\cw\mu$ recover the exact eigenvalue formulas of Theorem~\ref{thm:exact-eigenvalues-section5}.

The real analytic discriminant describes both the splitting for $0<\mu\ll\e$ and the return to the imaginary axis at $\mu\sim\e$. When $\mathrm{ind}_{\mathrm{BF}}>0$, the analytic implicit function theorem locates the positive collision point $\mu_*(\e)$ in \eqref{intro:BF-endpoint}; strict concavity in $\mu$ gives its uniqueness and the sign classification on either side. Monotonicity of the common imaginary part gives the simple upper lobe of the figure eight. The lower lobe follows from complex conjugation of the physical Bloch operators, as in Lemma~\ref{lem:negative-bloch-symmetry}. The analytic continuation from $\mu>0$ to negative $\mu$ is used only within the reduction. When $\mathrm{ind}_{\mathrm{BF}}<0$, the discriminant is negative near the origin and all three eigenvalues are purely imaginary. The critical case $\mathrm{ind}_{\mathrm{BF}}=0$ requires higher-order coefficients and is excluded from the present description.

Away from the origin, the first nonzero collision at $\bar\mu\ne0$ is treated by a two-dimensional Kato reduction, following \cite{BMV_isola}. The first amplitude derivative shifts Fourier indices by one, whereas the two resonant indices differ by two. Consequently, the term linear in $\e$ in the reduced matrix vanishes. The coefficients of $\e^2$ contain both the direct second-order contribution and the contribution through neighboring Fourier modes. Their explicit evaluation gives a strictly negative product of the off-diagonal coefficients for nonzero vorticity in the stated parameter range. The isola therefore opens at order $\e^2$. At zero vorticity, these quadratic off-diagonal coefficients vanish, consistently with the fourth-order opening obtained in \cite{BMV_isola}.

The implicit-function theorem determines the center and the two endpoints of the unstable interval. The quadratic formula then gives the eigenvalue expansion and the sign of the discriminant. The leading terms describe an ellipse, and the maximal growth rate is $\sqrt{\mathrm{ind}_{\mathrm{high}}}\e^2+\cO(|\e|^3)$.

For integer resonant wave numbers $\xi_*=m$, the collision occurs at $\bar\mu=0$. We apply the same reduction to the analytic extension from $\mu>0$, with reference operator $\mathcal L_{0^+,0}$. In the range $\cw^2\ge3g$, its zero Fourier block has eigenvalues $0$ and $-\im\omega$, both distinct from $\im s_*$, so $\im s_*$ remains an isolated double eigenvalue. Moreover, $m\ge2$, and the intermediate modes $m-1,m+1,m+3$ entering the quadratic coefficients are all positive. The previously computed coefficients therefore remain valid. Their diagonal difference determines the center $\mu_0(\e)=\mu_2\e^2+\cO(\e^3)$, and an explicit calculation gives $\mu_2>0$. Thus $0<\mu_0(\e)<\frac12$ for sufficiently small $\e>0$, and the quadratic formula yields an unstable eigenvalue of the physical Bloch operator, as stated in Proposition~\ref{prop:isola-integer}. Together with the noninteger isola, the low-frequency result, and the classical irrotational theorem \cite{NS,BMV1}, this proves spectral instability for every fixed constant vorticity, including the Benjamin--Feir threshold.

\paragraph{Comparison with zero vorticity}

At $\omega=0$, the water-wave equations reduce to the irrotational problem. The explicit leading coefficients in \eqref{intro:BF-pair}--\eqref{intro:DeltaBF} also recover the classical irrotational formulas. Indeed, evaluating \eqref{intro:speed} and \eqref{intro:ThetaBF} at $\omega=0$ gives
\begin{equation}\label{intro:zero-coefficients}
\cw=\sqrt g,\qquad
\mathrm{ind}_{\mathrm{BF}}=\frac g2.
\end{equation}
Therefore, the corresponding leading eigenvalue pair is
\begin{equation}\label{intro:zero-leading-pair}
\lambda^{\mathrm{lead}}_\pm(\mu,\e)
=\frac{\im\sqrt g}{2}\mu
\pm\frac{\sqrt g\,|\mu|}{8}
\sqrt{8\e^2-\mu^2}.
\end{equation}
With period $2\pi$ and leading surface profile $\eta_\e(x)=\e\cos (x)+\mathcal O(\e^2)$, this agrees with the leading figure-eight formula of Berti--Maspero--Ventura \cite[Theorem 1.1]{BMV1}, after restoring the factor $\sqrt g$ in the temporal eigenvalue. For this leading curve, the positive collision occurs at $\mu=2\sqrt2\,\e$, and the maximal real part is $\sqrt g\,\e^2/2$, attained at $|\mu|=2\e$.

Nguyen--Strauss \cite[Theorem 1.1]{NS} prove that, for each sufficiently small nonzero amplitude, instability occurs for all $0<|\mu|<\mu_0(\e)$, where the threshold may depend on the amplitude. In the same normalization, their unstable eigenvalue satisfies
\begin{equation}\label{intro:NS-comparison}
\lambda
=\frac{\im\sqrt g}{2}\mu
+\frac{\sqrt g}{2\sqrt2}|\mu\e|
+\mathcal O_\e(\mu^2)+\mathcal O(|\mu|\e^2).
\end{equation}
Here the subscript in $\mathcal O_\e(\mu^2)$ allows its constant to depend on the fixed amplitude. Expanding \eqref{intro:zero-leading-pair} for $|\mu|\ll\e$ gives the same leading drift and growth terms. The Nguyen--Strauss expansion describes the crossing near the origin, whereas the Berti--Maspero--Ventura formula follows the branches through their subsequent collision on the imaginary axis. Theorem~\ref{thm:main1} provides this complete local description for fixed nonzero vorticity in the Benjamin--Feir unstable regime.

This agreement concerns the explicit leading coefficients; it does not assert a uniform $\omega\to0$ limit of the three-dimensional reduction in Sections~\ref{sec:Kato-positive}--\ref{sec:characteristic-reduction-positive}. For fixed $\omega\ne0$, the constant Fourier block of the one-sided reference operator $\mathscr L_0=\mathcal L_{0^+,0}$ is
\[
\begin{bmatrix}
0&0\\
-g&-\im\omega
\end{bmatrix},
\]
with eigenvalues $0$ and $-\im\omega$. At $\omega=0$, this block becomes nilpotent, and the algebraic multiplicity of zero for the reference operator increases from three to four. Thus the spectral separation required for the three-dimensional Kato reduction is lost. The leading slope $\cw^2/\omega$ of the branch $\lambda_0$ in \eqref{intro:lambda0} also exhibits this singularity. Our reduction and the subsequent block triangularization are therefore carried out for fixed nonzero vorticity, with smallness thresholds and remainder bounds depending on $\omega$. At zero vorticity, the exact spectral conclusions follow from the four-dimensional irrotational analyses of \cite{NS,BMV1}; the calculation above verifies consistency of our leading Benjamin--Feir coefficients with those results.

There is a similar comparison for the first isola. At zero vorticity, \eqref{intro:isola-anchor} gives $\xi_*=1/4$ and $s_*=3\sqrt g/4$, the collision of \cite{BMV_isola}. Moreover,
\begin{equation}\label{intro:zero-isola-size}
\sqrt{\mathrm{ind}_{\mathrm{high}}}=\frac{3\sqrt3}{32g}|\omega|^3
\left(1+\mathcal O\!\left(\frac{|\omega|}{\sqrt g}\right)\right)
\qquad (\omega\to0).
\end{equation}
Thus the nonzero-vorticity interaction of order $\e^2$ vanishes at $\omega=0$. Berti--Maspero--Ventura \cite{BMV_isola} prove that the irrotational isola opens at order $\e^4$. This is consistent with our calculation; our fixed-vorticity estimate does not claim an amplitude bound uniform as $\omega\to0$. The second-order center shifts also agree with that work, as recorded in Remark~\ref{rem:isola-order}.

\paragraph{Organization}
Section~\ref{sec 2} introduces the equations and Stokes waves, derives the real-line linearization and its Bloch operators, and states the complete low-frequency theorem. Section~\ref{sec:Kato-positive} uses the Kato transformation to construct a three-dimensional reduction. Section~\ref{sec:matrix-representation-positive} computes its quadratic expansion, and Section~\ref{sec:characteristic-reduction-positive} proves the Benjamin--Feir spectrum and the figure-eight description. Section~\ref{sec:first-isola} treats the $p=2$ resonance, the isola, and the integer resonances, and combines the results for all vorticities. The appendices give the Stokes and conformal expansions and the detailed right and left basis calculations.

\section{The Benjamin--Feir instability and the full description of the spectrum}\label{sec 2}
In this section we present the complete spectral Theorem \ref{Complete BF thm}.
We first introduce the constant-vorticity water-wave equations and the Stokes waves.     \\
{\bf The water waves equations.} We consider the Euler equations of hydrodynamics for a 2-dimensional perfect and incompressible fluid  under the action of gravity at the free surface. The fluid fills the region
    \begin{align*}
       \mathcal{D}_\eta:=\{(x,y)\in \mathbb{R}\times \mathbb{R}: -\infty < y <\eta(t,x)\}\ ,
    \end{align*}
in deep water.

Imposing that the fluid particles at the free surface remain on it along the evolution (kinematic boundary condition) and that the pressure at the free surface is equal to the constant atmospheric pressure (dynamic boundary condition), we use the constant-vorticity surface formulation associated with the Zakharov--Craig--Sulem variables; see \cite{Zak68,CS} for the irrotational formulation and \cite{CIP08,Wahlen07} for its constant-vorticity Hamiltonian extension.  In the moving frame the equations read
\begin{equation} \label{Craig-Sulem formula}
\left\{\begin{aligned}
    \eta_t&=c\eta_x+G(\eta)\psi+\omega\eta\eta_x,\\
    \psi_t&=c\psi_x-g\eta-\frac{\psi_x^2}{2}+\frac{1}{2(1+\eta^2_x)}\left(G(\eta)\psi+\eta_x\psi_x\right)^2+\omega \eta\psi_x+\omega\pa_x^{-1} G(\eta)\psi,
\end{aligned}\right.
\end{equation}
where $g>0$ is the gravitational constant and $G(\eta)$ is the  Dirichlet-Neumann operator
$$
G(\eta)\psi := (- \Phi_x \eta_x + \Phi_y)\vert_{y = \eta(x)} \, . 
$$
For the periodic background problem in \eqref{Craig-Sulem formula}, $ (\pa_x^{-1} f) (x) $ denotes the unique primitive with  zero average  of a $ 2 \pi  $-periodic zero average function $ f(x) $. It turns out that $ G(\eta)\psi  $ has zero average. As consequence the average of $ \eta (x) $, $ \langle \eta \rangle := \eta_0 := \frac{1}{2  \pi } \int_\mathbb{ T} \eta (x) \, d x $ is conserved along solutions of \eqref{Craig-Sulem formula}. Note also that, since $ G(\eta) [1] = 0 $ vanishes on the constants, the vector field in the right hand side of \eqref{Craig-Sulem formula} does not depend on $ \psi_0 = \la \psi \ra $. For real-line perturbations, the variation of the primitive is instead normalized in $L^2(\mathbb R)$, as specified in \eqref{linearized primitive} below.

The local bifurcation of small periodic Stokes waves with constant vorticity is well established; see, in particular, \cite{BBMM26}.  We record here the analytic normalization and the Taylor coefficients required in the spectral calculation.

 \begin{theorem} [Stokes waves] \label{Thm: Stokes expansion} 
Let $g>0$ and $\omega \in \R $. There exist $\e_*:=\e_*(g,\omega) >0$ and a unique family  of real analytic solutions $(\eta_\e(x), \psi_\e(x), c_\e)$, parameterized by the amplitude $|\e| \leq \e_*$, of
\begin{equation}\label{travelingWWstokes}
c \, \eta_x+G(\eta)\psi+\omega\eta\eta_x = 0 \, , \quad 
c \, \psi_x -  g\eta - \dfrac{\psi_x^2}{2} + 
\dfrac{1}{2(1+\eta_x^2)} \big( G(\eta) \psi + \eta_x \psi_x \big)^2  +\omega\eta\psi_x+\omega\pa^{-1}_xG(\eta)\psi= 0 \, , 
\end{equation}
such that $ \eta_\e (x), \psi_\e (x) $ are $2\pi$-periodic;  $\eta_\e (x) $ is even and $\psi_\e (x) $ is odd.  We choose the analytic amplitude parameter by the normalization
\begin{equation}\label{amplitude-normalization}
\frac1\pi\int_0^{2\pi}\eta_\e(x)\cos(x)\,dx=\e.
\end{equation}
With this normalization they are of the form
 \begin{equation}\label{exp:Sto}
 \begin{aligned}
  & \eta_\e (x) =  \e\cos(x)+\e^2\left(\eta_2^{[0]}+\eta_2^{[2]}\cos(2x)\right) +\cO(\e^3), \\ 
  & \psi_\e (x)  = \e\cw\sin(x)+\e^2 \psi_2^{[2]}\sin(2x)+\cO(\e^3) ,\\
  & c_\e = \cw+\e^2 c_2+\cO(\e^3)\quad \text{where} \quad \cw :=\cw(\omega,g):=  \frac{\omega+\sqrt{\omega^2+4g}}{2}\,, 
   \end{aligned}
  \end{equation}
and
  \begin{align}\label{expcoef}
 &\eta_{2}^{[0]} := \frac{\omega\cw}{2g}, \qquad 
\eta_{2}^{[2]} := \frac{2g+\omega^2-2\omega\cw}{4\cw^2}, \qquad \psi_{2}^{[2]} :=  \frac{2g+\omega^2-\omega\cw}{4\cw}, \qquad 
\\  \label{expc2}
&c_2 := \frac{
-4\cw^8+11g\cw^6-7g^2\cw^4+3g^3\cw^2+g^4
}{
8g\cw^5
}.
\end{align}
More precisely for any  $ \sigma \geq  0 $ and $ {s > \frac52} $, there exists $ \e_*>0 $ such that the map $\e \mapsto (\eta_\e, \psi_\e, c_\e)$ is analytic from $B(\e_*) \to H^{\sigma,s}_{\mathtt{ev}} (\T)\times H^{\sigma,s}_{\mathtt{odd}}(\T)\times \R$, where $ H^{\sigma,s}_{\mathtt{ev}}(\T) $, respectively $ H^{\sigma,s}_{\mathtt{odd}}(\T) $, denote the  space of even, respectively odd, real valued $ 2 \pi $-periodic analytic functions $ u(x) = \sum_{k \in \mathbb{Z}} u_k e^{\im k x} $ such that $ \| u \|_{\sigma,s}^2 := \sum_{k \in \mathbb{Z}} |u_k|^2 \langle k \rangle^{2s} e^{2 \sigma |k|} < + \infty$.
\end{theorem}

The expansions in \eqref{exp:Sto}--\eqref{expc2} are proved in Appendix \ref{sec:App2}.

\paragraph{Linearization} 
We keep the $2\pi$-periodic Stokes wave $(\eta_\e,\psi_\e,c_\e)$ as the background and consider perturbations $(\hat\eta,\hat\psi)$ on $\mathbb R$. All spaces below are complexified. The real-line phase space and operator domain are
\begin{equation}\label{eq:real-line-spaces}
\mathcal X_{\mathbb R}:=L^2(\mathbb R;\mathbb C^2),\qquad
\mathcal Y_{\mathbb R}:=H^1(\mathbb R;\mathbb C^2).
\end{equation}
Only the perturbations belong to these spaces; the periodic background need not be square integrable. Let $G_{\e,\mathbb R}$ denote the Dirichlet--Neumann operator on the periodic fluid domain, acting on real-line boundary data. It maps $H^1(\mathbb R)$ boundedly into $L^2(\mathbb R)$.

We first linearize for smooth localized perturbations; the resulting operator extends boundedly from $\mathcal Y_{\mathbb R}$ to $\mathcal X_{\mathbb R}$. Using the shape derivative formula \cite{LD,LDbook},
\begin{equation}\label{shape derivative}
\mathrm d_\eta G(\eta_\e)[\hat\eta]\psi_\e
=-G_{\e,\mathbb R}(B\hat\eta)-\partial_x(V\hat\eta),
\end{equation}
where the real-valued periodic coefficients are
\begin{align}
V&:=V(x):=(\psi_\e)_x-B(\eta_\e)_x, \label{espV}\\
B&:=B(x):=
\frac{G(\eta_\e)\psi_\e+(\psi_\e)_x(\eta_\e)_x}
{1+(\eta_\e)_x^2}
=\frac{(\psi_\e)_x-c_\e-\omega\eta_\e}
{1+(\eta_\e)_x^2}(\eta_\e)_x
:=-\alpha_\e(\eta_\e)_x, \label{espB}
\end{align}
and
\begin{equation}\label{def:alpha}
\alpha_\e:=c_\e-V+\omega\eta_\e
=\frac{c_\e+\omega\eta_\e-(\psi_\e)_x}{1+(\eta_\e)_x^2}.
\end{equation}
Here $G(\eta_\e)\psi_\e$ is evaluated on the periodic background. The kinematic equation linearizes as
\begin{equation}\label{linearized kinematic}
\hat\eta_t
=-G_{\e,\mathbb R}(B\hat\eta)
+G_{\e,\mathbb R}\hat\psi
+\partial_x(\alpha_\e\hat\eta).
\end{equation}

For the nonlocal vorticity term, we use the bounded operator $K_\e:L^2(\mathbb R)\to L^2(\mathbb R)$ characterized by
\begin{equation}\label{eq:real-line-primitive-factorization}
G_{\e,\mathbb R}f=\partial_x(K_\e f),\qquad f\in H^1(\mathbb R).
\end{equation}
Its construction is given in \eqref{eq:real-line-primitive} below. We do not regard $\partial_x^{-1}$ as a bounded operator on arbitrary $L^2(\mathbb R)$ data. Instead, the variation of the primitive is fixed by requiring it to belong to $L^2(\mathbb R)$. Indeed,
\begin{equation}\label{linearized primitive first form}
\delta\big(G(\eta)\psi\big)
=\partial_x\big(K_\e(\hat\psi-B\hat\eta)-V\hat\eta\big).
\end{equation}
Thus the unique square-integrable primitive is
\begin{equation}\label{linearized primitive}
\delta\big(\partial_x^{-1}G(\eta)\psi\big)
=K_\e(\hat\psi-B\hat\eta)-V\hat\eta.
\end{equation}
Uniqueness follows because a distribution with zero derivative is constant, and a nonzero constant does not belong to $L^2(\mathbb R)$. Equivalently, one may first differentiate the dynamic equation in $x$, linearize, and then take its unique $L^2$ primitive. The zero-average periodic normalization used to construct the background is unchanged.

Combining the irrotational part of the linearization with the linearizations of $\omega\eta\psi_x$ and the vorticity primitive yields
\begin{subequations}
\begin{align}\label{linearized dynamic}
\hat\psi_t
={}&\Big(-g-B\alpha_\e\partial_x
+B\partial_x\circ\alpha_\e
-BG_{\e,\mathbb R}\circ B-\omega K_\e\circ B\Big)\hat\eta
\\ \label{linearized dynamic 2}
&+\Big(\alpha_\e\partial_x+BG_{\e,\mathbb R}
+\omega K_\e\Big)\hat\psi.
\end{align}
\end{subequations}
Here multiplication by $B$ and $\alpha_\e$ is understood as a multiplication operator, and
\begin{equation}\label{commutator identity}
-B\alpha_\e\partial_x+B\partial_x\circ\alpha_\e
=B(\alpha_\e)_x.
\end{equation}

We use the Alinhac good unknown \cite{Alinhac,LDbook}
\begin{equation}\label{alinhac}
u:=\hat\eta,\qquad v:=\hat\psi-B\hat\eta,\qquad
\begin{bmatrix}\hat\eta\\\hat\psi\end{bmatrix}
=Z\begin{bmatrix}u\\v\end{bmatrix},\qquad
Z:=\begin{bmatrix}1&0\\B&1\end{bmatrix}.
\end{equation}
Since the background Stokes wave is stationary in the moving frame, $B$ is independent of time and $v_t=\hat\psi_t-Bu_t$. Using \eqref{linearized kinematic}--\eqref{linearized dynamic 2} and
\[
B(\alpha_\e)_xu+\alpha_\e\partial_x(Bu)
-B\partial_x(\alpha_\e u)=\alpha_\e B_xu,
\]
we obtain
\begin{equation}\label{second linear eq}
\begin{bmatrix}u_t\\v_t\end{bmatrix}
=\widetilde{\mathcal A}_\e\begin{bmatrix}u\\v\end{bmatrix},\qquad
\widetilde{\mathcal A}_\e:=
\begin{bmatrix}
\partial_x\circ\alpha_\e&G_{\e,\mathbb R}\\[1mm]
-g+\alpha_\e B_x&\alpha_\e\partial_x+\omega K_\e
\end{bmatrix}:
\mathcal Y_{\mathbb R}\longrightarrow\mathcal X_{\mathbb R}.
\end{equation}

\paragraph{Conformal flattening} 
For $|\e|$ sufficiently small, there exists a unique small, real-valued, odd, and $2\pi$-periodic function $\mathfrak p_\e$ such that
\[
\phi_\e(x):=x+\mathfrak p_\e(x),\qquad x\in\mathbb R,
\]
is an orientation-preserving diffeomorphism satisfying
\[
\phi_\e(x+2\pi)=\phi_\e(x)+2\pi.
\]
We set $\rho_\e:=\phi_\e'=1+\partial_x\mathfrak p_\e$ and, after decreasing $|\e|$ if necessary, have $\inf_x\rho_\e(x)>0$. The function $\mathfrak p_\e$ is the unique small odd solution of
\begin{equation}\label{def:ttf}
\mathfrak p_\e(x)
=\mathfrak H\left[\eta_\e\big(x+\mathfrak p_\e(x)\big)\right],
\end{equation}
where $\mathfrak H$ denotes the Hilbert transform on $\mathbb T$; see \cite[formula (A.15)]{BBHM}. This equation constructs the conformal map from the periodic background. As proved in \cite{BMV3}, for every $s\geq0$ the map $\e\mapsto\mathfrak p_\e$ is analytic with values in $H^s_{\mathtt{odd}}(\mathbb T)$. The expansion in Appendix \ref{sec:App2} gives
\begin{equation}\label{expansion frak p}
\mathfrak p_\e(x)
=\e\sin(x)+\e^2\frac{(\cw^2+g)^2}{4\cw^4}\sin(2x)
+\mathcal O_{H^s}(\e^3).
\end{equation}

Define $(\mathfrak P_{\mathbb R}f)(x):=f(\phi_\e(x))$. For each fixed $s\geq0$, the operators $\mathfrak P_{\mathbb R}$ and $\mathfrak P_{\mathbb R}^{-1}$ are bounded on $H^s(\mathbb R)$. Let $\mathfrak H_{\mathbb R}$ be the Fourier multiplier with symbol $-\im\operatorname{sgn}(\xi)$, and let $|D|:=\partial_x\mathfrak H_{\mathbb R}$. The conformal representation of the real-line Dirichlet--Neumann operator is
\begin{equation}\label{eq:real-line-DNO}
G_{\e,\mathbb R}
=\partial_x\circ\mathfrak P_{\mathbb R}^{-1}
\circ\mathfrak H_{\mathbb R}\circ\mathfrak P_{\mathbb R}
:H^1(\mathbb R)\longrightarrow L^2(\mathbb R).
\end{equation}
This follows by pulling the harmonic extension back to the lower half-plane; compare \cite[Lemma A.5]{BBHM} and \cite[Section 3.2]{NS}. In particular,
\begin{equation}\label{eq:real-line-primitive}
K_\e:=\mathfrak P_{\mathbb R}^{-1}
\mathfrak H_{\mathbb R}\mathfrak P_{\mathbb R}
\end{equation}
is bounded on $L^2(\mathbb R)$ and satisfies \eqref{eq:real-line-primitive-factorization}.

We conjugate the good-unknown system by setting
\begin{equation}\label{LC}
h=\mathcal P_{\mathbb R}\begin{bmatrix}u\\v\end{bmatrix},\qquad
\mathcal P_{\mathbb R}:=
\begin{bmatrix}
\rho_\e\mathfrak P_{\mathbb R}&0\\
0&\mathfrak P_{\mathbb R}
\end{bmatrix}.
\end{equation}
This is a time-independent bounded isomorphism on both $\mathcal X_{\mathbb R}$ and $\mathcal Y_{\mathbb R}$, and
\begin{equation}\label{eq:real-line-primitive-flattening}
\mathfrak P_{\mathbb R}K_\e\mathfrak P_{\mathbb R}^{-1}
=\mathfrak H_{\mathbb R}.
\end{equation}

Define the real-valued, even, and $2\pi$-periodic coefficients
\begin{align}
\cw+p_\e(x)
&:=
\frac{\alpha_\e(\phi_\e(x))}{\rho_\e(x)}
,
\label{def:qeps}\qquad g+a_\e(x):=
\frac{
g-\alpha_\e(\phi_\e(x))B_x(\phi_\e(x))
}{\rho_\e(x)}.
\end{align}
By the analyticity of the Stokes waves and of the conformal map, the maps $\e\mapsto p_\e$ and $\e\mapsto a_\e$ are analytic with values in $H^s_{\mathtt{ev}}(\T)$ for every $s\geq0$. In Appendix \ref{sec:App2} we provide their Taylor expansions, which we collect in the following lemma.

\begin{lemma}\label{lem:pa-expansions}
The analytic functions $p_\e(x)$ and $a_\e(x)$ in \eqref{def:qeps} are even in $x$ and satisfy
\begin{equation}\label{exp:pa}
\begin{aligned}
p_\e(x)&=\e p_1(x)+\e^2p_2(x)+\cO_{H^s}(\e^3),\\
a_\e(x)&=\e a_1(x)+\e^2a_2(x)+\cO_{H^s}(\e^3),
\end{aligned}
\end{equation}
where
\begin{align}
p_1(x)&=p_1^{[1]}\cos (x),
&
 p_1^{[1]}&:=-\frac{\cw^2+g}{\cw},
\label{exp:p1}\\
p_2(x)&=p_2^{[0]}+p_2^{[2]}\cos(2x),
\label{exp:p2}
\end{align}
with
\begin{align}
p_2^{[0]}
:=
\frac{(\cw^2+g)(3\cw^4+2g\cw^2+g^2)}{8\cw^5},\qquad
p_2^{[2]}:=
\frac{(\cw^2+g)(\cw^4-4g\cw^2-g^2)}{4\cw^5}.
\label{exp:p22}
\end{align}
Moreover,
\begin{align}
a_1(x)&=a_1^{[1]}\cos (x),
&
 a_1^{[1]}&:=-(\cw^2+g),
\label{exp:a1}\\
a_2(x)&=a_2^{[0]}+a_2^{[2]}\cos(2x),
\label{exp:a2}
\end{align}
where
\begin{align}
a_2^{[0]}:=\cw^2+g,\qquad
a_2^{[2]}
:=
\frac{(\cw^2+g)(2\cw^4-3g\cw^2-g^2)}{2\cw^4}.
\label{exp:a22}
\end{align}
\end{lemma}

A direct conjugation of \eqref{second linear eq}, using \eqref{eq:real-line-DNO} and \eqref{eq:real-line-primitive-flattening}, yields
\begin{equation}\label{eq:flattened-real-line-operator}
\mathcal A_\e:=\mathcal P_{\mathbb R}\widetilde{\mathcal A}_\e
\mathcal P_{\mathbb R}^{-1}
=\begin{bmatrix}
\partial_x\circ(\cw+p_\e)&|D|\\[1mm]
-(g+a_\e)&(\cw+p_\e)\partial_x+\omega\mathfrak H_{\mathbb R}
\end{bmatrix},\qquad
D(\mathcal A_\e)=\mathcal Y_{\mathbb R}.
\end{equation}
For sufficiently small $|\e|$, one has $\inf_x(\cw+p_\e(x))>0$. The second row controls $\|\partial_xh_2\|_{L^2}$ by $C(\|\mathcal A_\e h\|_{L^2}+\|h\|_{L^2})$; the first row then controls $\|\partial_xh_1\|_{L^2}$ by the same quantity, since $\||D|h_2\|_{L^2}=\|\partial_xh_2\|_{L^2}$. Hence
\begin{equation}\label{eq:real-line-graph-norm}
\|h\|_{H^1}\leq C\big(\|\mathcal A_\e h\|_{L^2}+\|h\|_{L^2}\big),
\qquad h\in\mathcal Y_{\mathbb R},
\end{equation}
so $\mathcal A_\e$ is densely defined and closed on $\mathcal X_{\mathbb R}$.

\paragraph{Bloch--Floquet reduction.}
We use the unitary Bloch transform
\[
\mathscr B:\mathcal X_{\mathbb R}
\longrightarrow L^2([-\tfrac12,\tfrac12);L^2(\mathbb T;\mathbb C^2)).
\]
Its construction and intertwining identities are recalled in \cite[Section~2, ``Bloch--Floquet expansions'']{HLZZ2026}. Applying these identities to the periodic coefficients and Fourier multipliers of $\mathcal A_\e$ gives the direct-integral decomposition
\begin{equation}\label{eq:real-line-Bloch-decomposition}
\mathscr B\mathcal A_\e\mathscr B^{-1}
=\int_{[-1/2,1/2)}^\oplus\mathcal L_{\mu,\e}\,d\mu,
\qquad
\mathscr B\mathcal Y_{\mathbb R}
=L^2([-\tfrac12,\tfrac12);H^1(\mathbb T;\mathbb C^2)).
\end{equation}
For $\mu\ne0$, the fiber operator, formally written as $\mathcal L_{\mu,\e}=e^{-\im\mu x}\mathcal A_\e e^{\im\mu x}$, acts on periodic functions and is
\begin{equation}\label{Bloch operator nonzero mu}
\mathcal L_{\mu,\e}
=\begin{bmatrix}
(\partial_x+\im\mu)\circ(\cw+p_\e)&|D+\mu|\\[1mm]
-(g+a_\e)&(\cw+p_\e)(\partial_x+\im\mu)+\omega\mathfrak H_\mu
\end{bmatrix}
:H^1(\mathbb T;\mathbb C^2)\longrightarrow L^2(\mathbb T;\mathbb C^2),
\end{equation}
where
\[
|D+\mu|e^{\im kx}:=|k+\mu|e^{\im kx},
\qquad
\mathfrak H_\mu e^{\im kx}
:=-\im\operatorname{sgn}(k+\mu)e^{\im kx},
\qquad k\in\mathbb Z.
\]
The graph estimate \eqref{eq:real-line-graph-norm} holds uniformly for the fibers, with $\partial_x$ replaced by $\partial_x+\im\mu$, and identifies the domain of the direct-integral operator with $L^2([-\tfrac12,\tfrac12);H^1(\mathbb T;\mathbb C^2))$. Since $ \mathfrak H_\mu 1=-\im\operatorname{sgn}(\mu)$, the multiplier $\mathfrak H_\mu$ has distinct one-sided limits at $\mu=0$. The value assigned to the single fiber at $\mu=0$ does not affect \eqref{eq:real-line-Bloch-decomposition}. We use \eqref{Bloch operator nonzero mu} for $\mu\ne0$ and distinguish below its limit as $\mu\to0^+$ from the auxiliary periodic operator.

\begin{lemma}[From Bloch eigenvalues to real-line spectrum]
\label{lem:Bloch-to-real-line-spectrum}
Fix a sufficiently small $\e$. Let $I\subset(-\tfrac12,\tfrac12)\setminus\{0\}$ be an open interval and suppose that $\lambda:I\to\mathbb C$ and $v:I\to\mathcal Y$ are continuous, with
\[
\mathcal L_{\mu,\e}v(\mu)=\lambda(\mu)v(\mu),
\qquad \|v(\mu)\|_{\mathcal X}=1.
\]
Then $\lambda(\mu)\in\sigma_{\mathcal X_{\mathbb R}}(\mathcal A_\e)$ for every $\mu\in I$. In particular, any such branch with positive real part implies spectral instability on $\mathcal X_{\mathbb R}$.
\end{lemma}

\begin{proof}
Fix $\mu_0\in I$. Choose $b_n\in C_c^\infty(I)$ with $\|b_n\|_{L^2(I)}=1$ and support contained in $\{|\mu-\mu_0|<1/n\}$. Extend $b_n(\mu)v(\mu)$ by zero outside $I$ and set $h_n:=\mathscr B^{-1}(b_nv)$. Then $h_n\in\mathcal Y_{\mathbb R}$ and $\|h_n\|_{\mathcal X_{\mathbb R}}=1$. By \eqref{eq:real-line-Bloch-decomposition},
\[
\|(\mathcal A_\e-\lambda(\mu_0))h_n\|_{\mathcal X_{\mathbb R}}
\leq\sup_{\mu\in\operatorname{supp}b_n}
|\lambda(\mu)-\lambda(\mu_0)|\longrightarrow0.
\]
Thus $\lambda(\mu_0)$ belongs to the approximate point spectrum of $\mathcal A_\e$, and hence to its spectrum.
\end{proof}

\paragraph{The auxiliary periodic operator and the limit $\mu\to0^+$}
For the identities obtained by differentiating the periodic Stokes branch, we also use the linearization of the periodic problem with the zero-average primitive. Denote this auxiliary operator by $\mathcal L^{\mathrm{per}}_{\e}$. It is defined separately on $\mathcal Y\subset\mathcal X$ and is not needed to specify the single zero fiber in \eqref{eq:real-line-Bloch-decomposition}. Let $\mathfrak P f:=f\circ\phi_\e$ on $\mathbb T$. The periodic Dirichlet--Neumann operator satisfies $G(\eta_\e)=\partial_x\mathfrak P^{-1}\mathfrak H\mathfrak P$. Since the periodic primitive obeys $\partial_x^{-1}\partial_xf=f-\langle f\rangle$, its shape derivative contains the term $\langle V\hat\eta\rangle\mathbf1$, and
\[
\mathfrak P\partial_x^{-1}G(\eta_\e)\mathfrak P^{-1}w
=\mathfrak Hw-\langle\rho_\e\mathfrak Hw\rangle\mathbf1.
\]
Thus, defining
\[
\mathcal R_\e w:=\langle V(\phi_\e(\cdot))w\rangle\mathbf1,
\qquad
\mathcal M_\e w:=\langle\rho_\e\mathfrak Hw\rangle\mathbf1,
\]
the same good-unknown and conformal changes give
\begin{equation}\label{auxiliary periodic operator}
\mathcal L^{\mathrm{per}}_{\e}
=
\begin{bmatrix}
\partial_x\circ (\cw+p_\e)
&
|D|
\\[1mm]
-(g+a_\e)+\omega\mathcal R_{\e}
&
(\cw+p_\e)\partial_x
+\omega\big(\mathfrak H-\mathcal M_{\e}\big)
\end{bmatrix}.
\end{equation}

The one-sided limit of $\mathcal L_{\mu,\e}$ as $\mu\to0^+$ differs from this auxiliary periodic operator. As in \cite[Section 5.1]{NS}, we use the identity
\begin{align} \label{D+mu}
    |D+\mu|=|D|+\mu(\mathrm{sgn}(D)+\Pi_0),~~{\forall \mu\in
\left[0,\frac{1}{2}\right)},
\end{align}
where $\mathrm{sgn}(D)$ is the Fourier multiplier operator, acting on $2\pi$-periodic functions, with symbol
\begin{equation} \label{positive mu multipliers}
\mathrm{sgn}(k):=1 \ \ \forall k >0 \ , \ \
\mathrm{sgn}(0):=0  , \ \ \ 
\mathrm{sgn}(k):= -1 \ \ \forall k < 0 \  , 
\end{equation}
and $\Pi_0$ is the projector operator on the zero mode, $\Pi_0 f(x):=\frac{1}{2\pi}\int_\mathbb{T} f(x) dx$. Similarly, we have
\begin{equation}
\mathfrak H_\mu
=
\mathfrak H-\im\Pi_0,~~{\forall \mu\in\left(0,\frac{1}{2}\right)}.
\end{equation}
Consequently,
\begin{equation}\label{right zero fiber}
\mathcal L_{0^+,\e}
:=
\lim_{\mu\downarrow0}\mathcal L_{\mu,\e}
=
\begin{bmatrix}
\partial_x\circ (\cw+p_\e)
&
|D|
\\[1mm]
-(g+a_\e)
&
(\cw+p_\e)\partial_x+\omega\big(\mathfrak H-\im\Pi_0\big)
\end{bmatrix},
\end{equation}
where the limit holds in
\[
\mathcal L\left(
\big(H^1(\mathbb T;\mathbb C)\big)^2,
\big(L^2(\mathbb T;\mathbb C)\big)^2
\right).
\]
Indeed, for every $0<\mu<\frac12$,
\begin{equation}\label{analytic positive Bloch family}
\mathcal L_{\mu,\e}
=
\mathcal L_{0^+,\e}
+
\mu
\begin{bmatrix}
\im (\cw+p_\e)&\mathrm{sgn}(D)+\Pi_0\\
0&\im (\cw+p_\e)
\end{bmatrix}.
\end{equation}
Thus $\mathcal L_{\mu,\e}$ for $0<\mu<\frac12$ extends real-analytically, in fact affinely, to $\mu=0$, with value $\mathcal L_{0^+,\e}$.

Finally, comparing \eqref{auxiliary periodic operator} and \eqref{right zero fiber}, we obtain
\begin{equation}\label{finite rank jump}
\mathcal L_{0^+,\e}-\mathcal L^{\mathrm{per}}_{\e}
=
\omega
\begin{bmatrix}
0&0\\[1mm]
-\mathcal R_{\e}
&
-\im\Pi_0+\mathcal M_{\e}
\end{bmatrix}.
\end{equation}
In particular, the analytic continuation of $\mathcal L_{\mu,\e}$ from $\mu>0$ does not coincide with $\mathcal L^{\mathrm{per}}_{\e}$. Consequently, spectral perturbation for $\mu>0$ cannot be based at $\mathcal L^{\mathrm{per}}_{\e}$, but must instead be carried out from the one-sided limit $\mathcal L_{0^+,\e}$. The one-sided limit, rather than the separately normalized periodic linearization, is the base operator for the modulational spectral perturbation.

For the periodic identities below, let $\mathcal P$ denote the operator on $\mathcal Y$ given by the same formula as \eqref{LC}, with $\mathfrak P$ in place of $\mathfrak P_{\mathbb R}$, and let $Z$ be the periodic multiplication operator in \eqref{alinhac}. We write $\mathcal V_\e:=V\circ\phi_\e$, $m_\e:=\langle\mathcal V_\e\rangle$, and $\mathscr L_{0,\e}:=\mathcal L_{0^+,\e}$ (see \eqref{rmk mathscr Leps}). For $h=(h_1,h_2)^\top$, set
\begin{equation}
\label{def:jmath-section5}
\mathfrak j_\e(h)
:=
\omega\left(
-\langle\mathcal V_\e h_1\rangle
-\im\langle h_2\rangle
+\langle\rho_\e\mathfrak Hh_2\rangle
\right).
\end{equation}
Equation \eqref{finite rank jump} is equivalently
\begin{equation}
\label{jump-jmath-section5}
(\mathscr L_{0,\e}-\mathcal L^{\mathrm{per}}_{\e})h
=\mathfrak j_\e(h)\vet{0}{1}.
\end{equation}

\begin{lemma} [Generalized kernel of $\mathscr L_{0,\e}$]
\label{lem:one-sided-Jordan-section5}
Let $g>0$ and $\omega\neq0$. For $|\e|$ sufficiently small, the eigenvalue $0$ of $\mathscr L_{0,\e}$ has algebraic multiplicity three. More precisely, define
\begin{equation}
\label{adapted-zero-basis-section5}
\begin{aligned}
E_1(\e):={}&
\mathcal P Z^{-1}
\begin{bmatrix}
\partial_\e\eta_\e\\
\partial_\e\psi_\e
\end{bmatrix}-\frac{1}{\im\omega}
\left(
c_\e'm_\e-
\mathfrak j_\e\!\left(
\mathcal P Z^{-1}
\begin{bmatrix}
\partial_\e\eta_\e\\
\partial_\e\psi_\e
\end{bmatrix}
\right)
\right)\vet{0}{1},\\
E_2(\e):={}&
-\frac{\im}{\e}
\left(
\mathcal P Z^{-1}
\begin{bmatrix}
\partial_x\eta_\e\\
\partial_x\psi_\e
\end{bmatrix}
-m_\e \vet{0}{1}
\right),\\
E_3(\e):={}&
\mathcal P Z^{-1}
\begin{bmatrix}1\\0\end{bmatrix}
-\frac{1}{\im\omega}
\left(
g-\omega m_\e-
\mathfrak j_\e\!\left(
\mathcal P Z^{-1}
\begin{bmatrix}1\\0\end{bmatrix}
\right)
\right) \vet{0}{1}.
\end{aligned}
\end{equation}
Here $E_2(0)$ is defined by analytic continuation. Then
\begin{equation}
\label{exact-chain-section5}
\begin{aligned}
\mathscr L_{0,\e}E_2(\e)=0,\qquad\mathscr L_{0,\e}E_1(\e)=-\im\e c_\e' E_2(\e),\qquad
\mathscr L_{0,\e}E_3(\e)=\im\omega\e E_2(\e).
\end{aligned}
\end{equation}
\end{lemma}
\begin{proof}
We prove this Lemma in Appendix \ref{G kernel}.
\end{proof}

\paragraph{Complete local spectrum}
Recall the Benjamin--Feir coefficient $\mathrm{ind}_{\mathrm{BF}}$ defined in \eqref{intro:ThetaBF}. Since $\cw>0$ and $\omega=(\cw^2-g)/\cw$, we have
\begin{equation}
\label{equiv-main-threshold-section2}
\mathrm{ind}_{\mathrm{BF}}>0
\quad\Longleftrightarrow\quad
\cw^2<3g
\quad\Longleftrightarrow\quad
\omega<2\sqrt{\frac g3},
\end{equation}
with all inequalities reversed when $\mathrm{ind}_{\mathrm{BF}}<0$.

\begin{theorem}[Complete Benjamin--Feir spectrum]
\label{Complete BF thm}
Let $g>0$ and
\[
\omega\in\mathbb R\setminus
\left\{0,2\sqrt{\frac g3}\right\}.
\]
There exist $\e_0,\mu_0,r_*>0$ such that, for $0<\e<\e_0$ and $0<|\mu|<\mu_0$, the operator $\mathcal L_{\mu,\e}$ has exactly three eigenvalues in $\{|\lambda|<r_*\}$, counted with algebraic multiplicity.

In terms of the analytic blocks $A^{(3)}$ and $d^{(3)}$ in \eqref{def:M3-section5}, these eigenvalues are
\begin{equation}
\label{main-lambda0-section2}
\lambda_0(\mu,\e)
=
\im\cw\mu+\mu d^{(3)}(|\mu|,\e)
\in\im\mathbb R
\end{equation}
and the unordered pair
\begin{equation}
\label{main-lambda-pm-unstable-section2}
\lambda_\pm(\mu,\e)
=
\im\cw\mu
+\frac{\mu}{2}\operatorname{tr}A^{(3)}(|\mu|,\e)
\pm\frac{|\mu|}{2}
\sqrt{\Delta_{\mathrm{BF}}(|\mu|,\e)}.
\end{equation}
Here $\Delta_{\mathrm{BF}}$ is the real analytic discriminant in \eqref{exp:DeltaBF-section5}, with expansion
\begin{equation}
\label{main-Delta-exp-section2}
\Delta_{\mathrm{BF}}(\mu,\e)
=
\mathrm{ind}_{\mathrm{BF}}\e^2
-\frac{4\cw^6g^4}{(\cw^2+g)^6}\mu^2
+\mathcal O\big((|\mu|+|\e|)^3\big).
\end{equation}
For $\Delta_{\mathrm{BF}}<0$, the square root in \eqref{main-lambda-pm-unstable-section2} is understood as $\im\sqrt{-\Delta_{\mathrm{BF}}}$. The smallness thresholds and remainder bounds may depend on the fixed pair $(g,\omega)$.

If $\mathrm{ind}_{\mathrm{BF}}>0$, there is a real-analytic collision curve
\begin{equation}
\label{main-mustar-section2}
\mu_*(\e)
=
\frac{(\cw^2+g)^3}{2\cw^3g^2}
\sqrt{\mathrm{ind}_{\mathrm{BF}}}\,\e
+\mathcal O(\e^2).
\end{equation}
After decreasing $\e_0$ and $\mu_0$ if necessary, $0<\mu_*(\e)<\mu_0$ and
\begin{equation}
\label{main-sign-classification-section2}
\begin{cases}
\Delta_{\mathrm{BF}}(|\mu|,\e)>0,
&0<|\mu|<\mu_*(\e),
\\[1mm]
\Delta_{\mathrm{BF}}(|\mu|,\e)=0,
&|\mu|=\mu_*(\e),
\\[1mm]
\Delta_{\mathrm{BF}}(|\mu|,\e)<0,
&\mu_*(\e)<|\mu|<\mu_0.
\end{cases}
\end{equation}
Thus $\lambda_\pm$ have opposite nonzero real parts for $0<|\mu|<\mu_*(\e)$, collide on the imaginary axis at $|\mu|=\mu_*(\e)$, and are purely imaginary for $\mu_*(\e)<|\mu|<\mu_0$. The positive collision point is unique in $(0,\mu_0)$.

After defining $\lambda_\pm(0,\e)=0$ by continuity, the pair for $|\mu|\le\mu_*(\e)$ forms a figure eight: the branches for $0\le\mu\le\mu_*(\e)$ form a simple closed upper lobe, and those for $-\mu_*(\e)\le\mu\le0$ form its complex-conjugate lower lobe. The two lobes meet only at the origin.

If $\mathrm{ind}_{\mathrm{BF}}<0$, then $\Delta_{\mathrm{BF}}(|\mu|,\e)<0$ throughout the stated parameter range, and all three local eigenvalues are purely imaginary.
\end{theorem}

\begin{proof}
For $\mu>0$, the spectral reduction \eqref{reduced-spectrum-original-L} and Theorem~\ref{thm:exact-eigenvalues-section5} give the eigenvalue count and formulas \eqref{main-lambda0-section2}--%
\eqref{main-lambda-pm-unstable-section2}. The discriminant expansion is \eqref{exp:DeltaBF-section5}. Theorem~\ref{thm:exact-eigenvalues-section5} gives the collision curve, its uniqueness, the sign classification, and the figure-eight geometry, as well as the purely imaginary alternative.

For $\mu<0$, Lemma~\ref{lem:negative-bloch-symmetry} gives
\[
\sigma(\mathcal L_{\mu,\e})
=
\overline{\sigma(\mathcal L_{|\mu|,\e})}.
\]
Since $A^{(3)}$ and $d^{(3)}$ are purely imaginary for real parameters, this yields the stated formulas as an unordered eigenvalue pair and gives the complex-conjugate lower lobe. Finally, \eqref{equiv-main-threshold-section2} expresses the two alternatives in terms of the vorticity.
\end{proof}

\section{Kato transformation}
\label{sec:Kato-positive}
We work on the complex Hilbert space
\begin{equation}
\label{def:X-Y-section3}
\mathcal X
:=L^2(\mathbb T;\mathbb C^2),
\qquad
\mathcal Y
:=H^1(\mathbb T;\mathbb C^2),
\end{equation}
equipped with the scalar product
\begin{equation}
\label{def:scalar-product-section3}
(f,h)_{\mathcal X}
:=
\frac1{2\pi}
\int_0^{2\pi}
\left(
f_1(x)\overline{h_1(x)}
+
f_2(x)\overline{h_2(x)}
\right)\,dx,
\end{equation}
which is linear in the first argument.

We write, for $0<\mu<\frac12$,
\begin{equation}
\label{def:shifted-scriptL}
\mathcal L_{\mu,\e}
=
\im\cw\mu\,\mathrm{Id}+\mathscr L_{\mu,\e}.
\end{equation}
In view of \eqref{analytic positive Bloch family}, the shifted operator is
\begin{equation}
\label{eq:scriptL-family}
\mathscr L_{\mu,\e}
=
\mathscr L_{0,\e}
+
\mu\mathscr A_\e,
\end{equation}
where
\begin{equation}
\label{eq:scriptL-zero-plus}
\mathscr L_{0,\e}
:=
\begin{bmatrix}
\partial_x\circ(\cw+p_\e)
&
|D|
\\[1mm]
-(g+a_\e)
&
(\cw+p_\e)\partial_x
+\omega(\mathfrak H-\im\Pi_0)
\end{bmatrix}
\end{equation}
and
\begin{equation}
\label{eq:A-eps}
\mathscr A_\e
:=
\begin{bmatrix}
\im p_\e&\mathrm{sgn}(D)+\Pi_0\\
0&\im p_\e
\end{bmatrix}.
\end{equation}
Thus
\begin{equation}
\label{spectral-shift-scriptL}
\sigma(\mathcal L_{\mu,\e})
=
\im\cw\mu+\sigma(\mathscr L_{\mu,\e}).
\end{equation}

Although \eqref{eq:scriptL-family} has been obtained from $\mathcal L_{\mu,\e}$ for $\mu>0$, its right-hand side defines an analytic extension to $\mu$ in a neighborhood of the origin. We shall use this extension only as a perturbative device; for $\mu<0$ it is not identified with $\mathcal L_{\mu,\e}$. Here and below, $\mathscr L_{\mu,\e}$ denotes this analytic extension. Its value at $\mu=0$ is
\begin{align} \label{rmk mathscr Leps}
\mathscr L_{0,\e}=\mathcal L_{0^+,\e}.
\end{align}
The same convention applies to $P_{\mu,\e}$ and $U_{\mu,\e}$, defined below as well as their adjoints and the inverse of $U_{\mu,\e}$. All parameter derivatives at the origin are taken within this analytic family.

\begin{lemma}[The operators at $\mu$ and $-\mu$]
\label{lem:negative-bloch-symmetry}
Let $\e\in\mathbb R$ be sufficiently small and let $0<|\mu|<\frac12$.  Define complex conjugation by
\[
\mathfrak K:L^2(\mathbb T;\mathbb C^2)
\longrightarrow
L^2(\mathbb T;\mathbb C^2),
\qquad
\mathfrak K h:=\overline h.
\]
Then
\begin{equation}
\label{negative-bloch-conjugacy}
\mathfrak K\mathcal L_{\mu,\e}
=
\mathcal L_{-\mu,\e}\mathfrak K
\end{equation}
and therefore
\begin{equation}
\label{negative-bloch-spectrum-symmetry}
\sigma(\mathcal L_{-\mu,\e})
=
\overline{\sigma(\mathcal L_{\mu,\e})}.
\end{equation}
The equality in \eqref{negative-bloch-spectrum-symmetry} preserves algebraic multiplicity.  In particular, every local eigenvalue branch of $\mathcal L_{-\mu,\e}$ is the complex conjugate of the corresponding branch of $\mathcal L_{\mu,\e}$ for $0<\mu<\frac12$.
\end{lemma}

\begin{proof}
Since $p_\e$ and $a_\e$ are real-valued, ordinary complex conjugation commutes with multiplication by these coefficients.  Moreover, for every $f\in H^1(\mathbb T;\mathbb C)$, $\overline{(\partial_x+\im\mu)f} = (\partial_x-\im\mu)\bar f.$ If $f(x)=\sum_{k\in\mathbb Z}f_ke^{\im kx}$, then
\[
\overline{|D+\mu|f}
=
\sum_{k\in\mathbb Z}|k+\mu|\,\overline{f_k}e^{-\im kx}
=
|D-\mu|\,\overline f,
\]
because $|-k-\mu|=|k+\mu|$.  Similarly,
\begin{align*}
\overline{\mathfrak H_\mu f}
&=
\sum_{k\in\mathbb Z}
\im\,\operatorname{sign}(k+\mu)\,\overline{f_k}e^{-\im kx}=
\mathfrak H_{-\mu}\overline f,
\end{align*}
since $\operatorname{sign}(-k-\mu)=-\operatorname{sign}(k+\mu)$. Substituting these identities into \eqref{Bloch operator nonzero mu} proves \eqref{negative-bloch-conjugacy}.

If $(\mathcal L_{\mu,\e}-\lambda)^m h=0$, then \eqref{negative-bloch-conjugacy} gives
\[
(\mathcal L_{-\mu,\e}-\overline\lambda)^m\mathfrak K h=0.
\]
Since $\mathfrak K$ is bijective, it maps the full generalized eigenspace of $\lambda$ for $\mathcal L_{\mu,\e}$ by complex conjugation and bijectively onto the full generalized eigenspace of $\overline\lambda$ for $\mathcal L_{-\mu,\e}$.  This proves both the spectral symmetry and the statement about algebraic multiplicities.
\end{proof}

\paragraph{The subspace at $(0,0)$} At the flat surface one has $p_0=a_0=0$, and hence
\begin{equation}
\label{def:scriptL-star}
\mathscr L_0
:=
\mathscr L_{0,0}
=
\mathcal L_{0^+,0}
=
\begin{bmatrix}
\cw\partial_x&|D|
\\[1mm]
-g&
\cw\partial_x+\omega(\mathfrak H-\im\Pi_0)
\end{bmatrix}
:
\mathcal Y\subset\mathcal X\longrightarrow\mathcal X.
\end{equation}
The following proposition describes completely the spectral structure of $\mathscr L_0$ at the origin.

\begin{proposition}[The eigenspace at zero]
\label{prop:reference-spectrum}
The operator $\mathscr L_0$ has compact resolvent. Its eigenvalue $\lambda=0$ is isolated and semisimple, with
\begin{equation}
\label{alg-geo-mult-three}
\dim\ker\mathscr L_0
=
\dim\ker(\mathscr L_0)^m
=
3,
\qquad
m\geq1.
\end{equation}
More precisely,
\begin{equation}
\label{right-reference-basis}
\ker\mathscr L_0
=
\operatorname{span}_{\mathbb C}
\{f_+,f_-,f_0\},
\end{equation}
where the eigenvectors $f_+$, $f_-$, and $f_0$ are defined by
\begin{equation}
\label{def:right-reference-basis}
f_+(x)
:=
e^{\im x}
\begin{bmatrix}
1\\ -\im\cw
\end{bmatrix},
\qquad
f_-(x)
:=
e^{-\im x}
\begin{bmatrix}
1\\ \im\cw
\end{bmatrix},
\qquad
f_0
:=
\begin{bmatrix}
1\\ \displaystyle\frac{\im g}{\omega}
\end{bmatrix}.
\end{equation}
The adjoint operator is
\begin{equation}
\label{def:scriptL-star-adjoint}
\mathscr L_0^*
=
\begin{bmatrix}
-\cw\partial_x&-g
\\[1mm]
|D|&
-\cw\partial_x-\omega(\mathfrak H-\im\Pi_0)
\end{bmatrix},
\end{equation}
and
\begin{equation}
\label{left-reference-basis}
\ker\mathscr L_0^*
=
\operatorname{span}_{\mathbb C}
\{g_+,g_-,g_0\},
\end{equation}
where the eigenvectors $g_+$, $g_-$, and $g_0$ are defined by
\begin{equation}
\label{def:left-reference-basis}
g_+(x)
:=
\frac{e^{\im x}}{\cw^2+g}
\begin{bmatrix}
g\\ -\im\cw
\end{bmatrix},
\qquad
g_-(x)
:=
\frac{e^{-\im x}}{\cw^2+g}
\begin{bmatrix}
g\\ \im\cw
\end{bmatrix},
\qquad
g_0
:=
\begin{bmatrix}
1\\0
\end{bmatrix}.
\end{equation}
These eigenvectors satisfy
\begin{equation}
\label{biorth-reference}
(f_j,g_k)_{\mathcal X}
=
\delta_{jk},
\qquad
j,k\in\{+,-,0\}.
\end{equation}
\end{proposition}

\begin{proof}
The operator $\mathscr L_0$ is diagonal with respect to the Fourier decomposition. For $k\neq0$, its restriction to the subspace $e^{\im kx}\mathbb C^2$ is represented by
\begin{equation}
\label{Fourier-block-Ak}
A_k=
\begin{bmatrix}
\im\cw k&|k|
\\
-g&
\im\cw k-\im\omega\,\operatorname{sign}(k)
\end{bmatrix}.
\end{equation}
The two eigenvalues of $A_k$ are
\begin{equation}
\label{flat-eigenvalues-nonzero-k}
\lambda_{k}^{\pm}
=\im
\left(
\cw k-\frac{\omega}{2}\mathrm{sgn}(k)
\pm
\frac12\sqrt{\omega^2+4g |k|}
\right)\qquad k\neq 0.
\end{equation}
Indeed,
\begin{equation}
\label{det-Ak-zero}
\det A_k
=
|k|
\big(
g+\omega\cw-\cw^2|k|
\big)
=
\cw^2|k|(1-|k|),
\end{equation}
where we used \eqref{intro:speed}. Hence a zero eigenvalue can occur for $k\neq0$ only when $k=\pm1$. In each of these two Fourier blocks the zero eigenvalue is simple, and the corresponding eigenvectors are $f_+$ and $f_-$ in \eqref{def:right-reference-basis}.

On the constant Fourier mode one has
\begin{equation}
\label{constant-block-A0}
A_0
=
\begin{bmatrix}
0&0\\
-g&-\im\omega
\end{bmatrix}.
\end{equation}
Since $\omega\neq0$, its eigenvalues are $0$ and $-\im\omega$, and the zero eigenvalue is simple with eigenvector $f_0$. This proves \eqref{right-reference-basis}. Moreover, each of the three zero eigenvalues occurs in a distinct simple Fourier block, so zero is semisimple and \eqref{alg-geo-mult-three} follows.

Formula \eqref{flat-eigenvalues-nonzero-k}, together with \eqref{constant-block-A0}, also shows that the remainder of the spectrum is separated from the origin. In addition, the inverse Fourier blocks of $\mathscr L_0-\lambda$ are of order $|k|^{-1}$ as $|k|\to\infty$ whenever $\lambda$ belongs to the resolvent set. Hence the resolvent maps $\mathcal X$ boundedly into $\mathcal Y$, and is compact as an operator on $\mathcal X$ by the compact embedding $\mathcal Y\hookrightarrow\mathcal X$.

Since $\partial_x^*=-\partial_x$, $|D|^*=|D|$, $\mathfrak H^*=-\mathfrak H$, and $\Pi_0^*=\Pi_0$, the adjoint is \eqref{def:scriptL-star-adjoint}. A direct calculation in the Fourier blocks $k=1$, $k=-1$, and $k=0$ yields \eqref{def:left-reference-basis}. Finally, \eqref{biorth-reference} follows directly from the explicit formulas for $f_j$ and $g_j$.
\end{proof}
\paragraph{Kato projectors} Set
\begin{equation}
\label{def:spectral-gap-star}
d_0
:=
\operatorname{dist}
\left(
0,\,
\sigma(\mathscr L_0)\setminus\{0\}
\right)
>0
\end{equation}
and let
\begin{equation}
\label{def:Gamma-star}
\Gamma
:=
\left\{
\lambda\in\mathbb C:
|\lambda|=\frac{d_0}{2}
\right\},
\end{equation}
oriented counterclockwise.

\begin{lemma}[Kato projectors]
\label{lem:Riesz-Kato-positive}
There exist $\mu_0,\e_0>0$ such that, for $|\mu|<\mu_0$ and $|\e|<\e_0$, the following statements hold.

\begin{enumerate}[label=\textup{(\roman*)}]
\item
$\Gamma$ is contained in the resolvent set of $\mathscr L_{\mu,\e}$.

\item
The projectors
\begin{equation}
\label{def:Riesz-P-mu-eps}
P_{\mu,\e}
:=
-\frac{1}{2\pi\im}
\oint_\Gamma
(\mathscr L_{\mu,\e}-\lambda)^{-1}\,d\lambda
\end{equation}
are well-defined projectors commuting with $\mathscr{L}_{\mu,\e}$, i.e. $P^2_{\mu,\e}=P_{\mu,\e}$ and $P_{\mu,\e}\mathscr L_{\mu,\e} = \mathscr L_{\mu,\e}P_{\mu,\e}$. The map $(\mu,\e)\mapsto P_{\mu,\e}$ is analytic, both as an $\mathcal L(\mathcal X)$-valued map and as an $\mathcal L(\mathcal X,\mathcal Y)$-valued map. Set $P_0:=P_{0,0}$ and $\mathcal V_0:=\operatorname{Ran}P_0$. By Proposition \ref{prop:reference-spectrum},
\[
\mathcal V_0=\ker\mathscr L_0,
\qquad
\operatorname{Ran}P_0^*=\ker\mathscr L_0^*.
\]

\item
The subspace
\begin{equation}
\label{def:V-mu-eps}
\mathcal V_{\mu,\e}
:=
\operatorname{Ran}P_{\mu,\e}
\end{equation}
has dimension three, and
\[
\sigma(\mathscr L_{\mu,\e})
\cap\operatorname{int}(\Gamma)
=
\sigma\!\left(
\mathscr L_{\mu,\e}\big|_{\mathcal V_{\mu,\e}}
\right).
\]

\item
For $(\mu,\e)$ sufficiently small, define
\begin{equation}
\label{def:Kato-U}
U_{\mu,\e}
:=
\left(
\mathrm{Id}-(P_{\mu,\e}-P_0)^2
\right)^{-1/2}
\left[
P_{\mu,\e}P_0
+
(\mathrm{Id}-P_{\mu,\e})(\mathrm{Id}-P_0)
\right].
\end{equation}
Then $U_{\mu,\e}$ depends analytically on $(\mu,\e)$, is bounded and invertible on both $\mathcal X$ and $\mathcal Y$ with inverse
\begin{align*}
    U^{-1}_{\mu,\e}
=\left[
P_0 P_{\mu,\e}
+
(\mathrm{Id}-P_{0})(\mathrm{Id}-P_{\mu,\e})
\right] \left(
\mathrm{Id}-(P_{\mu,\e}-P_0)^2
\right)^{-1/2}.
\end{align*}
At the reference point, $U_{0,0}=\mathrm{Id}$. Moreover, $U_{\mu,\e}P_0U_{\mu,\e}^{-1} = P_{\mu,\e}$, as well as $U^{-1}_{\mu,\e}P_{\mu,\e}U_{\mu,\e} = P_{0}$ and hence
\begin{equation}
\label{Kato-subspaces}
\mathcal V_{\mu,\e}
=
U_{\mu,\e}\mathcal V_0.
\end{equation}
\end{enumerate}
\end{lemma}

\begin{proof}
We observe that $\mathscr L_{\mu,\e}$ is an analytic family with fixed domain $\mathcal Y$. Also, Proposition \ref{prop:reference-spectrum} shows that $\Gamma$ is contained in the resolvent set of $\mathscr L_0$. The resolvent identity and a Neumann-series argument imply that $\Gamma$ remains in the resolvent set for $(\mu,\e)$ sufficiently small, uniformly in $\lambda\in\Gamma$. The analyticity and the spectral properties of $P_{\mu,\e}$ then follow from the standard theory of spectral projectors; see \cite[Chapter I, \S\S 4--6 and Chapter II, \S 4]{Kato1966}. Since $P_{\mu,\e}$ is close to $P_0$ in operator norm, its rank is locally constant, and therefore equals $\operatorname{rank}P_0=3$.

Finally, after decreasing $\mu_0$ and $\e_0$ if necessary, one has $\|P_{\mu,\e}-P_0\|<1$. In \eqref{def:Kato-U} we use $\big(\mathrm{Id}-(P_{\mu,\e}-P_0)^2\big)^{-1/2} =\sum_{n=0}^\infty 4^{-n}\binom{2n}{n}(P_{\mu,\e}-P_0)^{2n}$. The remaining assertions are the standard properties of $U_{\mu,\e}$ in \eqref{def:Kato-U}.
\end{proof}

\paragraph{Bases and the reduced matrix} The absence of a symplectic structure on $\mathcal V_{\mu,\e}$ does not
affect the construction of the Kato transformation. What changes is the coordinate system used to represent the reduced operator.

For $j\in\{+,-,0\}$ define
\begin{equation}
\label{def:right-Kato-basis}
f_j(\mu,\e)
:=
U_{\mu,\e}f_j
\in\mathcal V_{\mu,\e}
\end{equation}
and
\begin{equation}
\label{def:left-Kato-basis}
g_j(\mu,\e)
:=
\left(U_{\mu,\e}^{-1}\right)^*g_j
\in\operatorname{Ran}P_{\mu,\e}^*.
\end{equation}

\begin{lemma}
\label{lem:biorthogonal-Kato}
For $(\mu,\e)$ sufficiently small,
\begin{align} \label{mathcal F}
\mathcal F_{\mu,\e}
:=
\{f_+(\mu,\e),f_-(\mu,\e),f_0(\mu,\e)\}
\end{align}
is a basis of $\mathcal V_{\mu,\e}$ and
\[
\mathcal G_{\mu,\e}
:=
\{g_+(\mu,\e),g_-(\mu,\e),g_0(\mu,\e)\}
\]
is a basis of $\operatorname{Ran}P_{\mu,\e}^*$. They satisfy
\begin{equation}
\label{moving-biorthogonality}
\big(
f_j(\mu,\e),
g_k(\mu,\e)
\big)_{\mathcal X}
=
\delta_{jk},
\qquad
j,k\in\{+,-,0\}.
\end{equation}
Moreover, all six vectors depend analytically on $(\mu,\e)$.
\end{lemma}

\begin{proof}
The basis property follows from \eqref{Kato-subspaces} and the invertibility of $U_{\mu,\e}$. Taking the adjoint of the identity $P_{\mu,\e}=U_{\mu,\e}P_0U_{\mu,\e}^{-1}$ in Lemma \ref{lem:Riesz-Kato-positive} gives
\[
P_{\mu,\e}^*
=
\left(U_{\mu,\e}^{-1}\right)^*
P_0^*
U_{\mu,\e}^*,
\]
which proves that the vectors in \eqref{def:left-Kato-basis} belong to $\operatorname{Ran}P_{\mu,\e}^*$. Finally,
\[
\big(
U_{\mu,\e}f_j,
(U_{\mu,\e}^{-1})^*g_k
\big)_{\mathcal X}
=
(f_j,g_k)_{\mathcal X}
=
\delta_{jk}.
\]
Analyticity follows from Lemma \ref{lem:Riesz-Kato-positive}.
\end{proof}

Every $h\in\mathcal V_{\mu,\e}$ therefore has the exact expansion
\begin{equation}
\label{biorthogonal-expansion}
h
=
\sum_{j\in\{+,-,0\}}
\big(
h,g_j(\mu,\e)
\big)_{\mathcal X}
f_j(\mu,\e).
\end{equation}
We define the reduced matrix
\begin{equation}
\label{def:reduced-matrix-M}
\mathbf M(\mu,\e)
=
\big(M_{jk}(\mu,\e)\big)_{j,k\in\{+,-,0\}}
\in\mathbb C^{3\times3}
\end{equation}
by
\begin{equation}
\label{def:entries-M}
M_{jk}(\mu,\e)
:=
\big(
\mathscr L_{\mu,\e}f_k(\mu,\e),
g_j(\mu,\e)
\big)_{\mathcal X},
\qquad
j,k\in\{+,-,0\}.
\end{equation}
The order of the basis is always understood to be $(+,-,0)$. By Lemma \ref{lem:Riesz-Kato-positive} and the matrix representation above,
\begin{equation}
\label{reduced-spectrum-scriptL}
\sigma(\mathscr L_{\mu,\e})
\cap\operatorname{int}(\Gamma)
=
\sigma(\mathbf M(\mu,\e)),
\end{equation}
counting algebraic multiplicities. For $0<\mu<\mu_0$, the spectral shift \eqref{spectral-shift-scriptL} consequently gives
\begin{equation}
\label{reduced-spectrum-original-L}
\sigma(\mathcal L_{\mu,\e})
\cap\left(\im\cw\mu+\operatorname{int}(\Gamma)\right)
=
\im\cw\mu+\sigma(\mathbf M(\mu,\e)).
\end{equation}



\begin{remark}[The case $\omega=0$]
\label{rem:omega-zero-singular}
The restriction $\omega\neq0$ is structural rather than technical. If $\omega=0$, the constant Fourier block \eqref{constant-block-A0} becomes
\[
\begin{bmatrix}
0&0\\-g&0
\end{bmatrix},
\]
so that zero has geometric multiplicity three but algebraic multiplicity four. Thus the three-dimensional semisimple reduction constructed above does not extend through $\omega=0$. At $\omega=0$ one recovers the familiar four-dimensional Jordan structure and the Hamiltonian reduction of the pure-gravity problem; see, for instance, \cite{BMV1}.
\end{remark}

\section{Matrix representation on $\mathcal V_{\mu,\e}$}
\label{sec:matrix-representation-positive}

In this section we compute the quadratic Taylor expansion of the matrix \(\mathbf M(\mu,\e)\) introduced in \eqref{def:reduced-matrix-M}.  The calculation is carried out in the bases constructed in Section \ref{sec:Kato-positive}.  Throughout the section we use the ordered right and left bases
\[
\big(f_+(\mu,\e),f_-(\mu,\e),f_0(\mu,\e)\big),
\qquad
\big(g_+(\mu,\e),g_-(\mu,\e),g_0(\mu,\e)\big).
\]
We write \(\mathcal O_3\) for a \(3\times3\) matrix-valued analytic function whose norm is bounded by \(C(|\mu|+|\e|)^3\) near \((\mu,\e)=(0,0)\). More generally, $\mathcal O_j$ denotes an analytic remainder of total order at least $j$ in $(\mu,\e)$, in the scalar or matrix space under consideration. We use a prime to denote $\partial_\e$ and a dot to denote $\partial_\mu$, both evaluated at $(\mu,\e)=(0,0)$. By \eqref{eq:scriptL-family}--\eqref{eq:A-eps} and Lemma \ref{lem:pa-expansions},
\begin{equation}
\label{first-jets-scriptL}
\mathscr L_0'
=
\begin{bmatrix}
\partial_x\circ p_1&0\\
-a_1&p_1\partial_x
\end{bmatrix},
\qquad
\dot{\mathscr L}_0
=
\begin{bmatrix}
0&\mathrm{sgn}(D)+\Pi_0\\
0&0
\end{bmatrix},
\qquad
\dot{\mathscr L}_0'
=
\begin{bmatrix}
\im p_1&0\\
0&\im p_1
\end{bmatrix},
\end{equation}
where $p_1$ and $a_1$ are defined in Lemma \ref{lem:pa-expansions}.

We first collect, in a form convenient for the matrix computation, the expansions obtained in Appendix \ref{sec:AppKato}.

\begin{lemma}[Expansion of the bases]
\label{lem:main-biorthogonal-Kato-expansion}
As \((\mu,\e)\to(0,0)\), the vectors $U_{\mu,\e}f_j$ satisfy
\begin{align}
 f_+(\mu,\e)
 ={}&
 e^{\im x}\begin{bmatrix}1\\-\im\cw\end{bmatrix}
 +\e e^{2\im x}
 \begin{bmatrix}
 \dfrac{(\cw^2+g)^2}{2\cw^4}\\[2mm]
 -\dfrac{\im g(\cw^2+g)}{2\cw^3}
 \end{bmatrix}
 +\mu e^{\im x}
 \begin{bmatrix}
 \dfrac{\cw^4}{(\cw^2+g)^2}\\[2mm]
 \dfrac{\im g\cw^3}{(\cw^2+g)^2}
 \end{bmatrix}
 +\mathcal O_{\mathcal Y}\big((|\mu|+|\e|)^2\big),
 \label{main:exp-fplus}
 \\
 f_-(\mu,\e)
 ={}&
 e^{-\im x}\begin{bmatrix}1\\\im\cw\end{bmatrix}
 +\e e^{-2\im x}
 \begin{bmatrix}
 \dfrac{(\cw^2+g)^2}{2\cw^4}\\[2mm]
 \dfrac{\im g(\cw^2+g)}{2\cw^3}
 \end{bmatrix}
 +\mu e^{-\im x}
 \begin{bmatrix}
 -\dfrac{\cw^4}{(\cw^2+g)^2}\\[2mm]
 \dfrac{\im g\cw^3}{(\cw^2+g)^2}
 \end{bmatrix}
 +\mathcal O_{\mathcal Y}\big((|\mu|+|\e|)^2\big),
 \label{main:exp-fminus}
 \\
 f_0(\mu,\e)
 ={}&
 \begin{bmatrix}1\\ \dfrac{\im g}{\omega}\end{bmatrix}
 +\e
 \begin{bmatrix}
 \dfrac{2\cw^2}{\cw^2+g}\cos(x)\\[2mm]
 -\dfrac{2\cw g}{\cw^2+g}\sin(x)
 \end{bmatrix}
 +\mu
 \begin{bmatrix}
 0\\[1mm]-\dfrac{\im g^2}{\omega^3}
 \end{bmatrix}
 +\mathcal O_{\mathcal Y}\big((|\mu|+|\e|)^2\big).
 \label{main:exp-fzero}
\end{align}
The vectors $(U_{\mu,\e}^{-1})^*g_j$ satisfy
\begin{align}
 g_+(\mu,\e)
 ={}&
 \frac{e^{\im x}}{\cw^2+g}
 \begin{bmatrix}g\\-\im\cw\end{bmatrix}
 +\e e^{2\im x}
 \begin{bmatrix}
 \dfrac{(\cw^2+g)^2-2\cw^4}{4\cw^4}\\[2mm]
 -\dfrac{\im(\cw^2+g)}{2\cw^3}
 \end{bmatrix}
 +\mu e^{\im x}
 \begin{bmatrix}
 -\dfrac{g^2\cw^2}{(\cw^2+g)^3}\\[2mm]
 -\dfrac{\im g^2\cw}{(\cw^2+g)^3}
 \end{bmatrix}
 +\mathcal O_{\mathcal X}\big((|\mu|+|\e|)^2\big),
 \label{main:exp-gplus}
 \\
g_-(\mu,\e)
 ={}&
 \frac{e^{-\im x}}{\cw^2+g}
 \begin{bmatrix}g\\\im\cw\end{bmatrix}
 +\e e^{-2\im x}
 \begin{bmatrix}
 \dfrac{(\cw^2+g)^2-2\cw^4}{4\cw^4}\\[2mm]
 \dfrac{\im(\cw^2+g)}{2\cw^3}
 \end{bmatrix}
 +\mu e^{-\im x}
 \begin{bmatrix}
 \dfrac{g^2\cw^2}{(\cw^2+g)^3}\\[2mm]
 -\dfrac{\im g^2\cw}{(\cw^2+g)^3}
 \end{bmatrix}
 +\mathcal O_{\mathcal X}\big((|\mu|+|\e|)^2\big),
 \label{main:exp-gminus}
 \\
 g_0(\mu,\e)
 ={}&
 \begin{bmatrix}1\\0\end{bmatrix}
 +\mu
 \begin{bmatrix}
 -\dfrac{g}{\omega^2}\\[2mm]
 \dfrac{\im}{\omega}
 \end{bmatrix}
 +\mathcal O_{\mathcal X}
 \big(|\mu|(|\mu|+|\e|)\big).
 \label{main:exp-gzero}
\end{align}
Moreover,
\begin{align} \label{g0 is zero}
    g_0(0,\e)=\vet{1}{0},
\end{align}
and
\begin{equation}
\label{main:moving-biorthogonality}
\big(f_j(\mu,\e),g_k(\mu,\e)\big)_{\mathcal X}=\delta_{jk},
\qquad j,k\in\{+,-,0\}.
\end{equation}
\end{lemma}

\begin{proof}
The long computations of the expansions are given in Appendix \ref{sec:AppKato}. The exact biorthogonality \eqref{main:moving-biorthogonality} is \eqref{moving-biorthogonality}. Finally, let us prove \eqref{g0 is zero}. Given $h=(h_1,h_2)=\in\mathcal{Y}$, we have
\begin{align}
    \left(\mathscr{L}_{0,\e} h,\vet{1}{0}\right)=\langle\pa_x\circ(\cw+p_\e)h_1+|D|h_2\rangle=0,
\end{align}
and, hence, $\mathscr{L}_{0,\e}^*\vet{1}{0}=0$. Recalling \eqref{def:Riesz-P-mu-eps}, we obtain
\begin{align*}
    P^*_{0,\e} \vet{1}{0}=P^*_{0} \vet{1}{0}=\vet{1}{0}, 
\end{align*}
and, hence,
\begin{align*}
    g_0(0,\e)=(U^{-1}_{0,\e})^* \vet{1}{0}=\left(
\mathrm{Id}-(P^*_{0,\e}-P^*_0)^2
\right)^{-1/2}
\left[
P^*_{0,\e}P^*_0
+
(\mathrm{Id}-P^*_{0,\e})(\mathrm{Id}-P^*_0)
\right]=\vet{1}{0}.
\end{align*}
\end{proof}

We now state the main result of this section.

\begin{lemma}\label{prop:quadratic-M-section4}
The matrix representing \(\mathscr L_{\mu,\e}|_{\mathcal V_{\mu,\e}}\) in the bases defined in \eqref{def:right-Kato-basis} and \eqref{def:left-Kato-basis} satisfies
\begin{equation}
\label{final:M-quadratic-section4}
\mathbf M(\mu,\e)
=
\im
\begin{pmatrix}
-\frac{\cw g}{\cw^2+g}\mu-\Theta\e^2+\frac{\cw^3g^2}{(\cw^2+g)^3}\mu^2
&-\Theta\e^2
&\dfrac{\omega}{2}\e+\beta_+\mu\e
\\[2mm]
\Theta\e^2
&-\frac{\cw g}{\cw^2+g}\mu+\Theta\e^2-\frac{\cw^3g^2}{(\cw^2+g)^3}\mu^2
&-\dfrac{\omega}{2}\e+\beta_-\mu\e
\\[2mm]
-\dfrac{\cw^2+g}{2\cw}\mu\e
&-\dfrac{\cw^2+g}{2\cw}\mu\e
&\frac{g}{\omega}\mu-\frac{g^2}{\omega^3}\mu^2
\end{pmatrix}
+\mathcal O_3,
\end{equation}
where the coefficients $\Theta$, $\beta_{\pm}$ are given as follows:
\begin{equation}
\label{def:coefficients-section4}
\begin{aligned}
\Theta&:=
\frac{3\cw^6-3g\cw^4+3g^2\cw^2+g^3}{8\cw^5}.
\end{aligned}
\end{equation}
and
\begin{equation}
\label{def:beta-pm-section4}
\begin{aligned}
\beta_+:=
\frac g2\left(
\frac1\omega-
\frac{\cw(2g+\omega^2)}{(\cw^2+g)^2}
\right),\qquad \beta_-:=
-\frac g2\left(
\frac1\omega+
\frac{\cw(2g+\omega^2)}{(\cw^2+g)^2}
\right).
\end{aligned}
\end{equation}
\end{lemma}

By \eqref{eq:scriptL-family}--\eqref{eq:A-eps} and \eqref{first-jets-scriptL}, we have
\begin{equation}
\label{decomp:scriptL-section4}
\mathscr L_{\mu,\e}
=
\mathscr L_{0,\e}
+\mu\dot{\mathscr L}_0
+\im\mu p_\e\,\mathrm{Id}_2.
\end{equation}
Accordingly, define the three matrices
\begin{align}
\big(\mathbf M_\e(\mu)\big)_{jk}
&:=
\big(\mathscr L_{0,\e}f_k(\mu,\e),g_j(\mu,\e)\big)_{\mathcal X},
\label{def:M-eps-section4}\\
\big(\mathbf M^\flat(\mu,\e)\big)_{jk}
&:=
\mu\big(\dot{\mathscr L}_0 f_k(\mu,\e),g_j(\mu,\e)\big)_{\mathcal X},
\label{def:M-flat-section4}\\
\big(\mathbf M^\sharp(\mu,\e)\big)_{jk}
&:=
\mu\big(\im p_\e f_k(\mu,\e),g_j(\mu,\e)\big)_{\mathcal X}.
\label{def:M-sharp-section4}
\end{align}
Thus,
\begin{equation}
\label{decomp:M-section4}
\mathbf M(\mu,\e)
=
\mathbf M_\e(\mu)
+\mathbf M^\flat(\mu,\e)
+\mathbf M^\sharp(\mu,\e).
\end{equation}

\paragraph{Expansion of matrices}

The amplitude derivatives satisfy
\begin{equation}
\label{amplitude-derivatives-section4}
\frac{1}{r!}\left.\partial_\e^r\mathscr L_{0,\e}\right|_{\e=0}
=
\begin{bmatrix}
\partial_x\circ p_r&0\\
-a_r&p_r\partial_x
\end{bmatrix},
\qquad r=1,2,
\end{equation}
so that
\begin{equation}
\label{exp:L0plus-eps-section4}
\mathscr L_{0,\e}
=
\mathscr L_0+\e\mathscr L_0'+\frac{\e^2}{2}\mathscr L_0''
+\mathcal O_{\mathcal L(\mathcal Y,\mathcal X)}(\e^3).
\end{equation}

\begin{lemma}[Expansion of \(\mathbf M_\e\)]
\label{lem:exp-M-eps-section4}
The matrix \(\mathbf M_\e(\mu)\) defined in \eqref{def:M-eps-section4} has the expansion
\begin{equation}
\label{exp:M-eps-section4}
\mathbf M_\e(\mu)
=
\im
\begin{pmatrix}
-\Theta\e^2-\dfrac{\cw^3g^2}{(\cw^2+g)^3}\mu^2
&-\Theta\e^2
&\dfrac{\omega}{2}\e\\[2mm]
\Theta\e^2
&\Theta\e^2+\dfrac{\cw^3g^2}{(\cw^2+g)^3}\mu^2
&-\dfrac{\omega}{2}\e\\[2mm]
0&0&\dfrac{g^2}{\omega^3}\mu^2
\end{pmatrix}
+\mathcal O_3.
\end{equation}
In particular, the coefficient of \(\mu\e\) in \(\mathbf M_\e(\mu)\) vanishes identically.
\end{lemma}

\begin{proof}
In this proof, all basis vectors and their derivatives are evaluated at $(\mu,\e)=(0,0)$; a prime denotes $\partial_\e$ and a dot denotes $\partial_\mu$. Since \(\mathscr L_0f_k=0\) and \(\mathscr L_0^*g_j=0\), the terms containing second derivatives of $U_{\mu,\e}f_k$ and $(U_{\mu,\e}^{-1})^*g_j$ do not contribute to the quadratic Taylor polynomial.  Expanding \eqref{def:M-eps-section4} by means of \eqref{exp:L0plus-eps-section4} gives the following four coefficient matrices. First,
\begin{equation}
\label{calc:M-eps-order-eps}
\Big( (\mathscr L_0'f_k,g_j)_{\mathcal X}\Big)_{j,k}
=
\im
\begin{pmatrix}
0&0&\dfrac{\omega}{2}\\
0&0&-\dfrac{\omega}{2}\\
0&0&0
\end{pmatrix}.
\end{equation}
Second,
\begin{equation}
\label{calc:M-eps-order-mu2}
\Big( (\mathscr L_0\dot f_k,\dot g_j)_{\mathcal X}\Big)_{j,k}
=
\im
\begin{pmatrix}
-\dfrac{\cw^3g^2}{(\cw^2+g)^3}&0&0\\[2mm]
0&\dfrac{\cw^3g^2}{(\cw^2+g)^3}&0\\[2mm]
0&0&\dfrac{g^2}{\omega^3}
\end{pmatrix}.
\end{equation}
For the pure \(\e^2\)-coefficient, one obtains
\begin{align}
&\Big(
(\mathscr L_0f_k',g_j')_{\mathcal X}
+(\tfrac12\mathscr L_0''f_k,g_j)_{\mathcal X}
+(\mathscr L_0'f_k',g_j)_{\mathcal X}
+(\mathscr L_0'f_k,g_j')_{\mathcal X}
\Big)_{j,k}=
\im\Theta
\begin{pmatrix}
-1&-1&0\\
1&1&0\\
0&0&0
\end{pmatrix}.
\label{calc:M-eps-order-eps2}
\end{align}
Finally, the mixed coefficient is
\begin{align}
&\Big(
(\mathscr L_0f_k',\dot g_j)_{\mathcal X}
+(\mathscr L_0\dot f_k,g_j')_{\mathcal X}
+(\mathscr L_0'\dot f_k,g_j)_{\mathcal X}
+(\mathscr L_0'f_k,\dot g_j)_{\mathcal X}
\Big)_{j,k}
=0.
\label{calc:M-eps-order-mueps}
\end{align}
All four identities follow by substituting the six expansions from Lemma \ref{lem:main-biorthogonal-Kato-expansion}, using \eqref{exp:p1}--\eqref{exp:a22}, and applying Fourier orthogonality.  No uncomputed basis coefficient enters these formulas.  Summing \eqref{calc:M-eps-order-eps}--\eqref{calc:M-eps-order-mueps} proves \eqref{exp:M-eps-section4}.
\end{proof}

\begin{lemma}[Expansion of \(\mathbf M^\flat\)]
\label{lem:exp-M-flat-section4}
The matrix \(\mathbf M^\flat(\mu,\e)\) defined in \eqref{def:M-flat-section4} satisfies
\begin{equation}
\label{exp:M-flat-section4}
\mathbf M^\flat(\mu,\e)
=
\im
\begin{pmatrix}
-\frac{\cw g}{\cw^2+g}\mu+\frac{2\cw^3g^2}{(\cw^2+g)^3}\mu^2
&0
&\frac{\cw g^2}{(\cw^2+g)^2}\mu\e\\[1mm]
0
&-\frac{\cw g}{\cw^2+g}\mu-\frac{2\cw^3g^2}{(\cw^2+g)^3}\mu^2
&\frac{\cw g^2}{(\cw^2+g)^2}\mu\e\\[1mm]
0&0&\frac{g}{\omega}\mu-\frac{2g^2}{\omega^3}\mu^2
\end{pmatrix}
+\mathcal O_3.
\end{equation}
\end{lemma}

\begin{proof}
Because the prefactor \(\mu\) is explicit in \eqref{def:M-flat-section4}, only the expansions of the two bases to first order are needed.  Taylor expansion gives
\begin{align*}
\mathbf M^\flat(\mu,\e)
=
\mu\Big((\dot{\mathscr L}_0f_k,g_j)_{\mathcal X}\Big)_{j,k}+\mu^2\Big(
(\dot{\mathscr L}_0\dot f_k,g_j)_{\mathcal X}
+(\dot{\mathscr L}_0f_k,\dot g_j)_{\mathcal X}
\Big)_{j,k}+\mu\e\Big(
(\dot{\mathscr L}_0f_k',g_j)_{\mathcal X}
+(\dot{\mathscr L}_0f_k,g_j')_{\mathcal X}
\Big)_{j,k}
+\mathcal O_3.
\end{align*}
The three coefficient matrices are, respectively,
\begin{align}
\Big((\dot{\mathscr L}_0f_k,g_j)_{\mathcal X}\Big)_{j,k}
&=
\im
\begin{pmatrix}
-\frac{\cw g}{\cw^2+g}&0&0\\
0&-\frac{\cw g}{\cw^2+g}&0\\
0&0&\frac{g}{\omega}
\end{pmatrix},
\label{calc:M-flat-order-mu}\\[1mm]
\Big(
(\dot{\mathscr L}_0\dot f_k,g_j)_{\mathcal X}
+(\dot{\mathscr L}_0f_k,\dot g_j)_{\mathcal X}
\Big)_{j,k}
&=
2\im
\begin{pmatrix}
\frac{\cw^3g^2}{(\cw^2+g)^3}&0&0\\
0&-\frac{\cw^3g^2}{(\cw^2+g)^3}&0\\
0&0&-\frac{g^2}{\omega^3}
\end{pmatrix},
\label{calc:M-flat-order-mu2}\\[1mm]
\Big(
(\dot{\mathscr L}_0f_k',g_j)_{\mathcal X}
+(\dot{\mathscr L}_0f_k,g_j')_{\mathcal X}
\Big)_{j,k}
&=
\im\frac{\cw g^2}{(\cw^2+g)^2}
\begin{pmatrix}
0&0&1\\
0&0&1\\
0&0&0
\end{pmatrix}.
\label{calc:M-flat-order-mueps}
\end{align}
In conclusion, \eqref{exp:M-flat-section4} is proven.
\end{proof}

\begin{lemma}[Expansion of \(\mathbf M^\sharp\)]
\label{lem:exp-M-sharp-section4}
The matrix \(\mathbf M^\sharp(\mu,\e)\) defined in \eqref{def:M-sharp-section4} satisfies
\begin{equation}
\label{exp:M-sharp-section4}
\mathbf M^\sharp(\mu,\e)
=
\im\mu\e
\begin{pmatrix}
0&0&\dfrac{g^2}{2\cw(\cw^2-g)}\\[3mm]
0&0&-\dfrac{g(2\cw^2-g)}{2\cw(\cw^2-g)}\\[3mm]
-\dfrac{\cw^2+g}{2\cw}&-\dfrac{\cw^2+g}{2\cw}&0
\end{pmatrix}
+\mathcal O_3.
\end{equation}
\end{lemma}

\begin{proof}
By \eqref{exp:pa}, $p_\e=\e p_1+\mathcal O_{H^s}(\e^2)$. Since \(\mathbf M^\sharp\) already carries the prefactor \(\mu\), the quadratic Taylor coefficient only uses the eigenvectors:
\begin{equation}
\label{calc:M-sharp-reduction}
\mathbf M^\sharp(\mu,\e)
=
\mu\e
\Big((\im p_1f_k,g_j)_{\mathcal X}\Big)_{j,k}
+\mathcal O_3.
\end{equation}
Recalling Lemma~\ref{lem:pa-expansions}, \eqref{def:right-reference-basis}, and \eqref{biorth-reference}, a direct calculation gives
\begin{equation}
\label{calc:M-sharp-coefficient}
\Big((\im p_1f_k,g_j)_{\mathcal X}\Big)_{j,k}
=
\im
\begin{pmatrix}
0&0&\dfrac{g^2}{2\cw(\cw^2-g)}\\[3mm]
0&0&-\dfrac{g(2\cw^2-g)}{2\cw(\cw^2-g)}\\[3mm]
-\dfrac{\cw^2+g}{2\cw}&-\dfrac{\cw^2+g}{2\cw}&0
\end{pmatrix}.
\end{equation}
This is precisely \eqref{exp:M-sharp-section4}.
\end{proof}
Summing \eqref{exp:M-eps-section4}, \eqref{exp:M-flat-section4}, and \eqref{exp:M-sharp-section4} proves Lemma \ref{prop:quadratic-M-section4}.

\section{Block triangularization and Benjamin--Feir spectrum}
\label{sec:characteristic-reduction-positive}

\paragraph{Properties of $\mathbf{M}(0,\e)$} Our goal is to triangularize of the matrix $\mathbf M_{\mu,\e}$ defined in Lemma \ref{prop:quadratic-M-section4}. Let us denote
\begin{align}
     \mathbf M(\mu,\e)=\begin{pmatrix}
   A & B\\
    C & d
\end{pmatrix}, \quad\text{with}\quad A:=\begin{pmatrix}
   \mathbf{M}_{11} & \mathbf{M}_{12}\\
    \mathbf{M}_{21} & \mathbf{M}_{22}
\end{pmatrix},\quad B:=\pvet{\mathbf{M}_{13}}{\mathbf{M}_{23}},\quad C:=\begin{pmatrix}
    \mathbf{M}_{31} & \mathbf{M}_{32}
\end{pmatrix},\quad d:=\mathbf{M}_{33}.
\end{align}
\begin{lemma}
\label{lem:reversibility-section5}
For real $\mu$ and $\e$ sufficiently small, the matrix $\mathbf{M}(\mu,\e)$ is purely imaginary.
\end{lemma}

\begin{proof}
For real $(\mu,\e)$, define the antiunitary involution
\[
(\mathfrak S h)(x):=
\begin{bmatrix}
\overline{h_1(-x)}\\
-\overline{h_2(-x)}
\end{bmatrix}.
\]
Since $p_\e$ and $a_\e$ are real-valued and even, we have $\mathfrak S\mathscr L_{\mu,\e} =-\mathscr L_{\mu,\e}\mathfrak S$. The invariance of $\Gamma$ under $\lambda\mapsto-\overline\lambda$ implies that $P_{\mu,\e}$, and hence $U_{\mu,\e}$, commute with $\mathfrak S$. The reference vectors $f_j,g_j$ are fixed by $\mathfrak S$; therefore the moving bases satisfy
\[
\mathfrak S f_j(\mu,\e)=f_j(\mu,\e),
\qquad
\mathfrak S g_j(\mu,\e)=g_j(\mu,\e).
\]
Consequently, using $(\mathfrak S u,\mathfrak S v)_{\mathcal X} =\overline{(u,v)_{\mathcal X}}$,
\[
\overline{\mathbf M_{jk}(\mu,\e)}
=
\big(\mathfrak S\mathscr L_{\mu,\e}f_k(\mu,\e),
\mathfrak Sg_j(\mu,\e)\big)_{\mathcal X}
=
-\mathbf M_{jk}(\mu,\e),
\]
which proves that every entry of $\mathbf M(\mu,\e)$ is purely imaginary.
\end{proof}

\begin{lemma}
The matrix $\mathbf{M}(\mu,\e)$ satisfies $(\mathbf{M}(0,\e))^2=0$, $C(0,\e)=0$, and $d(0,\e)=0$.
\end{lemma}
\begin{proof}
By \eqref{exact-chain-section5}, each vector is annihilated by $\mathscr L_{0,\e}^{2}$ on $\mathcal V_{0,\e}$. By definition of $\mathbf{M}(\mu,\e)$ and recalling \eqref{mathcal F}, we know $\mathscr{L}_{0,\e} \mathcal{F}_{0,\e}=\mathcal{F}_{0,\e}\mathbf{M}(0,\e)$, and, hence $\mathscr{L}^2_{0,\e} \mathcal{F}_{0,\e}=\mathcal{F}_{0,\e}(\mathbf{M}(0,\e))^2$. Therefore, we prove $(\mathbf{M}(0,\e))^2=0$. Next, recalling the fact $g_0(0,\e)=\vet{1}{0}$ (see \eqref{g0 is zero}), we obtain
\begin{align*}
\mathbf{M}_{3k}(0,\e)=(\mathscr{L}_{0,\e} f_k(0,\e),g_0(0,\e))=\left(\mathscr{L}_{0,\e} f_k(0,\e),\vet{1}{0}\right)=0.
\end{align*}
This proves that $C(0,\e)=0$ and $d(0,\e)=0$.
\end{proof}

\paragraph{A change of coordinates at $\mu=0$}
We introduce
\begin{equation}
\label{def:R-epsilon-section5}
R_\e:=
\begin{pmatrix}
\dfrac12&
\dfrac{\mathbf M_{13}(0,\e)}{\im\omega\e}&0\\[2mm]
\dfrac12&
\dfrac{\mathbf M_{23}(0,\e)}{\im\omega\e}&0\\[2mm]
0&0&1
\end{pmatrix},
\end{equation}
where the quotients at $\e=0$ are defined by analytic continuation, and set
\begin{equation}
\label{def:M1-section5}
\mathbf M^{(1)}(\mu,\e)
:=R_\e^{-1}\mathbf M(\mu,\e)R_\e.
\end{equation}

\begin{lemma}
\label{lem:M1-properties-section5}
For $|\e|$ sufficiently small, the map $\e\mapsto R_\e$ is analytic and $R_\e$ is invertible. There exists an analytic scalar function $a(\e)$ such that
\begin{equation}
\label{eq:M1-at-zero-section5}
\mathbf M^{(1)}(0,\e)
=
\begin{pmatrix}
0&0&0\\
a(\e)&0&\im\omega\e\\
0&0&0
\end{pmatrix}.
\end{equation}
Moreover,
\begin{equation}
\label{eq:M1-mu-derivatives-section5}
\partial_\mu\mathbf M^{(1)}_{12}(0,\e)
=
\partial_\mu\mathbf M^{(1)}_{32}(0,\e)
=0.
\end{equation}
\end{lemma}

The proof of Lemma~\ref{lem:M1-properties-section5} is given in Appendix~\ref{M1-properties}. Combining Lemma~\ref{prop:quadratic-M-section4}, \eqref{def:R-epsilon-section5}, and Lemma~\ref{lem:M1-properties-section5}, we write
\begin{equation}
\label{def:M1-blocks-section5}
\mathbf M^{(1)}(\mu,\e)
=
\begin{pmatrix}
A^{(1)}(\mu,\e)&B^{(1)}(\mu,\e)\\
C^{(1)}(\mu,\e)&d^{(1)}(\mu,\e)
\end{pmatrix},
\end{equation}
with the following expansions:
\begin{align}
\label{eq:A1-expansion-section5}
A^{(1)}(\mu,\e)
&=\im
\begin{pmatrix}
-\dfrac{\cw g}{\cw^2+g}\mu
+\cO\big(|\mu|(|\mu|+|\e|)^2\big)
&
\dfrac{\cw^3g^2}{(\cw^2+g)^3}\mu^2
+\cO\big(\mu^2(|\mu|+|\e|)\big)
\\[3mm]
\dfrac{\cw^3g^2}{(\cw^2+g)^3}\mu^2-2\Theta\e^2
+\cO\big((|\mu|+|\e|)^3\big)
&
-\dfrac{\cw g}{\cw^2+g}\mu
+\cO\big(|\mu|(|\mu|+|\e|)^2\big)
\end{pmatrix},
\\[3mm]
\label{eq:B1-expansion-section5}
B^{(1)}(\mu,\e)
&=\im
\begin{pmatrix}
-\dfrac{g\cw(2g+\omega^2)}{(\cw^2+g)^2}\mu\e
+\cO\big(|\mu|(|\mu|+|\e|)^2\big)
\\[3mm]
\omega\e+\dfrac g\omega\mu\e
+\cO\big((|\mu|+|\e|)^3\big)
\end{pmatrix},
\\[3mm]
\label{eq:C1-expansion-section5}
C^{(1)}(\mu,\e)
&=\im
\begin{pmatrix}
-\dfrac{\cw^2+g}{2\cw}\mu\e
+\cO\big(|\mu|(|\mu|+|\e|)^2\big)
&
\cO\big(\mu^2(|\mu|+|\e|)\big)
\end{pmatrix},
\\[3mm]
\label{eq:d1-expansion-section5}
d^{(1)}(\mu,\e)
&=\im
\dfrac g\omega\mu-\im\dfrac{g^2}{\omega^3}\mu^2
+\im\cO\big(|\mu|(|\mu|+|\e|)^2\big)
.
\end{align}
Here $\Theta$ is defined in \eqref{def:coefficients-section4}, and all remainders are real analytic for real $(\mu,\e)$ and satisfy the indicated bounds.

Indeed, $R_\e=R_0+\cO(\e^2)$, so the quadratic Taylor coefficients follow from Lemma~\ref{prop:quadratic-M-section4}. The factors of $\mu$ in the remainders follow from \eqref{eq:M1-at-zero-section5}; for the $(1,2)$ and $(3,2)$ entries, \eqref{eq:M1-mu-derivatives-section5} gives the additional factor of $\mu$.

\paragraph{Exact block triangularization} We first perform a singular change of coordinates. For $\mu\neq0$,
define
\begin{equation}
\label{def:Q-M2-section5}
Q(\mu):=\operatorname{diag}(\mu,1,\mu),
\qquad
\mathbf M^{(2)}(\mu,\e)
:=\frac1\mu Q(\mu)^{-1}\mathbf M^{(1)}(\mu,\e)Q(\mu).
\end{equation}
The singular rescaling in \eqref{def:Q-M2-section5} is analogous to the one used by Sun and Wahl\'en \cite[Eq.~(5.18)]{SunWahlen2025} in their modulation analysis of capillary--gravity waves. By Lemma~\ref{lem:M1-properties-section5}, this matrix extends analytically to $\mu=0$. We write
\begin{equation}
\label{def:M2-blocks-section5}
\mathbf M^{(2)}(\mu,\e)
=
\begin{pmatrix}
A^{(2)}(\mu,\e)&B^{(2)}(\mu,\e)\\
C^{(2)}(\mu,\e)&d^{(2)}(\mu,\e)
\end{pmatrix},
\end{equation}
By \eqref{eq:A1-expansion-section5}--\eqref{eq:d1-expansion-section5} and \eqref{def:Q-M2-section5}, the blocks satisfy
\begin{align}
\label{eq:A2-expansion-section5}
A^{(2)}(\mu,\e)
&=\im
\begin{pmatrix}
-\dfrac{\cw g}{\cw^2+g}
+\cO\big((|\mu|+|\e|)^2\big)
&
\dfrac{\cw^3g^2}{(\cw^2+g)^3}
+\cO(|\mu|+|\e|)
\\[3mm]
\dfrac{\cw^3g^2}{(\cw^2+g)^3}\mu^2-2\Theta\e^2
+\cO\big((|\mu|+|\e|)^3\big)
&
-\dfrac{\cw g}{\cw^2+g}
+\cO\big((|\mu|+|\e|)^2\big)
\end{pmatrix},
\\[3mm]
\label{eq:B2-expansion-section5}
B^{(2)}(\mu,\e)
&=\im
\begin{pmatrix}
-\dfrac{g\cw(2g+\omega^2)}{(\cw^2+g)^2}\e
+\cO\big((|\mu|+|\e|)^2\big)
\\[3mm]
\omega\e+\cO\big((|\mu|+|\e|)^2\big)
\end{pmatrix},
\\[3mm]
\label{eq:C2-expansion-section5}
C^{(2)}(\mu,\e)
&=\im
\begin{pmatrix}
-\dfrac{\cw^2+g}{2\cw}\e
+\cO\big((|\mu|+|\e|)^2\big)
&
\cO(|\mu|+|\e|)
\end{pmatrix},
\\[3mm]
\label{eq:d2-expansion-section5}
d^{(2)}(\mu,\e)
&=\im
\dfrac g\omega+\im\cO(|\mu|+|\e|).
\end{align}
All remainders are real analytic for real $(\mu,\e)$ and satisfy the indicated bounds near $(0,0)$. In particular, we obtain
\begin{equation}
\label{eq:M2-origin-section5}
\mathbf M^{(2)}(0,0)
=
\im
\begin{pmatrix}
-\dfrac{\cw g}{\cw^2+g}
&
\dfrac{\cw^3g^2}{(\cw^2+g)^3}
&0\\[2mm]
0&-\dfrac{\cw g}{\cw^2+g}&0\\[2mm]
0&0&\dfrac g\omega
\end{pmatrix}.
\end{equation}
Also, this matrix is nonsingular:
\[
\det\mathbf M^{(2)}(0,0)
=-\frac{\im\cw^2g^3}{\omega(\cw^2+g)^2}\neq0.
\]

The block triangularization below is analogous to the similarity transformation of Sun and Wahl\'en \cite[Eqs.~(5.31)--(5.32)]{SunWahlen2025}, which eliminates an off-diagonal block of their $2+2$ matrix.

\begin{lemma}[Exact block triangularization]
\label{lem:exact-block-triangularization}
There exists a real row vector $X$, analytic in $(\mu,\e)$ near $(0,0)$, with $X(0,0)=0$ and
\begin{equation}
\label{eq:X-leading-section5}
\begin{aligned}
X(\mu,\e)
&:=\begin{pmatrix}x_1(\mu,\e)&x_2(\mu,\e)\end{pmatrix},\\
x_1(\mu,\e)
&=\frac{\omega(\cw^2+g)^2}{4\cw^3g}\e
+\cO\big((|\mu|+|\e|)^2\big),\\
x_2(\mu,\e)&=\cO(|\mu|+|\e|),
\end{aligned}
\end{equation}
such that the invertible matrix
\begin{equation}
\label{def:T-section5}
T(\mu,\e):=
\begin{pmatrix}
I_2&0\\
X(\mu,\e)&1
\end{pmatrix}
\end{equation}
conjugates $\mathbf M^{(2)}$ to the block upper triangular matrix
\begin{equation}
\label{def:M3-section5}
\mathbf M^{(3)}(\mu,\e)
:=T(\mu,\e)^{-1}\mathbf M^{(2)}(\mu,\e)T(\mu,\e)
=
\begin{pmatrix}
A^{(3)}(\mu,\e)&B^{(3)}(\mu,\e)\\
0&d^{(3)}(\mu,\e)
\end{pmatrix},
\end{equation}
where
\begin{equation}
\label{eq:M3-blocks-section5}
A^{(3)}=A^{(2)}+B^{(2)}X,
\qquad
B^{(3)}=B^{(2)},
\qquad
d^{(3)}=d^{(2)}-XB^{(2)}.
\end{equation}
The solution $X$ is locally unique among row vectors close to zero.
\end{lemma}

\begin{proof}
Since
\[
T^{-1}=
\begin{pmatrix}
I_2&0\\
-X&1
\end{pmatrix},
\]
a direct multiplication gives
\[
T^{-1}\mathbf M^{(2)}T
=
\begin{pmatrix}
A^{(2)}+B^{(2)}X&B^{(2)}\\
C^{(2)}+d^{(2)}X-XA^{(2)}-XB^{(2)}X&
d^{(2)}-XB^{(2)}
\end{pmatrix}.
\]
Thus, exact block triangularization amounts to solving
\begin{equation}
\label{eq:Riccati-section5}
F(X,\mu,\e)
:=
C^{(2)}(\mu,\e)+d^{(2)}(\mu,\e)X
-XA^{(2)}(\mu,\e)-XB^{(2)}(\mu,\e)X
=0.
\end{equation}
The entries of $\mathbf M^{(2)}$ are purely imaginary for real $(\mu,\e)$, so $-\im F$ is a real analytic map for real $X$. By \eqref{eq:M2-origin-section5}, $F(0,0,0)=0$, and
\begin{equation}
\label{eq:Riccati-linearization-section5}
D_XF(0,0,0)[H]
=
\im H
\begin{pmatrix}
\dfrac{2\cw^2g}{\omega(\cw^2+g)}
&
-\dfrac{\cw^3g^2}{(\cw^2+g)^3}\\[2mm]
0&
\dfrac{2\cw^2g}{\omega(\cw^2+g)}
\end{pmatrix},
\qquad H\in\mathbb R^{1\times2}.
\end{equation}
Here we used $\cw^2-\omega\cw-g=0$ to simplify
\[
\frac g\omega+\frac{\cw g}{\cw^2+g}
=\frac{2\cw^2g}{\omega(\cw^2+g)}\neq0.
\]
Consequently, $D_X(-\im F)(0,0,0)$ is invertible. The analytic implicit function theorem gives a unique small real analytic solution $X=X(\mu,\e)$ with $X(0,0)=0$. In particular,
\[
X(\mu,\e)=\cO(|\mu|+|\e|).
\]

Using \eqref{eq:A2-expansion-section5}--\eqref{eq:d2-expansion-section5} and this bound, the first component of \eqref{eq:Riccati-section5}, multiplied by $-\im$, gives
\[
0=
-\frac{\cw^2+g}{2\cw}\e
+\frac{2\cw^2g}{\omega(\cw^2+g)}x_1(\mu,\e)
+\cO\big((|\mu|+|\e|)^2\big).
\]
Hence
\[
x_1(\mu,\e)
=
\frac{\omega(\cw^2+g)^2}{4\cw^3g}\e
+\cO\big((|\mu|+|\e|)^2\big),
\]
which proves \eqref{eq:X-leading-section5}. Substitution into the block multiplication above proves \eqref{def:M3-section5}--\eqref{eq:M3-blocks-section5}.
\end{proof}

\paragraph{The Benjamin--Feir instability} Since $X$ is real and $\mathbf M^{(2)}$ is purely imaginary,
both $A^{(3)}$ and $d^{(3)}$ are purely imaginary. By \eqref{eq:M3-blocks-section5} and the preceding expansions,
\begin{equation}
\label{eq:trace-A3-section5}
\operatorname{tr}A^{(3)}(\mu,\e)
=
-\frac{2\im\cw g}{\cw^2+g}
+\im\cO\big((|\mu|+|\e|)^2\big).
\end{equation}
We define
\begin{equation}
\label{ThetaBF-section5}
\mathrm{ind}_{\mathrm{BF}}
:=
\frac{8\cw^3g^2}{(\cw^2+g)^3}\Theta
-\frac{\omega^2g}{\cw^2+g}
=
\frac{g(3g-\cw^2)(\cw^4+g^2)}{(\cw^2+g)^3}.
\end{equation}
Then the real analytic discriminant satisfies
\begin{equation}
\label{exp:DeltaBF-section5}
\begin{aligned}
\Delta_{\mathrm{BF}}(\mu,\e)
&:=
\big(\operatorname{tr}A^{(3)}(\mu,\e)\big)^2
-4\det A^{(3)}(\mu,\e)\\
&=
\mathrm{ind}_{\mathrm{BF}}\e^2
-\frac{4\cw^6g^4}{(\cw^2+g)^6}\mu^2
+\cO\big((|\mu|+|\e|)^3\big).
\end{aligned}
\end{equation}
Indeed, using the first-order term of $x_1$,
\[
A^{(3)}_{21}
=
\im\left[
\frac{\cw^3g^2}{(\cw^2+g)^3}\mu^2
-2\Theta\e^2
+\frac{\omega^2(\cw^2+g)^2}{4\cw^3g}\e^2
\right]
+\cO\big((|\mu|+|\e|)^3\big).
\]
Together with the expansions of $A^{(3)}_{12}$ and the diagonal entries, the identity
\[
\Delta_{\mathrm{BF}}
=
\big(A^{(3)}_{11}-A^{(3)}_{22}\big)^2
+4A^{(3)}_{12}A^{(3)}_{21}
\]
gives \eqref{exp:DeltaBF-section5}.

\begin{theorem}
\label{thm:exact-eigenvalues-section5}
For sufficiently small $0<\mu<\mu_0$ and $|\e|<\e_0$, the three eigenvalues of $\mathcal L_{\mu,\e}$ close to zero, counting algebraic multiplicities, are
\begin{equation}
\label{lambda0-exact-section5}
\lambda_0(\mu,\e)
=
\im\cw\mu+\mu d^{(3)}(\mu,\e)
\in\im\mathbb R,
\end{equation}
and
\begin{equation}
\label{lambda-pm-unstable-form-section5}
\lambda_\pm(\mu,\e)
=
\im\cw\mu
+\frac{\mu}{2}\operatorname{tr}A^{(3)}(\mu,\e)
\pm\frac{\mu}{2}\sqrt{\Delta_{\mathrm{BF}}(\mu,\e)}.
\end{equation}
When $\Delta_{\mathrm{BF}}<0$, the square root is understood as $\im\sqrt{-\Delta_{\mathrm{BF}}}$, and both eigenvalues are purely imaginary.
\end{theorem}

\begin{proof}
The block upper triangular matrix $\mathbf M^{(3)}$ has eigenvalues $d^{(3)}$ and the two roots
\[
\frac12\operatorname{tr}A^{(3)}
\pm\frac12\sqrt{\Delta_{\mathrm{BF}}}.
\]
By the definition of $\mathbf M^{(2)}$ and the similarity transformations, $\sigma(\mathbf M)=\mu\sigma(\mathbf M^{(3)})$. Adding the shift $\im\cw\mu$ gives the stated formulas.
\end{proof}

\section{The first isola away from the origin}
\label{sec:first-isola}

We now consider the spectrum near a nonzero imaginary eigenvalue of $\mathcal L_{\mu,0}$. Two Fourier modes whose indices differ by two meet at this eigenvalue; this is the $p=2$ resonance. We shall show that, for nonzero vorticity in the range specified below, their interaction produces a closed curve of unstable eigenvalues. Both its width and its maximal growth rate are of order $\e^2$.

The argument follows the two-dimensional reduction used by Berti, Maspero and Ventura \cite{BMV_isola}; see also \cite{BCMV_isolas}. We use the biorthogonal coordinates of the preceding sections and determine the two eigenvalues from the trace and determinant of the resulting matrix. All wave and operator coefficients needed here are already contained in Lemma \ref{lem:pa-expansions}.

In this section, we only consider the parameter region
\begin{equation} \label{isola:range}
\cw^2>\frac{2g}{3},
\qquad\text{equivalently}\qquad
\omega>-\sqrt{\frac g6}.
\end{equation}

\paragraph{The eigenvalues at $\e=0$} Recall that $\cw>0$ and $\cw^2-\omega \cw-g=0$. For $\mu\ne0$,
\eqref{Bloch operator nonzero mu} reads
\begin{equation}\label{isola:operator}
\mathcal L_{\mu,\e}=
\begin{bmatrix}
(\partial_x+\im\mu)\circ(\cw+p_\e)&|D+\mu|\\
-(g+a_\e)&(\cw+p_\e)(\partial_x+\im\mu)+\omega\mathfrak H_\mu
\end{bmatrix}.
\end{equation}
Throughout this section, $\mu\in[-\tfrac12,\tfrac12)$. Local intervals are understood on $\mathbb R/\mathbb Z$, with $-\tfrac12$ and $\tfrac12$ identified.

At $\e=0$, each Fourier mode is invariant. If $\xi=k+\mu\ne0$ is its physical wave number, the corresponding matrix and its two frequencies are
\begin{equation}\label{isola:flat-block}
\begin{bmatrix}
\im \cw\xi&|\xi|\\-g&\im \cw\xi-\im\omega\operatorname{sign}\xi
\end{bmatrix},
\qquad
s_\pm(\xi):=\cw\xi-\frac\omega2\operatorname{sign}\xi
\pm\sqrt{\frac{\omega^2}{4}+g|\xi|}.
\end{equation}
The temporal eigenvalues are $\im s_\pm(\xi)$. We seek a collision between the upper frequency at a positive wave number and the lower frequency at the wave number two units above it.

\begin{lemma}[The first nonzero collision]\label{lem:isola-flat}
Suppose \eqref{isola:range} holds. Set
\begin{equation}\label{isola:anchor}
\xi_*:=\frac{3\cw^2-2g}{4g},\quad
s_*:=\frac{3\cw^3}{4g},\quad
j_*:=\left\lfloor\xi_*+\frac12\right\rfloor,\quad
\bar\mu:=\xi_*-j_*\in\left[-\frac12,\frac12\right).
\end{equation}
If $\bar\mu\ne0$, then $\im s_*$ is an isolated semisimple eigenvalue of $\mathcal L_{\bar\mu,0}$ of algebraic multiplicity two. Its two physical wave numbers are $\xi_*$ and $\xi_*+2$, and
\begin{equation}\label{isola:collision}
s_+(\xi_*)=s_-(\xi_*+2)=s_*.
\end{equation}
\end{lemma}

\begin{proof}
For $\xi>0$, equation $s_+(\xi)=s_-(\xi+2)$ is equivalent to
\[
\sqrt{\omega^2/4+g\xi}+
\sqrt{\omega^2/4+g(\xi+2)}=2\cw.
\]
The difference of the two square roots is $g/\cw$. Hence the smaller square root is $\cw-g/(2\cw)$, and
\[
g\xi=\left(\cw-\frac{g}{2\cw}\right)^2-\frac{\omega^2}{4}
=\frac{3\cw^2}{4}-\frac g2.
\]
Condition \eqref{isola:range} makes this wave number positive and verifies the signs of both square roots. Substitution gives
\begin{align*}
    s_*=\cw\xi_*-\frac \omega2+\sqrt{\frac{\omega^2}{4}+g\xi_*}=\frac{3\cw^3}{4 g}
\end{align*}
and, hence, proves \eqref{isola:anchor}--\eqref{isola:collision}.

For a fixed real frequency $s$, recalling \eqref{isola:flat-block}, the positive-wave-number equation is
\begin{equation}\label{isola:dispersion-polynomial}
s=\cw\xi-\frac\omega2
\pm\sqrt{\frac{\omega^2}{4}+g\xi},\qquad\text{equivalently}\qquad (s-\cw\xi)^2+\omega(s-\cw\xi)-g\xi=0.
\end{equation}
At $s=s_*$ its two roots are exactly $\xi_*$ and $\xi_*+2$. There are no further modes at this frequency. To exclude negative wave numbers, suppose that $\xi=-\beta<0$, with $\beta>0$, also produces the frequency $s_*$. Then
\[
s_*=-\cw\beta+\frac{\omega}{2}
\pm\sqrt{\frac{\omega^2}{4}+g\beta},
\]
and hence
\begin{align*}
0
&=(s_*+\cw\beta)^2-\omega(s_*+\cw\beta)-g\beta=\cw^2\beta^2+(2s_*\cw-\omega\cw-g)\beta+s_*(s_*-\omega).
\end{align*}
However, $\beta>0$, $\xi_*>0$, $s_*>0$, and
\[
2s_*\cw-\omega\cw-g=2\cw s_*-\cw^2=\cw^2\left(\frac{3\cw^2}{2g}-1\right)>0,\qquad
s_*-\omega
=\frac{3\cw^4-4g\cw^2+4g^2}{4g\cw}
=\frac{(\cw^2-2g)^2+2\cw^4}{4g\cw}>0.
\]
Thus every term in the preceding sum is positive, a contradiction. Consequently, no negative wave number produces the frequency $s_*$. Since $\bar\mu\ne0$, no Fourier mode has zero physical wave number. Therefore, the only modes producing $\im s_*$ are those with physical wave numbers $\xi_*$ and $\xi_*+2$.

Each of the two resonant Fourier blocks has two distinct eigenvalues, so $\im s_*$ is simple within each block. It follows that $\im s_*$ is semisimple and has algebraic multiplicity two.

Finally, all Fourier-block eigenvalues are purely imaginary. Hence every Fourier block of $\mathcal L_{\bar\mu,0}-\mathrm{Id}$ is invertible, and its inverse is $\cO((1+|k|)^{-1})$ as $|k|\to\infty$. Consequently, $1$ belongs to the resolvent set and $(\mathcal L_{\bar\mu,0}-\mathrm{Id})^{-1}$ maps $\mathcal X$ boundedly into $\mathcal Y$, as defined in \eqref{def:X-Y-section3}. The compact embedding $\mathcal Y\hookrightarrow\mathcal X$ then implies that $\mathcal L_{\bar\mu,0}$ has compact resolvent. In particular, $\im s_*$ is an isolated eigenvalue.
\end{proof}

\begin{remark}\label{rem:isola-range}
The restriction \eqref{isola:range} comes from the location of the collision. At $\omega=-\sqrt{g/6}$ the smaller wave number reaches zero, where the multiplier $\operatorname{sign}\xi$ requires separate treatment. For more negative vorticity this positive-wave-number collision has no solution. At zero vorticity, \eqref{isola:anchor} gives $\xi_*=1/4$ and $s_*=3\sqrt g/4$, recovering the collision studied in \cite{BMV_isola}.
\end{remark}

\paragraph{Kato expansions} Assume \eqref{isola:range} and $\bar\mu\ne0$. Recalling \eqref{isola:anchor}, we set
\begin{equation}\label{isola:reference-parameters}
k_1=j_*,\qquad k_2=j_*+2,\qquad
\xi_1=\xi_*,\qquad \xi_2=\xi_*+2.
\end{equation}
The two eigenvectors of $\mathcal L_{\bar\mu,0}$ at $\im s_*$ are
\begin{equation}\label{isola:flat-right-basis}
f_1(x)=
\begin{bmatrix}1\\ \dfrac{2\im\cw g}{3\cw^2-2g}\end{bmatrix}e^{\im j_*x},
\qquad
f_2(x)=
\begin{bmatrix}1\\-\dfrac{2\im\cw g}{\cw^2+2g}\end{bmatrix}e^{\im(j_*+2)x}.
\end{equation}
Their adjoint eigenvectors are
\begin{equation}\label{isola:flat-left-basis}
\begin{aligned}
g_1(x)&=
\begin{bmatrix}
\dfrac{3\cw^2-2g}{2(2\cw^2-g)}\\[1mm]
\dfrac{\im\cw(3\cw^2-2g)}{4g(2\cw^2-g)}
\end{bmatrix}e^{\im j_*x},\qquad
g_2(x)&=
\begin{bmatrix}
\dfrac{\cw^2+2g}{2(2\cw^2+g)}\\[1mm]
-\dfrac{3\im\cw(\cw^2+2g)}{4g(2\cw^2+g)}
\end{bmatrix}e^{\im(j_*+2)x}.
\end{aligned}
\end{equation}
With the scalar product \eqref{def:scalar-product-section3}, they satisfy
\begin{equation}\label{isola:flat-biorthogonality}
\mathcal L_{\bar\mu,0}f_a=\im s_*f_a,\qquad
\mathcal L_{\bar\mu,0}^*g_a=-\im s_*g_a,\qquad
(f_b,g_a)_{\mathcal X}=\delta_{ab}.
\end{equation}

Write $\mu=\bar\mu+\delta$, with $|\delta|$ small enough that this local interval contains no integer. The multipliers $\operatorname{sign}(k+\bar\mu+\delta)$ are then fixed, so $\mathcal L_{\bar\mu+\delta,\e}$ is analytic with the common domain $\mathcal Y$. Let $\Gamma$ be a small positively oriented circle around $\im s_*$ containing no other reference eigenvalue, and set
\begin{equation}\label{isola:projection}
P(\delta,\e):=-\frac1{2\pi\im}\oint_\Gamma
(\mathcal L_{\bar\mu+\delta,\e}-\lambda)^{-1}\,d\lambda.
\end{equation}
Write $P_0:=P(0,0)$. Use the Kato construction of Lemma \ref{lem:Riesz-Kato-positive} for this projector, and let $U(\delta,\e)$ be the Kato map \eqref{def:Kato-U} with these $P$ and $P_0$. The transported right and dual bases are
\begin{equation}\label{isola:transported-bases}
f_a(\delta,\e):=U(\delta,\e)f_a,
\qquad
g_a(\delta,\e):=(U(\delta,\e)^{-1})^*g_a.
\end{equation}
At $(\delta,\e)=(0,0)$, these bases satisfy $f_a(0,0)=f_a$ and $g_a(0,0)=g_a$. Moreover, $(f_b(\delta,\e),g_a(\delta,\e))_{\mathcal X}=\delta_{ab}$, and $f_1(\delta,\e),f_2(\delta,\e)$ span the two-dimensional spectral subspace $\operatorname{Ran}P(\delta,\e)$.

\begin{lemma}[The Kato basis expansions]\label{lem:isola-kato-expansions}
For $a=1,2$,
\begin{equation}\label{isola:basis-expansions}
\begin{aligned}
f_a(\delta,\e)=f_a+\delta\dot f_a+\e f_a'
+\cO_{\mathcal Y}((|\delta|+|\e|)^2),\qquad
g_a(\delta,\e)=g_a+\delta\dot g_a+\e g_a'
+\cO_{\mathcal X}((|\delta|+|\e|)^2).
\end{aligned}
\end{equation}
where
\begin{equation}\label{app:isola-right-delta}
\begin{aligned}
\dot f_1(x)&=
\frac{\cw^4g}{(3\cw^2-2g)(2\cw^2-g)^2}
\begin{bmatrix}1\\-2\im g/\cw\end{bmatrix}e^{\im j_*x},\quad
\dot f_2(x)=
\frac{3\cw^4g}{(\cw^2+2g)(2\cw^2+g)^2}
\begin{bmatrix}1\\2\im g/(3\cw)\end{bmatrix}e^{\im(j_*+2)x},
\end{aligned}
\end{equation}
\begin{equation}\label{app:isola-left-delta}
\begin{aligned}
\dot g_1(x)&=
\frac{\cw^2(3\cw^2-2g)}{2(2\cw^2-g)^3}
\begin{bmatrix}-g\\\im(3\cw^2-2g)/(2\cw)\end{bmatrix}e^{\im j_*x},\quad
\dot g_2(x)=
\frac{\cw^2(\cw^2+2g)}{2(2\cw^2+g)^3}
\begin{bmatrix}-g\\-\im(\cw^2+2g)/(2\cw)\end{bmatrix}e^{\im(j_*+2)x},
\end{aligned}
\end{equation}
\begin{equation}\label{app:isola-right-eps}
\begin{aligned}
f_1'(x)
={}&\frac{\cw^2+g}{2\cw^3}
\begin{bmatrix}
(\cw+\omega)(\xi_1+1)\\[1mm]
\dfrac{\im g}{2}
\end{bmatrix}
e^{\im(j_*+1)x}+\begin{cases}
\displaystyle
\frac{\cw^2+g}{3\cw^3}
\begin{bmatrix}
\omega(\xi_1-1)\\[1mm]
-\dfrac{3\im}{2}(2\cw^2-g)
\end{bmatrix}
e^{\im(j_*-1)x},
&\displaystyle\frac{2g}{3}<\cw^2<2g,
\\[5mm]
\displaystyle
-\frac{\cw^2+g}{6\cw^3}
\begin{bmatrix}
(3\cw+\omega)(\xi_1-1)\\[1mm]
\dfrac{3\im}{2}(2\cw^2-g)
\end{bmatrix}
e^{\im(j_*-1)x},
&\cw^2>2g,
\end{cases}
\\[4mm]
f_2'(x)
={}&\frac{\cw^2+g}{2\cw^3}
\begin{bmatrix}
(\omega-\cw)(\xi_2-1)\\[1mm]
\dfrac{\im}{2}(2\cw^2+g)
\end{bmatrix}
e^{\im(j_*+1)x}-\frac{\cw^2+g}{6\cw^3}
\begin{bmatrix}
(\omega-3\cw)(\xi_2+1)\\[1mm]
\dfrac{\im}{2}(4\cw^2+5g)
\end{bmatrix}
e^{\im(j_*+3)x}.
\end{aligned}
\end{equation}
\end{lemma}

The calculation is given in Appendix \ref{sec:AppIsolaKato}. We pause to remark that we don't need the closed forms of $g'_1$ and $g'_2$ for the matrix representation, so we omit the formulas here.

\paragraph{The two-dimensional matrix} In the basis $(f_1(\delta,\e),f_2(\delta,\e))$, the matrix representation of
$\mathcal L_{\bar\mu+\delta,\e}|_{\operatorname{Ran}P(\delta,\e)}$ is
\begin{equation}\label{isola:H-definition}
\begin{gathered}
\mathbf H(\bar\mu+\delta,\e)
:=\big((\mathcal L_{\bar\mu+\delta,\e}f_b(\delta,\e),g_a(\delta,\e))_{\mathcal X}\big)_{a,b\in\{1,2\}}
=\im\begin{pmatrix}h_{11}&h_{12}\\h_{21}&h_{22}\end{pmatrix}.
\end{gathered}
\end{equation}
Thus
\begin{equation}\label{isola:H-action}
\mathcal L_{\bar\mu+\delta,\e}f_b(\delta,\e)
=\im\sum_{a=1}^2h_{ab}f_a(\delta,\e).
\end{equation}
The numbers $h_{ab}$ are real for real parameters, by the reversibility argument of Lemma \ref{lem:reversibility-section5}.

\begin{lemma}\label{lem:isola-coefficients} The matrix $\mathbf H$ is analytic in $(\delta,\e)$ and its expansion at $(\delta,\e)=(0,0)$ is
\begin{equation}\label{isola:quadratic-matrix}
\begin{aligned}
\mathbf H(\bar\mu+\delta,\e)
={}&\im s_*I_2+\im\delta\begin{pmatrix}
\dfrac{2\cw^3}{2\cw^2-g}&0\\[2mm]
0&\dfrac{2\cw^3}{2\cw^2+g}
\end{pmatrix}+\im\e^2\begin{pmatrix}C_{11}&C_{12}\\C_{21}&C_{22}\end{pmatrix}
+\im\cO(\delta^2+|\delta|\e^2+\e^3),
\end{aligned}
\end{equation}
where
\begin{align}
C_{12}=-\frac{3(\cw^2-g)^3(\cw^2+g)(3\cw^2-2g)}
{64\cw^3g^2(2\cw^2-g)},\qquad C_{21}=\frac{9(\cw^2-g)^3(\cw^2+g)(\cw^2+2g)}
{64\cw^3g^2(2\cw^2+g)},\label{isola:C21}
\end{align}
\begin{equation}\label{isola:C22}
C_{22}=
\frac{(\cw^2+g)(\cw^2+2g)
(3\cw^8-11g\cw^6-4g^2\cw^4-7g^3\cw^2+g^4)}
{32\cw^5g^2(2\cw^2+g)},
\end{equation}
and
\begin{equation}\label{isola:C11}
C_{11}=
\begin{cases}
-\displaystyle\frac{(\cw^2+g)(3\cw^2-2g)
(14\cw^6-18g\cw^4-g^2\cw^2+3g^3)}
{64\cw^3g^2(2\cw^2-g)},&\displaystyle\frac{2g}{3}<\cw^2<2g,\\[3mm]
-\displaystyle\frac{(\cw^2+g)(3\cw^2-2g)
(\cw^8+3g\cw^6+4g^2\cw^4-9g^3\cw^2+3g^4)}
{32\cw^5g^2(2\cw^2-g)},&\cw^2>2g.
\end{cases}
\end{equation}
In particular, if \eqref{isola:range} holds and $\omega\ne0$, then
\begin{equation}\label{isola:negative-product}
\mathrm{ind}_{\mathrm{high}}:=-C_{12}C_{21}=
\frac{27(\cw^2-g)^6(\cw^2+g)^2(\cw^2+2g)(3\cw^2-2g)}
{4096\cw^6g^4(2\cw^2-g)(2\cw^2+g)}>0,
\end{equation}
or equivalently
\begin{align*}
    \mathrm{ind}_{\mathrm{high}}=\frac{27\omega^6(\cw^2+g)^2}{4096 g^4}
\frac{(\cw^2+2g)(3\cw^2-2g)}{(2\cw^2-g)(2\cw^2+g)}.
\end{align*}
\end{lemma}
\begin{proof}
The derivation of Lemma \ref{lem:isola-coefficients} is given in Appendix \ref{sec:AppIsolaKato}
\end{proof}

\paragraph{The unstable spectral curve}
\label{subsec:isola-curve}
The eigenvalues of $\mathbf H$ are obtained directly from
\[
\lambda^2-\im(h_{11}+h_{22})\lambda
-(h_{11}h_{22}-h_{12}h_{21})=0.
\]
The quadratic formula gives
\begin{equation}\label{isola:exact-eigenvalues}
\begin{aligned}
\lambda^{(2)}_\pm(\mu,\e)
&=\frac{\im}{2}(h_{11}+h_{22})
\pm\sqrt{\Delta_2(\mu,\e)},\qquad\Delta_2(\mu,\e):=-\frac14(h_{11}-h_{22})^2-h_{12}h_{21}.
\end{aligned}
\end{equation}
Thus $\Delta_2>0$ gives opposite nonzero real parts, while $\Delta_2<0$ gives two purely imaginary eigenvalues. Assume
\begin{equation}\label{isola:theorem-assumptions}
\omega>-\sqrt{g/6},\qquad \omega\ne0,\qquad \xi_*\notin\mathbb Z.
\end{equation}
Inspired by \cite{NRS2023,JRSY2025}, we first solve $h_{11}(\mu,\e)=h_{22}(\mu,\e)$. Since $\partial_\mu(h_{11}-h_{22})(\bar\mu,0)=\frac{4g\cw^3}{4\cw^4-g^2}>0$, the implicit-function theorem gives a real analytic function
\begin{equation}\label{isola:mu-center}
\mu_0(\e)=\bar\mu+\mu_2\e^2+\cO(\e^3),\qquad
\mu_2=-\frac{(4\cw^4-g^2)(C_{11}-C_{22})}{4g\cw^3}.
\end{equation}
Also, recalling \eqref{isola:quadratic-matrix}, \eqref{isola:negative-product} and using \eqref{isola:mu-center}, we notice that
\begin{align}
    -(h_{12}h_{21})(\mu_0(\e),\e)=\e^4\left(\mathrm{ind}_{\mathrm{high}}+\cO(|\delta|+|\e|)\right)\overset{\eqref{isola:mu-center}}{=}\e^4\left(\mathrm{ind}_{\mathrm{high}}+\cO(|\e|)\right).
\end{align}
For the reminder, we apply Lemma \ref{final lemma} so that we can pull out the factor $\e^4$. Therefore, we obtain
\begin{align} \label{D2}
    \Delta_2(\mu_0(\e),\e)=\e^4\left(\mathrm{ind}_{\mathrm{high}}+\cO(|\e|)\right).
\end{align}

Next, the endpoints satisfy $\Delta_2=0$, or equivalently
\[
h_{11}(\mu,\e)-h_{22}(\mu,\e)
=\pm2\sqrt{-h_{12}(\mu,\e)h_{21}(\mu,\e)}.
\]
Here $\sqrt{-h_{12}h_{21}}=\e^2(\sqrt{\mathrm{ind}_{\mathrm{high}}}+\cO(|\mu-\bar\mu|+\e))$ is analytic in $(\mu,\e)$, since Lemma \ref{final lemma}. Moving the right-hand side to the left, both equations have $\mu$-derivative $\frac{4g\cw^3}{4\cw^4-g^2}$ at $(\bar\mu,0)$. The implicit-function theorem therefore gives two real analytic functions
\begin{equation}\label{isola:endpoints}
\mu_\pm(\e)=\mu_0(\e)
\pm\frac{\sqrt{\mathrm{ind}_{\mathrm{high}}}(4\cw^4-g^2)}{2g\cw^3}\e^2+\cO(\e^3).
\end{equation}
We now substitute the known coefficients into the quadratic formula on the interval between these endpoints.

\begin{theorem}[The first isola]\label{thm:first-isola}
Assume \eqref{isola:theorem-assumptions}. For sufficiently small $0<|\e|<\e_0$, the two eigenvalues near $\im s_*$ have opposite nonzero real parts precisely when
\begin{equation}\label{isola:interval}
\mu_-(\e)<\bar\mu+\delta<\mu_+(\e).
\end{equation}
They collide at the endpoints and are purely imaginary outside this interval in a fixed neighborhood of $\bar\mu$. On the closed interval,
\begin{equation}\label{isola:quadratic-eigenvalues}
\begin{gathered}
\lambda^{(2)}_\pm(\bar\mu+\delta,\e)
=\im\left(s_*+\frac{4\cw^5}{4\cw^4-g^2}\delta
+\frac{C_{11}+C_{22}}2\e^2+\cO(|\e|^3)\right)\\
\qquad\pm\sqrt{
\mathrm{ind}_{\mathrm{high}}\e^4-
\left(\dfrac{2g\cw^3}{4\cw^4-g^2}\delta+\frac{C_{11}-C_{22}}{2}\e^2\right)^2
+\cO(|\e|^5)}.
\end{gathered}
\end{equation}
\end{theorem}

\begin{proof}
The two endpoint equations above are strictly increasing in $\mu$. The factorization
\[
\Delta_2=-\frac{1}{4}\big(h_{11}-h_{22}-2\sqrt{-h_{12}h_{21}}\big)
\big(h_{11}-h_{22}+2\sqrt{-h_{12}h_{21}}\big)
\]
therefore makes $\Delta_2$ positive between $\mu_-$ and $\mu_+$, and negative outside. This proves the instability assertion by \eqref{isola:exact-eigenvalues}. For $\mu_-(\e)\le\bar\mu+\delta\le\mu_+(\e)$, \eqref{isola:mu-center}--\eqref{isola:endpoints} give $\delta=\cO(\e^2)$. Substituting \eqref{isola:quadratic-matrix} into the quadratic formula gives \eqref{isola:quadratic-eigenvalues}.
\end{proof}

\begin{remark}[The order of the expansion]\label{rem:isola-order}
For each fixed vorticity satisfying \eqref{isola:theorem-assumptions}, the second-order coefficients determine the formulas above. The estimate is not uniform as $\omega\to0$, since $\sqrt{\mathrm{ind}_{\mathrm{high}}}=\cO(|\omega|^3)$. At $\omega=0$, the off-diagonal coefficients of $\e^2$ vanish and the isola opens at fourth order \cite{BMV_isola}. Our diagonal coefficients at this value give $\mu_2=-57/64$ consistent with that calculation when $g=1$.
\end{remark}

\begin{remark} [Krein Signature for $\mathrm{index}_{\mathrm{high}}$ in \eqref{D2}]
Suppose $\mu\notin \mathbb Z$. Recalling \eqref{isola:operator}, by Wahlén \cite{Wahlen07}, we have
\begin{equation} \label{mathcal L mu e}
\begin{aligned}
    \mathcal{L}_{\mu,\e}:=&\begin{bmatrix}
(\partial_x+\im\mu)\circ(\cw+p_\e)&|D+\mu|\\
-(g+a_\e)&(\cw+p_\e)(\partial_x+\im\mu)+\omega\mathfrak H_\mu
\end{bmatrix}\\
    =&\underbrace{\begin{bmatrix}
        0& \mathrm{Id}\\
        -\mathrm{Id}& \omega (\pa_x+\im \mu)^{-1}
    \end{bmatrix}}_{=:\mathcal{J}_\mu}\underbrace{\begin{bmatrix}
 g+a_\epsilon+\omega(\cw+p_\e)& ~~-(\cw+p_\epsilon)(\partial_x+\im\mu)\\
 (\partial_x+\im\mu)\circ(\cw+p_\epsilon)& ~~|D+\mu| 
    \end{bmatrix}}_{=:\mathcal{B}_{\mu,\e}},
\end{aligned}    
\end{equation}
where
\begin{align}
    \mathcal{J}^*_{\mu}=-\mathcal{J}_{\mu},\qquad \mathcal{B}^*_{\mu,\e}=\mathcal{B}_{\mu,\e},\qquad (\mathcal{J}_{\mu})^{-1}=\begin{bmatrix}
        \omega (\pa_x+\im \mu)^{-1}& -\mathrm{Id}\\
        \mathrm{Id}& 0
    \end{bmatrix}.
\end{align}
We notice that
\begin{align} \label{LJ JL}
    (\mathcal{J}_{\mu})^{-1}\mathcal{L}_{\mu,\e}+\mathcal{L}^*_{\mu,\e}(\mathcal{J}_\mu)^{-1}=0.
\end{align}
Set $\mu=\bar\mu+\delta$. Recalling \eqref{isola:transported-bases}, we denote $\mathfrak F(\delta,\e):=\begin{bmatrix}
        f_1(\delta,\e) & f_2(\delta,\e)
    \end{bmatrix}$. Also, the matrix representation \eqref{isola:H-definition} means
    \begin{align} \label{LF FH}
        \mathcal{L}_{\mu,\e}\mathfrak F(\delta,\e)=\mathfrak F(\delta,\e) \mathbf{H}(\mu,\e).
    \end{align}
Multiplies $\mathfrak F^*$ on the left hand side of \eqref{LJ JL} and $\mathfrak F$ on the right hand side of \eqref{LJ JL} and combining \eqref{LF FH}, we may obtain
\begin{align} \label{KH HK}
    K \mathbf{H}+\mathbf{H}^* K=0, \qquad K:=\im\mathfrak F^*(\mathcal{J}_\mu)^{-1}\mathfrak F.
\end{align}
The element in $(1,2)$ of \eqref{KH HK} satisfies
\begin{align}
    K_{11} h_{12}-K_{22} h_{21}=K_{12}(h_{11}-h_{22}).
\end{align}
In particular, at $(\mu,\e)=(\mu_0(\e),\e)$, we have
\begin{align}
    K_{11} h_{12}=K_{22} h_{21}.
\end{align}
Therefore, recalling \eqref{isola:exact-eigenvalues} and \eqref{D2}, we arrive at
\begin{align}
    \Delta_2(\mu_0(\e),\e)=-\frac{K_{11}(\mu_0(\e),\e)}{K_{22}(\mu_0(\e),\e)}|h_{12}(\mu_0(\e),\e)|^2.
\end{align}
Also, recalling \eqref{KH HK}, we have
\begin{align}
    K(\bar\mu,0)=\im\mathfrak F^*(0,0)(\mathcal{J}_{\bar{\mu}})^{-1}\mathfrak F(0,0),
\end{align}
where, by \eqref{isola:flat-right-basis},
\begin{align}
    \mathfrak F(0,0)=\begin{bmatrix}
        f_1& f_2
    \end{bmatrix}=
        \begin{bmatrix} e^{\im j_*x} & e^{\im(j_*+2)x}\\ \dfrac{2\im\cw g}{3\cw^2-2g} e^{\im j_*x} & -\dfrac{2\im\cw g}{\cw^2+2g} e^{\im(j_*+2)x}
    \end{bmatrix}.
\end{align}
A straightforward calculation reveals that
\begin{align}
    K_{11}(\bar\mu,0)=\frac{4g(2\cw^2-g)}{\cw(3\cw^2-2g)}, \qquad K_{22}(\bar\mu,0)=-\frac{4g(2\cw^2+g)}{3\cw(\cw^2+2g)}, \qquad K_{12}(\bar\mu,0)=K_{21}(\bar\mu,0)=0.
\end{align}
This implies
\begin{align}
    \Delta_2(\mu_0(\e),\e)=\left(\frac{3(2\cw^2-g)(\cw^2+2g)}{(3\cw^2-2g)(2\cw^2+g)}+\cO(|\e|)\right)|h_{12}(\mu_0(\e),\e)|^2,
\end{align}
and hence, using the fact $h_{12}(\mu_0(\e),\e)=\e^2 (C_{12}+\cO(|\e|)$,
\begin{align}
    \Delta_2(\mu_0(\e),\e)=\e^4\left(\frac{3(2\cw^2-g)(\cw^2+2g)}{(3\cw^2-2g)(2\cw^2+g)}C^2_{12}+\cO(|\e|)\right)=\e^4 \left(\mathrm{ind}_{\mathrm{high}}+\cO(|\e|)\right).
\end{align}
\end{remark}

\paragraph{Integer resonances and instability for every vorticity}
\label{subsec:isola-all-vorticities}

It remains to treat integer resonances when $\cw^2\ge3g$. In this case $\bar\mu=0$ and $0<\mu_0(\e)<\tfrac12$, where the same quadratic formula for $\mathcal L_{\mu_0(\e),\e}$ gives instability.

\begin{proposition}[Integer resonances]\label{prop:isola-integer}
Suppose $\cw^2\ge3g$ and $\xi_*=m\in\mathbb Z$. Then, for every sufficiently small $\e>0$, there is a Bloch exponent $\mu(\e)\in(0,\tfrac12)$, with $\mu(\e)=\cO(\e^2)$, for which $\mathcal L_{\mu(\e),\e}$ has an eigenvalue satisfying
\begin{equation}\label{isola:integer-growth}
\operatorname{Re}\lambda=\sqrt{\mathrm{ind}_{\mathrm{high}}}\e^2+\cO(|\e|^3)>0.
\end{equation}
\end{proposition}

\begin{proof}
Here $m\ge2$. The one-sided limit of $\mathcal L_{\mu,0}$ as $\mu\to0^+$ has the same isolated eigenvalue $\im s_*$, of multiplicity two, at wave numbers $m$ and $m+2$; its zero block has eigenvalues $0,-\im\omega$, distinct from $\im s_*$. Applying the reduction above to the analytic extension from $\mu>0$ therefore gives the matrix of $\mathcal L_{\mu,\e}$ for $0<\mu<\tfrac12$. The coefficients \eqref{isola:C22} and the second line of \eqref{isola:C11} are unchanged, since $m-1,m+1,m+3>0$.

Solving $h_{11}=h_{22}$ gives $\mu_0(\e)=\mu_2\e^2+\cO(\e^3)$ with $\mu_2$ as in \eqref{isola:mu-center}. Substitution of the diagonal coefficients, with $x=\cw^2/g\ge3$, gives
\begin{align}\label{isola:positive-center}
\mu_2=\frac{x+1}{32x^4}\big[&
3(x-3)^6+55(x-3)^5+413(x-3)^4\notag\\
&+1609(x-3)^3+3386(x-3)^2+3584(x-3)+1444\big]>0.
\end{align}
Thus $0<\mu_0(\e)<\tfrac12$ for small $\e>0$, and \eqref{isola:exact-eigenvalues} gives the eigenvalues of $\mathcal L_{\mu_0(\e),\e}$. At this exponent, the quadratic formula gives
\[
\lambda^{(2)}_+(\mu_0(\e),\e)
=\frac{1}{2}\operatorname{tr}\mathbf H(\mu_0(\e),\e)
+\sqrt{-h_{12}h_{21}}
=\frac{1}{2}\operatorname{tr}\mathbf H(\mu_0(\e),\e)
+\sqrt{\mathrm{ind}_{\mathrm{high}}}\e^2+\cO(\e^3),
\]
which has positive real part.
\end{proof}

\begin{remark}\label{rem:isola-integer-geometry}
The unstable exponent in Proposition \ref{prop:isola-integer} is strictly positive. Its proof therefore uses the positive Bloch family and does not require the auxiliary periodic operator $\mathcal L^{\mathrm{per}}_{\e}$ in \eqref{auxiliary periodic operator}.
\end{remark}

Combining this result with the low-frequency instability gives the following conclusion for every constant vorticity.

\begin{theorem}[Spectral instability for every constant vorticity]
\label{thm:all-vorticities-unstable}
For every $g>0$ and every fixed $\omega\in\mathbb R$, there exists $\e_0=\e_0(g,\omega)>0$ such that each Stokes wave in the family \eqref{exp:Sto} with $0<\e<\e_0$ is spectrally unstable under longitudinal perturbations in $L^2(\mathbb R;\mathbb C^2)$. More precisely, $\mathcal L_{\mu,\e}$ has an eigenvalue with positive real part for some $\mu\ne0$.
\end{theorem}

\begin{proof}
If $\omega<2\sqrt{g/3}$ and $\omega\ne0$, Theorem \ref{Complete BF thm} supplies an unstable Benjamin--Feir pair for $\mathcal L_{\mu,\e}$ with sufficiently small $\mu>0$. At $\omega=0$, the same conclusion is the classical deep-water result \cite{NS,HY,BMV1}.

Suppose next that $\omega\ge2\sqrt{g/3}$, equivalently $\cw^2\ge3g$. Then \eqref{isola:range} holds and $\omega\ne0$. If $\xi_*\notin\mathbb Z$, Theorem \ref{thm:first-isola} provides an unstable eigenvalue of $\mathcal L_{\mu,\e}$ for some $\mu\ne0$. If $\xi_*\in\mathbb Z$, Proposition \ref{prop:isola-integer} provides one. This covers the parameter range. In each case, the unstable eigenvalue can be chosen in the interior of an instability interval, with a nonzero Bloch exponent away from the endpoints of the chosen coordinate chart. It is simple in its isolated spectral cluster and admits a local continuous eigenvector in $\mathcal Y$, normalized in $\mathcal X$. Lemma~\ref{lem:Bloch-to-real-line-spectrum} therefore gives a spectral point of $\mathcal A_\e$ with positive real part. The smallness threshold is chosen for the fixed pair $(g,\omega)$.
\end{proof}

\appendix
\section{Expansion of the Stokes waves with constant vorticity in deep water}\label{sec:App2}

In this Appendix we provide the expansions \eqref{exp:Sto}--\eqref{expc2}, \eqref{expansion frak p}, and prove Lemma \ref{lem:pa-expansions}.
\\[1mm]
\noindent
{\bf Proof  of \eqref{exp:Sto}-\eqref{expcoef}.} Writing
 \begin{equation}\label{etapsic}
\begin{aligned}
 & \eta_\e(x) = \e \eta_1(x) + \e^2 \eta_2(x)+\e^3 \eta_3(x) + \cO(\e^4) \, , \\
 &  \psi_\e(x) = \e \psi_1(x) + \e^2 \psi_2(x)+\e^3 \psi_3(x) + \cO(\e^4) \, ,  
  \end{aligned}
\qquad \quad c_\e = \cw + \e c_1 + \e^2 c_2+ \cO(\e^3) \, ,  
\end{equation}
where  $\eta_i$ is $even(x)$ and $\psi_i$ is $odd(x)$ for $i=1,2$, we solve order by order in $ \e $ the equations \eqref{travelingWWstokes}, that we rewrite  as
\begin{equation}
\label{Sts}
\begin{cases}
-c \, \psi_x +  g\eta   + \dfrac{\psi_x^2}{2} - 
\dfrac{\eta_x^2}{2(1+\eta_x^2)} ( c+\omega\eta  -  \psi_x )^2-\omega\eta\psi_x-\omega\pa^{-1}_x G(\eta)\psi  = 0, \\
c \, \eta_x+G(\eta)\psi+\omega\eta\eta_x = 0 	\, ,
\end{cases}
\end{equation}
having substituted $G(\eta)\psi $ with $-c \, \eta_x-\omega\eta\eta_x $  in the first equation. We expand the Dirichlet-Neumann operator $ G(\eta)=  G_0+ G_1(\eta) + G_2(\eta) + \cO(\eta^3)  $ where, according to \cite[formula (2.14)]{CS},
\begin{equation} \label{expDiriNeu}
\begin{aligned}
 G_0 & := |D| \,, \\
 G_1(\eta) & :=  -\pa_x \eta \pa_x - |D|\eta|D| , \\
 G_2(\eta) & := -\frac12 |D| \Big( 
 |D|{\eta}^2  +\eta^2 |D| - 2\eta |D|\eta \Big) |D|\, .
\end{aligned}
\end{equation}
{\bf First order in $ \e $.}  Substituting the expansions in \eqref{etapsic} into \eqref{Sts}, we get the linear system
\begin{equation}\label{cB0}
 \left\{\begin{matrix} -\cw (\psi_1)_x + g\eta_1 -\omega\pa^{-1}_xG_0\psi_1= 0 \\
 \cw (\eta_1)_x + G_0\psi_1 =0 \, ,  \end{matrix}\right.
 \quad  \text{i.e.} \, \vet{\eta_1}{\psi_1} \in \text{Ker }\cB_0  \text{ with } \cB_0 := \begin{bmatrix} g & -\cw\pa_x-\omega\pa^{-1}_xG_0 \\ \cw\pa_x & G_0  \end{bmatrix},  
\end{equation}
where $\eta_1$ is $\mathrm{even}(x)$ and $\psi_1$ is $\mathrm{odd}(x)$.
\begin{lemma}\label{lem:B0}
If $g>0$ and $\omega \in \R $,  the kernel of the linear operator $\cB_0$ in \eqref{cB0} is one dimensional and given by
\begin{equation}\label{chk}
\mathrm{Ker} \,\cB_0= \mathrm{span}\,\left\{\vet{\cos(x)}{\cw\sin(x)} \right\}.
\end{equation}
\end{lemma}

\begin{proof} The action of $\cB_0$ on each subspace spanned by $\footnotesize{\,\left\{\vet{\cos(kx)}{0}, \vet{0}{\sin(kx)}\right\}} $, $k\in \N$, is represented by the $2\times 2$ matrix $\footnotesize{ \begin{bmatrix} g & -\cw k+\omega \\ -\cw k & k  \end{bmatrix}}$. Its determinant is given by
 $$ gk  - \cw^2 k^2+\omega\cw k\stackrel{\eqref{exp:Sto}}{=} \cw^2k(1-k)$$
and, provided $g>0$ and $\omega \in \R$, vanishes  if and only if $k=1$. For $k=1$ we obtain the kernel of $\cB_0$ given in \eqref{chk}. For $k=0$ it has no kernel since $\psi_1(x)$ is odd.
\end{proof}
We set $ \eta_1(x) := \cos(x)$, $\psi_1(x) := \cw\sin(x) $ in agreement with \eqref{exp:Sto}.
\\[1mm]
{\bf Second order in $\e$.} Collecting the terms of order $\e^2$ in \eqref{Sts}, and using $\cw^2(\eta_1)_x^2=(G_0\psi_1)^2$, we obtain
\begin{equation}\label{syslin2}
\cB_0
\vet{\eta_2}{\psi_2}
=\vet{R^{[1]}_2}{R^{[2]}_2},
\end{equation}
where
\begin{align*}
R^{[1]}_2:={}&
c_1(\psi_1)_x
-\frac12(\psi_1)_x^2
+\frac12\cw^2(\eta_1)_x^2
+\omega\eta_1(\psi_1)_x
+\omega\partial_x^{-1}G_1(\eta_1)\psi_1,
\\
R^{[2]}_2:={}&
-c_1(\eta_1)_x
-G_1(\eta_1)\psi_1
-\omega\eta_1(\eta_1)_x.
\end{align*}

Notice that the operator $\cB_0$ in \eqref{cB0} is not self-adjoint. Its adjoint kernel is given by
\begin{equation}\label{ker:B0adj}
    \mathrm{Ker}\,\cB_0^*
    =
    \mathrm{span}
    \left\{
    \vet{\cw\cos (x)}{g\sin (x)}
    \right\}.
\end{equation}
Hence, system \eqref{syslin2} admits a solution if and only if
\begin{equation}\label{orth1}
\left(
\vet{R^{[1]}_2}{R^{[2]}_2},\vet{\cw\cos(x)}{g\sin(x)}
\right)
=0.
\end{equation}

By the first-order expansion, $\eta_1(x)=\cos (x)$, $\psi_1(x)=\cw\sin (x)$. Moreover, using \eqref{expDiriNeu}, we have
\begin{align*}
G_1(\eta_1)\psi_1
&=
-\partial_x\big(\eta_1\partial_x\psi_1\big)
-|D|\big(\eta_1|D|\psi_1\big)\\
&=
-\partial_x\big(\cw\cos^2(x)\big)
-|D|\big(\cw\cos (x)\sin (x)\big)\\
&=
\cw\sin(2x)-\cw\sin(2x)=0.
\end{align*}
Therefore,
\begin{align*}
R^{[1]}_2
&=
c_1\cw\cos (x)
-\frac{\cw^2}{2}\cos^2(x)
+\frac{\cw^2}{2}\sin^2(x)
+\omega\cw\cos^2(x)\\
&=
c_1\cw\cos(x)
+\frac{\omega\cw}{2}
+\frac{\omega\cw-\cw^2}{2}\cos(2x)\\
&=
c_1\cw\cos (x)
+\frac{\omega\cw}{2}
-\frac{g}{2}\cos(2x),
\end{align*}
where we used the dispersion relation $\cw^2-\omega\cw-g=0$. Similarly,
\begin{align*}
  R^{[2]}_2
    =
    c_1\sin (x)+\frac{\omega}{2}\sin(2x).
\end{align*}
Substituting these expressions into \eqref{orth1}, we obtain $\frac12(\cw^2+g)c_1=0$, and hence $c_1=0$. Consequently, \eqref{syslin2} reduces to
\begin{equation}
\begin{bmatrix}\label{sisto2}
g&
-\cw\partial_x-\omega\partial_x^{-1}G_0\\
\cw\partial_x&
G_0
\end{bmatrix}
    \vet{\eta_2}{\psi_2}
=
\vet{\dfrac{\omega\cw}{2}-\dfrac{g}{2}\cos(2x)}{\dfrac{\omega}{2}\sin(2x)}.
\end{equation}

We set
\[
    \eta_2(x)
    =
    \eta_2^{[0]}+\eta_2^{[2]}\cos(2x),
    \qquad
    \psi_2(x)
    =
    \psi_2^{[2]}\sin(2x).
\]
Since
\[
    \partial_x^{-1}G_0\sin(2x)=-\cos(2x),
\]
system \eqref{sisto2} is equivalent to
\begin{equation*}
\left\{
\begin{aligned}
g\eta_2^{[0]}&=\frac{\omega\cw}{2},\\
g\eta_2^{[2]}+(\omega-2\cw)\psi_2^{[2]}&=-\frac{g}{2},\\
-2\cw\eta_2^{[2]}+2\psi_2^{[2]}&=\frac{\omega}{2}.
\end{aligned}
\right.
\end{equation*}
Solving the last two equations and using $\cw^2-\omega\cw-g=0$, gives $\eta_2^{[0]}$, $\eta_2^{[2]}$, and $\psi_2^{[2]}$ as in \eqref{expcoef}.
\\[1 mm]
{\bf Third order in $ \e $.} The third-order expansion of the Dirichlet--Neumann operator is
\[
\begin{aligned}
G(\eta_\e)\psi_\e
={}&
\e G_0\psi_1
+\e^2\big(
G_0\psi_2+G_1(\eta_1)\psi_1
\big)\\
&+\e^3\big(
G_0\psi_3
+G_1(\eta_1)\psi_2
+G_1(\eta_2)\psi_1
+G_2(\eta_1)\psi_1
\big)
+\cO(\e^4).
\end{aligned}
\]
Set
\[
\mathcal N_3
:=
G_1(\eta_1)\psi_2
+G_1(\eta_2)\psi_1
+G_2(\eta_1)\psi_1.
\]
Collecting the terms of order $\e^3$ in \eqref{Sts}, we obtain
\begin{equation}\label{syslin3}
\cB_0
\vet{\eta_3}{\psi_3}
=
\vet{R_3^{[1]}}{R_3^{[2]}},
\end{equation}
where
\begin{equation*} 
\begin{aligned}
R_3^{[1]}
={}&
c_2(\psi_1)_x
-(\psi_1)_x(\psi_2)_x
+\cw^2(\eta_1)_x(\eta_2)_x
+\cw(\eta_1)_x^2
\big(\omega\eta_1-(\psi_1)_x\big)
\notag\\
&\quad
+\omega\eta_1(\psi_2)_x
+\omega\eta_2(\psi_1)_x
+\omega\partial_x^{-1}\mathcal N_3,
\\
R_3^{[2]}
={}&
-c_2(\eta_1)_x
-\mathcal N_3
-\omega\big(
\eta_1(\eta_2)_x+\eta_2(\eta_1)_x
\big).
\end{aligned}
\end{equation*}
A direct calculation using \eqref{expDiriNeu} gives
\begin{align*}
G_1(\eta_1)\psi_2=0,
\qquad
G_1(\eta_2)\psi_1=\cw \eta_2^{[2]}\sin (x),
\qquad
G_2(\eta_1)\psi_1=\frac{\cw}{4}\sin (x).
\end{align*}
Thus
\begin{equation}\label{def:N3}
\mathcal N_3
=
\cw\left(\eta_2^{[2]}+\frac14\right)\sin (x),
\qquad
\partial_x^{-1}\mathcal N_3
=
-\cw\left(\eta_2^{[2]}+\frac14\right)\cos (x).
\end{equation}

Denote by $\Pi_1$ the projection onto the first Fourier harmonics. From \eqref{def:N3}, we find
\begin{align*}
\Pi_1 R_3^{[1]}
=
(c_2\cw+\rho_1)\cos (x),
\qquad
\Pi_1 R_3^{[2]}
=
(c_2+\rho_2)\sin (x),
\end{align*}
where
\[
\rho_1
=\eta_2^{[0]}\cw\omega+\eta_2^{[2]}\cw^2-\frac12 \eta_2^{[2]}\cw\omega-\frac{\cw^2}{4}-\cw \psi_2^{[2]}+\omega \psi_2^{[2]},\quad\rho_2=\eta_2^{[0]}\omega
-\eta_2^{[2]}\cw
+\frac12 \eta_2^{[2]}\omega
-\frac{\cw}{4}.
\]
By \eqref{ker:B0adj}, the Fredholm solvability condition for \eqref{syslin3} is
\begin{align*}
\left(
\vet{R_3^{[1]}}{R_3^{[2]}},
\vet{\cw\cos (x)}{g\sin (x)}
\right)=0.
\end{align*}
Consequently,
\begin{equation}\label{solv:c2}
(\cw^2+g)c_2+\cw\rho_1+g\rho_2=0.
\end{equation}

Using \eqref{expcoef} and the dispersion relation $\cw^2-\omega\cw-g=0$, a direct simplification yields
\[
\cw\rho_1+g\rho_2
=
\frac{\cw^2+g}{8g\cw^5}
\left(
4\cw^8-11g\cw^6+7g^2\cw^4-3g^3\cw^2-g^4
\right).
\]
Substituting this identity into \eqref{solv:c2}, we obtain
\[
c_2
=
\frac{
-4\cw^8+11g\cw^6-7g^2\cw^4+3g^3\cw^2+g^4
}{
8g\cw^5
},
\]
which proves \eqref{expc2}.

\vspace{1mm}
\noindent
{\bf Expansion of the conformal map and proof of Lemma \ref{lem:pa-expansions}.} We first expand the odd function $\mathfrak p_\e$ defined by \eqref{def:ttf} as
\begin{equation}\label{app:frakp-expand}
\mathfrak p_\e(x)
=
\e\mathfrak p_1(x)+\e^2\mathfrak p_2(x)+\cO_{H^s}(\e^3).
\end{equation}
Using
\[
\eta_\e\big(x+\mathfrak p_\e(x)\big)
=
\e\eta_1
+
\e^2\big(\eta_2+(\eta_1)_x\mathfrak p_1\big)
+
\cO_{H^s}(\e^3)
\]
in \eqref{def:ttf}, and recalling that $\mathfrak H\cos(kx)=\sin(kx)$ for $k\in\N$ and $\mathfrak H[1]=0$, we obtain
\begin{align}
\mathfrak p_1(x)
&=
\mathfrak H\eta_1
=
\sin (x),
\label{app:frakp1}\\
\mathfrak p_2(x)
&=
\mathfrak H\big(\eta_2+(\eta_1)_x\mathfrak p_1\big)
=
\left(\eta_2^{[2]}+\frac12\right)\sin(2x)=
\frac{(\cw^2+g)^2}{4\cw^4}\sin(2x).
\label{app:frakp2}
\end{align}
This proves \eqref{expansion frak p}. In particular,
\begin{equation}\label{app:rho-expand}
\rho_\e(x)
=
1+\e\mathfrak p_1'(x)+\e^2\mathfrak p_2'(x)
+\cO_{H^s}(\e^3).
\end{equation}

In view of \eqref{exp:Sto}--\eqref{expc2}, the functions $B$ and $V$ defined in \eqref{espV}--\eqref{espB} expand as
\begin{align}
B
&=:
\e B_1(x)+\e^2B_2(x)+\cO_{H^s}(\e^3)
\notag\\
&=
\e\cw\sin (x)
+
\e^2
\frac{-\cw^4+g\cw^2+g^2}{2\cw^3}\sin(2x)
+
\cO_{H^s}(\e^3),
\label{app:B-expand}
\end{align}
and
\begin{align}
V
&=:
\e V_1(x)+\e^2V_2(x)+\cO_{H^s}(\e^3)
\notag\\
&=
\e\cw\cos (x)
+
\e^2\left(
\frac{\cw}{2}
+
\frac{-\cw^4+g\cw^2+g^2}{2\cw^3}\cos(2x)
\right)
+
\cO_{H^s}(\e^3).
\label{app:V-expand}
\end{align}

We next compute $p_\e$. By \eqref{def:alpha} and \eqref{def:qeps},
\begin{equation}\label{app:def-p}
\cw+p_\e(x)
=
\frac{
c_\e
-V\big(x+\mathfrak p_\e(x)\big)
+\omega\eta_\e\big(x+\mathfrak p_\e(x)\big)
}{1+\partial_x\mathfrak p_\e(x)}.
\end{equation}
Denoting derivatives with respect to $x$ by primes and suppressing the $x$-dependence when no confusion arises, we have
\begin{align*}
V(x+\mathfrak p_\e)
&=
\e V_1
+
\e^2\big(V_2+V_1'\mathfrak p_1\big)
+
\cO(\e^3),\\
\eta_\e(x+\mathfrak p_\e)
&=
\e\eta_1
+
\e^2\big(\eta_2+\eta_1'\mathfrak p_1\big)
+
\cO(\e^3),\\
\frac{1}{1+\partial_x\mathfrak p_\e}
&=
1-\e\mathfrak p_1'
+
\e^2\big((\mathfrak p_1')^2-\mathfrak p_2'\big)
+
\cO(\e^3).
\end{align*}
Therefore
\begin{align}
p_1
&=
-V_1+\omega\eta_1-\cw\mathfrak p_1',
\label{app:p1-formula}\\
p_2
&=
c_2-V_2-V_1'\mathfrak p_1
+\omega\big(\eta_2+\eta_1'\mathfrak p_1\big)
+(V_1-\omega\eta_1)\mathfrak p_1'-\cw\mathfrak p_2'
+\cw(\mathfrak p_1')^2.
\label{app:p2-formula}
\end{align}
Using \eqref{app:frakp1}--\eqref{app:V-expand}, together with $\omega=(\cw^2-g)/\cw$, we obtain
\begin{align}
p_1(x)
&=
-\frac{\cw^2+g}{\cw}\cos (x),
\label{app:p1-final}\\
p_2(x)
&=
\frac{(\cw^2+g)(3\cw^4+2g\cw^2+g^2)}{8\cw^5}+
\frac{(\cw^2+g)(\cw^4-4g\cw^2-g^2)}{4\cw^5}
\cos(2x).
\label{app:p2-final}
\end{align}
This proves \eqref{exp:p1}--\eqref{exp:p22}.

Finally, by \eqref{def:qeps}, we may write the exact identity
\begin{equation}\label{app:def-a}
g+a_\e(x)
=
\frac{g}{1+\partial_x\mathfrak p_\e(x)}
-
(\cw+p_\e(x))
B_x\big(x+\mathfrak p_\e(x)\big).
\end{equation}
Moreover,
\[
B_x(x+\mathfrak p_\e)
=
\e B_1'
+
\e^2\big(B_2'+B_1''\mathfrak p_1\big)
+
\cO(\e^3).
\]
Hence
\begin{align}
a_1
&=
-g\mathfrak p_1'-\cw B_1',
\label{app:a1-formula}\\
a_2
&=
g\big((\mathfrak p_1')^2-\mathfrak p_2'\big)
-\cw B_2'
-\cw B_1''\mathfrak p_1
-p_1B_1'.
\label{app:a2-formula}
\end{align}
Substituting \eqref{app:frakp1}--\eqref{app:frakp2}, \eqref{app:B-expand}, and \eqref{app:p1-final}, we find
\begin{align}
a_1(x)
&=
-(\cw^2+g)\cos (x),
\label{app:a1-final}\\
a_2(x)
&=
\cw^2+g
+
\frac{(\cw^2+g)(2\cw^4-3g\cw^2-g^2)}{2\cw^4}\cos(2x).
\label{app:a2-final}
\end{align}
This proves \eqref{exp:a1}--\eqref{exp:a22}, and hence Lemma \ref{lem:pa-expansions}.

\section{Generalized kernel of $\mathscr{L}_{0,\e}$} \label{G kernel}
\begin{proof}
Let $\mathbf F$ denote the stationary vector field in \eqref{travelingWWstokes} on $2\pi$-periodic functions, with $\partial_x^{-1}$ denoting the zero-average periodic primitive. The variations below are used only in the periodic space. Explicitly,
\[
\mathbf F(\eta,\psi,c):=
\begin{bmatrix}
c\eta_x+G(\eta)\psi+\omega\eta\eta_x\\
\begin{aligned}
&c\psi_x-g\eta-\dfrac{\psi_x^2}{2}+\dfrac{\big(G(\eta)\psi+\eta_x\psi_x\big)^2}
       {2(1+\eta_x^2)}
+\omega\eta\psi_x+\omega\partial_x^{-1}G(\eta)\psi
\end{aligned}
\end{bmatrix}.
\]
Thus $\mathbf F(\eta_\e,\psi_\e,c_\e)=0$. All derivatives of $\mathbf F$ below are evaluated at $(\eta_\e,\psi_\e,c_\e)$. By translation invariance, differentiation with respect to $x$ gives
\[
D_{(\eta,\psi)}\mathbf F
\begin{bmatrix}
\partial_x\eta_\e\\
\partial_x\psi_\e
\end{bmatrix}
=0.
\]
Since $D_c\mathbf F=(\partial_x\eta_\e,\partial_x\psi_\e)^\top$, differentiation with respect to $\e$ gives
\[
D_{(\eta,\psi)}\mathbf F
\begin{bmatrix}
\partial_\e\eta_\e\\
\partial_\e\psi_\e
\end{bmatrix}
=
-c_\e'
\begin{bmatrix}
\partial_x\eta_\e\\
\partial_x\psi_\e
\end{bmatrix}.
\]
For the vertical variation, $G(\eta_\e+a)=G(\eta_\e)$ in deep water, and the changes in $c$ and $\omega\eta$ cancel. Hence
\[
\mathbf F(\eta_\e+a,\psi_\e,c_\e-\omega a)
=
\begin{bmatrix}0\\-ga\end{bmatrix}.
\]
Differentiating at $a=0$ yields
\[
D_{(\eta,\psi)}\mathbf F
\begin{bmatrix}1\\0\end{bmatrix}
=
\omega
\begin{bmatrix}
\partial_x\eta_\e\\
\partial_x\psi_\e
\end{bmatrix}
-g\begin{bmatrix}0\\1\end{bmatrix}.
\]
Gauge invariance also gives
\[
D_{(\eta,\psi)}\mathbf F
\begin{bmatrix}0\\1\end{bmatrix}=0.
\]

At the fixed background wave, the periodic good-unknown and conformal changes of variables give
\[
\mathcal L^{\mathrm{per}}_{\e}\mathcal P Z^{-1}
=
\mathcal P Z^{-1}D_{(\eta,\psi)}\mathbf F.
\]
We therefore define the transformed variational vectors by
\begin{align*}
T_\e:=\mathcal P Z^{-1}
\begin{bmatrix}
\partial_x\eta_\e\\\partial_x\psi_\e
\end{bmatrix},\qquad A_\e:=\mathcal P Z^{-1}
\begin{bmatrix}
\partial_\e\eta_\e\\
\partial_\e\psi_\e
\end{bmatrix},\qquad W_\e:=\mathcal P Z^{-1}
\begin{bmatrix}1\\0\end{bmatrix},
\qquad G:=\begin{bmatrix}0\\1\end{bmatrix}.
\end{align*}
Since $\mathcal P Z^{-1}G=G$, the preceding identities become
\begin{align*}
\mathcal L^{\mathrm{per}}_{\e}T_\e=0,\qquad
\mathcal L^{\mathrm{per}}_{\e}G=0,\qquad
\mathcal L^{\mathrm{per}}_{\e}A_\e=-c_\e'T_\e,\qquad \mathcal L^{\mathrm{per}}_{\e}W_\e=\omega T_\e-gG.
\end{align*}

By \eqref{jump-jmath-section5},
\[
(\mathscr L_{0,\e}-\mathcal L^{\mathrm{per}}_{\e})h
=\mathfrak j_\e(h)G.
\]
Since $T_\e=((\eta_\e\circ\phi_\e)_x,V\circ\phi_\e)^\top$, its components are real odd and real even, respectively. The definition of $\mathfrak j_\e$ therefore gives
\[
\mathfrak j_\e(G)=-\im\omega,
\qquad
\mathfrak j_\e(T_\e)=-\im\omega m_\e.
\]
Consequently,
\begin{align*}
\mathscr L_{0,\e}G=-\im\omega G,\qquad \mathscr L_{0,\e}T_\e=-\im\omega m_\e G,\qquad
\mathscr L_{0,\e}A_\e
=-c_\e'T_\e+\mathfrak j_\e(A_\e)G,\qquad
\mathscr L_{0,\e}W_\e=\omega T_\e+\big(\mathfrak j_\e(W_\e)-g\big)G.
\end{align*}
Substituting the corrections in \eqref{adapted-zero-basis-section5} now yields
\begin{align*}
\mathscr L_{0,\e}E_2(\e)=0,\qquad
\mathscr L_{0,\e}E_1(\e)=-c_\e'(T_\e-m_\e G)
=-\im\e c_\e'E_2(\e),\qquad
\mathscr L_{0,\e}E_3(\e)=\omega(T_\e-m_\e G)
=\im\omega\e E_2(\e).
\end{align*}
This proves \eqref{exact-chain-section5}.
\end{proof}

\section{Expansion of the right and left bases}
\label{sec:AppKato}

In this appendix we compute the first-order expansions of $U_{\mu,\e}f_j$ and $(U_{\mu,\e}^{-1})^*g_j$, $j\in\{+,-,0\}$, defined in \eqref{def:right-Kato-basis}--\eqref{def:left-Kato-basis}. All parameter-dependent objects refer to the analytic family introduced in Section \ref{sec:Kato-positive}. All derivatives in this appendix are evaluated at $(\mu,\e)=(0,0)$. We write $U_0:=U_{0,0}=\mathrm{Id}$, and denote $\partial_\e$ by a prime and $\partial_\mu$ by a dot. The next lemma follows as \cite[Lemma A.1]{BMV1}

\begin{lemma} \label{U mu epsilon P00}
The first jets of $U_{\mu,\epsilon} P_{0}$ are
\begin{align} \label{A1}
    &U_{0} P_{0} = P_{0}, \quad\quad U'_{0} P_{0} = P'_{0} P_{0}, \quad\quad \dot{U}_{0} P_{0} = \dot{P}_{0} P_{0},
\end{align}
where
\begin{align} \label{A3}
    P'_{0} &= \frac{1}{2\pi \im} \oint_{\Gamma} (\mathscr{L}_{0} - \lambda)^{-1} \mathscr{L}'_{0} (\mathscr{L}_{0} - \lambda)^{-1} d\lambda\,,  \\ \label{A4}
    \dot{P}_{0} &= \frac{1}{2\pi \im} \oint_{\Gamma} (\mathscr{L}_{0} - \lambda)^{-1} \dot{\mathscr{L}}_{0} (\mathscr{L}_{0} - \lambda)^{-1} d\lambda\,.
\end{align}
The operators $\mathscr{L}'_{0}$ and $\dot{\mathscr{L}}_{0}$ are
\[
\mathscr L_0'
=
\begin{bmatrix}
\partial_x\circ p_1&0\\
-a_1&p_1\partial_x
\end{bmatrix},
\qquad
\dot{\mathscr L}_0
=
\begin{bmatrix}
0&\mathrm{sgn}(D)+\Pi_0\\
0&0
\end{bmatrix},
\]
where $p_1$ and $a_1$ are defined in Lemma \ref{lem:pa-expansions}.
\end{lemma}

By the Lemma \ref{U mu epsilon P00}, we have the Taylor expansion
\begin{equation} \label{the expandsion of f sigma mu}
    \begin{aligned}
    f_j (\mu, \epsilon) &= f_j + \epsilon P'_{0} f_j + \mu \dot{P}_{0} f_j  + \mathcal{O}((|\mu|+|\e|)^2),\qquad j\in\{+,-,0\}
\end{aligned}
\end{equation}
To compute the vectors $P'_{0} f_j$ and $\dot{P}_{0} f_j$ using \eqref{A3} and \eqref{A4}, it is useful to know the action of $(\mathscr{L}_{0} - \lambda)^{-1}$ on the vectors
\begin{equation} \label{A10}
    \begin{aligned}
        f^+_j:=e^{\im j x}\vet{1}{-\im \cw},~f^-_j:=e^{\im j x}\vet{1}{\im g/\cw},~f^+_{-j}:=e^{-\im j x}\vet{1}{\im \cw},~f^-_{-j}:=e^{-\im j x}\vet{1}{-\im g/\cw},~j\geq 1.
    \end{aligned}
\end{equation}
We pause to remark that
\begin{align}
    f^+_1=f_+,\qquad f^+_{-1}=f_-.
\end{align}

\begin{lemma} \label{Lemma B2}
The space $\mathcal Y=(H^1(\mathbb{T};\mathbb{C}))^2$ decomposes as $\mathcal Y = \mathcal{V}_{0} \oplus \mathcal{U} \oplus \mathcal{W}_{\mathcal Y}$, with
\begin{equation*}
    \mathcal{W}_{\mathcal Y} = \overline{\bigoplus_{j=2}^{\infty} \mathcal{W}_j}^{\mathcal Y}
\end{equation*}
where the subspaces $\mathcal{V}_{0}, \mathcal{U}$, and $\mathcal{W}_j$, defined below, are invariant under $\mathscr{L}_{0}$ and the following properties hold:

\begin{itemize}
\item[(i)] $\mathcal{V}_{0} = \text{span}\{f_+, f_-, f_0\}$ is the kernel of $\mathscr{L}_{0}$. For any $\lambda \neq 0$ the operator $\mathscr{L}_{0} - \lambda : \mathcal{V}_{0} \to \mathcal{V}_{0}$ is invertible and
    \begin{equation} \label{A11-0}
        (\mathscr{L}_{0} - \lambda)^{-1} f_j = -\frac{1}{\lambda} f_j, \quad j\in\{+,-,0\}.
    \end{equation}
\item[(ii)] $\mathcal{U} := \text{span}\{ f^-_1, f^-_{-1}, \vet{0}{1} \}$. For any $\lambda \notin \{\pm \im(\cw^2+g)/\cw, -\im \omega\}$ the operator $\mathscr{L}_{0} - \lambda : \mathcal{U} \to \mathcal{U}$ is invertible and
    \begin{equation} \label{A13}
    \begin{aligned}
        (\mathscr{L}_{0} - \lambda)^{-1} f^-_{\pm 1} = \frac{1}{\pm\im (\cw^2+g)/\cw-\lambda} f^-_{\pm 1},\quad (\mathscr{L}_{0} - \lambda)^{-1} \vet{0}{1}=\frac{1}{-\im \omega-\lambda}\vet{0}{1}.
    \end{aligned}
    \end{equation}

\item[(iii)] For $j\geq 2$, each subspace $\mathcal{W}_j := \text{span}\{ f_j^+, f_j^-, f_{-j}^+, f_{-j}^- \}$ is invariant under $\mathscr{L}_{0}$. Let
    \begin{equation*}
        \mathcal{W}_{\mathcal{X}} = \overline{\bigoplus_{j=2}^{\infty} \mathcal{W}_j}^{ \mathcal{X}}.
    \end{equation*}
For any $|\lambda| < \delta(\omega,g)$ small enough, the operator $\mathscr{L}_{0} - \lambda : \mathcal{W}_{\mathcal Y} \to \mathcal{W}_{\mathcal X}$ is invertible. Moreover, for every $f \in \mathcal{W}_{\mathcal X}$,
    \begin{equation} \label{A14}
        (\mathscr{L}_{0} - \lambda)^{-1} f = \left[\cw^2(\partial_x^2+|D|)\right]^{-1}
\begin{bmatrix}
\cw\partial_x+\omega\mathfrak H & -|D|\\
g & \cw\partial_x
\end{bmatrix}f
+\lambda\varphi_f(\lambda),
    \end{equation}
for some analytic function $\lambda \mapsto \varphi_f(\lambda, \cdot) \in \mathcal Y=(H^1(\mathbb{T};\mathbb{C}))^2$.
\end{itemize}
\end{lemma}

\begin{proof}
By inspection the spaces $\mathcal{V}_{0}, \mathcal{U}$ and $\mathcal{W}_j$ are invariant under $\mathscr{L}_{0}$ and, by Fourier series, they decompose $\mathcal Y=(H^1(\mathbb{T};\mathbb{C}))^2$. Formula \eqref{A11-0} follows using that $f_+, f_-, f_0$ are in the kernel of $\mathscr{L}_{0}$. Formula \eqref{A13} follows because $f_1^-$, $f_{-1}^-$, and $\vet{0}{1}$ are eigenvectors of $\mathscr L_0$ with respective eigenvalues $\im(\cw^2+g)/\cw$, $-\im(\cw^2+g)/\cw$, and $-\im\omega$. Let us prove item $(iii)$. Let $\mathcal{W} := \mathcal{W}_{\mathcal Y}$. The operator $(\mathscr{L}_{0} - \lambda \mathrm{Id})|_{\mathcal{W}}$ is invertible for any
$$\lambda \notin \left\{ \pm\im\left(j\cw-\frac{\omega}{2}\right)\pm\frac{\im}{2}\sqrt{\omega^2+4g j}, \quad j \geq 2, j \in \mathbb{N} \right\}$$
and
$$
(\mathscr{L}_{0}|_{\mathcal{W}})^{-1} = \left( \cw^2 (\partial^2_{x} + |D|) \right)^{-1} \begin{bmatrix} \cw \partial_x+\omega\mathfrak{H} & -|D|  \\ g & \cw \partial_x \end{bmatrix} \Big|_{\mathcal{W}}.
$$
By Neumann series, for any $\lambda$ such that
$$|\lambda| \| (\mathscr{L}_{0}|_{\mathcal{W}})^{-1} \|_{\mathcal{L}(\mathcal{W}, \mathcal Y)} < 1$$
we have
$$
(\mathscr{L}_{0}|_{\mathcal{W}} - \lambda)^{-1} = (\mathscr{L}_{0}|_{\mathcal{W}})^{-1} (\text{Id} - \lambda (\mathscr{L}_{0}|_{\mathcal{W}})^{-1})^{-1} 
= (\mathscr{L}_{0}|_{\mathcal{W}})^{-1} \sum_{j \geq 0} ((\mathscr{L}_{0}|_{\mathcal{W}})^{-1} \lambda)^j.
$$
Formula \eqref{A14} follows with
$$\varphi_f(\lambda, x) := (\mathscr{L}_{0}|_{\mathcal{W}})^{-1} \sum_{j \geq 1} \lambda^{j-1} [(\mathscr{L}_{0}|_{\mathcal{W}})^{-1}]^j f.$$ 
\end{proof}

To expand the right basis, recalling $f_+$, $f_-$, and $f_0$ in Proposition \ref{prop:reference-spectrum}, one may first observe that
\begin{align}
\mathscr L_0'f_+
&=
-\frac{\im(\cw^2+g)}{\cw}
e^{2\im x}
\begin{bmatrix}1\\0\end{bmatrix},
\qquad\mathscr L_0'f_-=
\frac{\im(\cw^2+g)}{\cw}
e^{-2\im x}
\begin{bmatrix}1\\0\end{bmatrix},\qquad\mathscr L_0'f_0
=
\begin{bmatrix}
\dfrac{\cw^2+g}{\cw}\sin(x)\\
(\cw^2+g)\cos(x)
\end{bmatrix},
\label{app:Lprime-f0}\\
\dot{\mathscr L}_0f_+
&=
e^{\im x}\begin{bmatrix}-\im\cw\\0\end{bmatrix},\qquad\dot{\mathscr L}_0f_-=
e^{-\im x}\begin{bmatrix}-\im\cw\\0\end{bmatrix},\qquad\dot{\mathscr L}_0f_0=
\begin{bmatrix}\dfrac{\im g}{\omega}\\0\end{bmatrix}.
\label{app:Ldot-f}
\end{align}
We are ready to calculate $P'_{0} f_j$ and $\dot{P}_{0} f_j$.

\begin{lemma}
\label{lem:app-right-Kato-jets}
One has
\begin{align}
P'_0 f_+&=
e^{2\im x}
\begin{bmatrix}
\dfrac{(\cw^2+g)^2}{2\cw^4}\\[2mm]
-\dfrac{\im g(\cw^2+g)}{2\cw^3}
\end{bmatrix},\qquad P'_0 f_-=
e^{-2\im x}
\begin{bmatrix}
\dfrac{(\cw^2+g)^2}{2\cw^4}\\[2mm]
\dfrac{\im g(\cw^2+g)}{2\cw^3}
\end{bmatrix},\qquad
P'_0 f_0=
\begin{bmatrix}
\dfrac{2\cw^2}{\cw^2+g}\cos(x)\\[2mm]
-\dfrac{2\cw g}{\cw^2+g}\sin(x)
\end{bmatrix},
\label{app:right-eps-zero}\\
\dot{P}_{0} f_+
&=
e^{\im x}
\begin{bmatrix}
\dfrac{\cw^4}{(\cw^2+g)^2}\\[2mm]
\dfrac{\im g\cw^3}{(\cw^2+g)^2}
\end{bmatrix},\qquad \dot{P}_{0} f_-=
e^{-\im x}
\begin{bmatrix}
-\dfrac{\cw^4}{(\cw^2+g)^2}\\[2mm]
\dfrac{\im g\cw^3}{(\cw^2+g)^2}
\end{bmatrix},
\qquad \dot{P}_{0} f_0=
\begin{bmatrix}
0\\[1mm]
-\dfrac{\im g^2}{\omega^3}
\end{bmatrix}.
\label{app:right-mu-zero}
\end{align}
\end{lemma}

\begin{proof}
We first calculate $P'_{0} f_+$. By \eqref{A3}, \eqref{A11-0}, and \eqref{app:Lprime-f0} we deduce
    \begin{align*}
        P'_{0} f_+=\frac{1}{2\pi \im}\oint_{\Gamma} \frac{1}{\lambda} \left(\mathscr{L}_{0}-\lambda\right)^{-1} \vet{\frac{\im(\cw^2+g)}{\cw}e^{2\im x}}{0} d\lambda.
    \end{align*}
Therefore by \eqref{A14} there is an analytic function $\lambda \mapsto \varphi(\lambda,\cdot)\in \mathcal Y$ so that we obtain
\begin{align*}
    P'_{0} f_+&=\frac{1}{2\pi \im} \oint_\Gamma \frac{1}{\lambda}\left(\left( \cw^2 (\partial^2_{x} + |D|) \right)^{-1} 
    \begin{bmatrix} \cw \partial_x+\omega\mathfrak H & -|D|  \\ g & \cw \partial_x \end{bmatrix} \vet{\frac{\im(\cw^2+g)}{\cw}e^{2\im x}}{0}
        +\lambda\varphi(\lambda)\right) d\lambda\\
        &=\frac{1}{2\pi \im} \oint_\Gamma \frac{1}{\lambda}\left(  e^{2\im x}\vet{\frac{(\cw^2+g)^2}{2\cw^4}}{-\frac{\im(\cw^2+g)g}{2\cw^3}}      +\lambda\varphi(\lambda)\right) d\lambda.
\end{align*}
Thus, by means of the residue theorem we obtain the first identity in \eqref{app:right-eps-zero}. Similarly, one may calculate $P'_{0} f_-$. Next, we calculate $P'_{0}f_0$. By \eqref{A11-0}, we have
\begin{align*}
        P'_{0} f_0=-\frac{1}{2\pi \im}\oint_{\Gamma} \frac{1}{\lambda} \left(\mathscr{L}_{0}-\lambda\right)^{-1} \vet{\frac{\cw^2+g}{\cw}\sin(x)}{(\cw^2+g)\cos(x)} d\lambda.
    \end{align*}
We note that
\begin{equation} \label{BBB}
   \vet{\frac{\cw^2+g}{\cw}\sin(x)}{(\cw^2+g)\cos(x)}=\frac{\im \omega}{2}\left(f_+-f_-\right)-\im \cw\left(f^-_1-f^-_{-1}\right).
\end{equation}
Therefore, combining \eqref{BBB}, \eqref{A11-0} and \eqref{A13} with the residue theorem, one may obtain
\begin{align*}
        P'_{0} f_0&=-\frac{1}{2\pi \im}\oint_{\Gamma} \frac{-\im \cw}{\lambda} \left(\mathscr{L}_{0}-\lambda\right)^{-1} \left(f^-_1-f^-_{-1}\right) d\lambda\\
        &=-\frac{1}{2\pi \im}\oint_{\Gamma} \frac{-\im \cw}{\lambda} \left(\frac{1}{\im(\cw^2+g)/\cw-\lambda}f^-_1-\frac{1}{-\im(\cw^2+g)/\cw-\lambda}f^-_{-1}\right) d\lambda\\
        &=\frac{\cw^2}{\cw^2+g}\left(f^-_1+f^-_{-1}\right).
    \end{align*}
This verifies the third identity in \eqref{app:right-eps-zero}.

Next, let us prove \eqref{app:right-mu-zero}. We notice that
\begin{align*}
    \dot{\mathscr{L}}_0 f_+=-\frac{\im \cw g}{\cw^2+g} f_+-\frac{\im \cw^3}{\cw^2+g} f^-_1.
\end{align*}
We deduce
\begin{align*}
    \dot{P}_0 f_+=\frac{1}{2\pi \im} \oint_{\Gamma} \frac{1}{\lambda}(\mathscr{L}_0-\lambda)^{-1} \left(\frac{\im \cw g}{\cw^2+g} f_+ +\frac{\im \cw^3}{\cw^2+g} f^-_1\right) d\lambda.
\end{align*}
By \eqref{A11-0} and \eqref{A13} with the residue theorem, we verify the first identity in \eqref{app:right-mu-zero}. We also observe that
\begin{align*}
    \dot{\mathscr{L}}_0 f_-=-\frac{\im \cw g}{\cw^2+g} f_- -\frac{\im \cw^3}{\cw^2+g} f^-_{-1}, \qquad \dot{\mathscr{L}}_0 f_0=\frac{\im g}{\omega} f_0+\frac{g^2}{\omega^2} \vet{0}{1}.
\end{align*}
Therefore, following a similar argument, we complete the proof of \eqref{app:right-mu-zero}. We have obtained the expansions \eqref{main:exp-fplus}, \eqref{main:exp-fminus}, and \eqref{main:exp-fzero}.
\end{proof}

Next, we compute the first-order expansions of $(U_{\mu,\e}^{-1})^*g_j$, $j\in\{+,-,0\}$, defined in \eqref{def:left-Kato-basis}. We have a similar lemma as follows.

\begin{lemma} \label{U inverse mu epsilon P00}
The first jets of $(U^{-1}_{\mu,\epsilon})^* P^*_{0}$ are
\begin{align} 
    &(U^{-1}_{0})^* P^*_{0} = P^*_{0}, \quad\quad \pa_\e(U^{-1}_{\mu,\e})^*\big|_{0} P^*_{0} = (P'_{0})^* P^*_{0}, \quad\quad \pa_\mu(U^{-1}_{\mu,\e})^*\big|_{0} P^*_{0} = (\dot{P}_{0})^* P^*_{0},
\end{align}
where
\begin{align} 
    (P'_{0})^* &= \frac{1}{2\pi \im} \oint_{\Gamma} (\mathscr{L}_{0}^* - \lambda)^{-1} (\mathscr{L}'_{0})^* (\mathscr{L}_{0}^* - \lambda)^{-1} d\lambda\,,  \\ 
    (\dot{P}_{0})^* &= \frac{1}{2\pi \im} \oint_{\Gamma} (\mathscr{L}_{0}^* - \lambda)^{-1} (\dot{\mathscr{L}}_{0})^* (\mathscr{L}_{0}^* - \lambda)^{-1} d\lambda\,.
\end{align}
The operators $(\mathscr{L}'_{0})^*$ and $(\dot{\mathscr{L}}_{0})^*$ are
\[
(\mathscr L_0')^*
=
\begin{bmatrix}
-p_1\partial_x &-a_1\\
0&-\partial_x\circ p_1
\end{bmatrix},
\qquad
(\dot{\mathscr L}_0)^*
=
\begin{bmatrix}
0&0\\
\mathrm{sgn}(D)+\Pi_0&0
\end{bmatrix},
\]
where $p_1$ and $a_1$ are defined in Lemma \ref{lem:pa-expansions}.
\end{lemma}

By the Lemma \ref{U inverse mu epsilon P00}, we have the Taylor expansion
\begin{equation} 
    \begin{aligned}
    g_j (\mu, \epsilon) &= g_j + \epsilon (P'_{0})^* g_j + \mu (\dot{P}_{0})^* g_j  + \mathcal{O}((|\mu|+|\e|)^2),\qquad j\in\{+,-,0\}.
\end{aligned}
\end{equation}
To expand the right basis, let us introduce the following lemma.
\begin{lemma} \label{Lemma B2-star}
The space $\mathcal Y=(H^1(\mathbb{T};\mathbb{C}))^2$ decomposes as $\mathcal Y = \mathcal{V}_{0}^* \oplus \mathcal{U}^* \oplus \mathcal{W}_{\mathcal Y}$, with $\mathcal{W}_{\mathcal Y}$ defined in Lemma \ref{Lemma B2}, where the subspaces $\mathcal{V}_{0}^*$, $\mathcal{U}^*$, and $\mathcal{W}_j$, defined below. We have the properties as follows.
\begin{itemize}
\item[(i)] $\mathcal{V}_{0}^* = \text{span}\{g_+, g_-, g_0\}$ is the kernel of $\mathscr{L}_{0}^*$. For any $\lambda \neq 0$ the operator $\mathscr{L}_{0}^* - \lambda : \mathcal{V}_{0}^* \to \mathcal{V}_{0}^*$ is invertible and
    \begin{equation} \label{A11-0-star}
        (\mathscr{L}_{0}^* - \lambda)^{-1} g_j = -\frac{1}{\lambda} g_j, \quad j\in\{+,-,0\}.
    \end{equation}

\item[(ii)] $\mathcal{U}^* := \text{span}\Big\{ g^-_1, g^-_{-1}, \vet{\im g/\omega}{1} \Big\}$, where
    \[
        g^-_{\pm1}:=e^{\pm\im x}\vet{1}{\pm\im/\cw}.
    \]
For any $\lambda \notin \{\pm \im(\cw^2+g)/\cw, \im \omega\}$ the operator $\mathscr{L}_{0}^* - \lambda : \mathcal{U}^* \to \mathcal{U}^*$ is invertible and
    \begin{equation} \label{A13-star}
    \begin{aligned}
        (\mathscr{L}_{0}^* - \lambda)^{-1} g^-_{\pm1}
        = \frac{1}{\mp\im(\cw^2+g)/\cw-\lambda} g^-_{\pm1},\qquad (\mathscr{L}_{0}^* - \lambda)^{-1}
        \vet{\im g/\omega}{1}
        = \frac{1}{\im\omega-\lambda}
        \vet{\im g/\omega}{1}.
    \end{aligned}
    \end{equation}

\item[(iii)] For $j\geq 2$, each subspace $\mathcal{W}_j := \text{span}\{ f_j^+, f_j^-, f_{-j}^+, f_{-j}^- \}$ is invariant under $\mathscr{L}_{0}^*$. For any $|\lambda| < \delta(\omega,g)$ small enough, the operator $\mathscr{L}_{0}^* - \lambda : \mathcal{W}_{\mathcal Y} \to \mathcal{W}_{\mathcal X}$ is invertible. Moreover, for every $f \in \mathcal{W}_{\mathcal X}$,
    \begin{equation} \label{A14-star}
        (\mathscr{L}_{0}^* - \lambda)^{-1} f
        =
        \left[\cw^2(\partial_x^2+|D|)\right]^{-1}
        \begin{bmatrix}
        -\cw\partial_x-\omega\mathfrak H & g\\
        -|D| & -\cw\partial_x
        \end{bmatrix}f
        +\lambda\varphi_f(\lambda),
    \end{equation}
for some analytic function $\lambda \mapsto \varphi_f(\lambda, \cdot) \in \mathcal Y=(H^1(\mathbb{T};\mathbb{C}))^2$.
\end{itemize}
\end{lemma}
Recalling $g_+$, $g_-$, and $g_0$ in Proposition \ref{prop:reference-spectrum}, one may first observe that
\begin{align}
(\mathscr L_0')^*g_+&=
-\frac{\im\omega}{2}g_0
+
e^{2\im x}
\begin{bmatrix}
-\dfrac{\im\omega}{2}\\[1mm]
1
\end{bmatrix},\qquad (\mathscr L_0')^*g_-=
\frac{\im\omega}{2}g_0
+
e^{-2\im x}
\begin{bmatrix}
\dfrac{\im\omega}{2}\\[1mm]
1
\end{bmatrix},\qquad(\mathscr L_0')^*g_0=0,
\label{app:Lprime-star-gzero}\\
(\dot{\mathscr L}_0)^*g_+
&=
\frac{\im g\cw}{\cw^2+g}g_+
-
\frac{\im g^2\cw}{(\cw^2+g)^2}g_1^-,
\quad
(\dot{\mathscr L}_0)^*g_-
=
\frac{\im g\cw}{\cw^2+g}g_-
-
\frac{\im g^2\cw}{(\cw^2+g)^2}g_{-1}^-,
\quad
(\dot{\mathscr L}_0)^*g_0
=
-\frac{\im g}{\omega}g_0
+
\vet{\im g/\omega}{1}.
\label{app:Ldot-star-g}
\end{align}
Using Lemma \ref{Lemma B2-star}, and applying the same contour argument as in Lemma \ref{lem:app-right-Kato-jets}, we obtain the following identities.
\begin{lemma}
\label{lem:app-left-Kato-jets}
One has
\begin{align}
(P'_{0})^*g_+
&=
e^{2\im x}
\begin{bmatrix}
\dfrac{(\cw^2+g)^2-2\cw^4}{4\cw^4}\\[2mm]
-\dfrac{\im(\cw^2+g)}{2\cw^3}
\end{bmatrix},\qquad (P'_{0})^* g_-=
e^{-2\im x}
\begin{bmatrix}
\dfrac{(\cw^2+g)^2-2\cw^4}{4\cw^4}\\[2mm]
\dfrac{\im(\cw^2+g)}{2\cw^3}
\end{bmatrix},\qquad (P'_{0})^* g_0=0,
\label{app:left-eps-zero}\\
(\dot{P}_{0})^* g_+
&=
e^{\im x}
\begin{bmatrix}
-\dfrac{g^2\cw^2}{(\cw^2+g)^3}\\[2mm]
-\dfrac{\im g^2\cw}{(\cw^2+g)^3}
\end{bmatrix},\qquad(\dot{P}_{0})^* g_-=
e^{-\im x}
\begin{bmatrix}
\dfrac{g^2\cw^2}{(\cw^2+g)^3}\\[2mm]
-\dfrac{\im g^2\cw}{(\cw^2+g)^3}
\end{bmatrix},\qquad(\dot{P}_{0})^* g_0=
\begin{bmatrix}
-\dfrac{g}{\omega^2}\\[2mm]
\dfrac{\im}{\omega}
\end{bmatrix}.
\label{app:left-mu-zero}
\end{align}
\end{lemma}

\section{Proof of Lemma \ref{lem:M1-properties-section5}}\label{M1-properties}
\begin{proof}
By Lemma~\ref{prop:quadratic-M-section4},
\[
\begin{pmatrix}
\mathbf M_{13}(0,\e)\\
\mathbf M_{23}(0,\e)
\end{pmatrix}
=\frac{\im\omega\e}{2}
\begin{pmatrix}1\\-1\end{pmatrix}
+\cO(\e^3).
\]
Since these entries are analytic and vanish at $\e=0$, the quotients in \eqref{def:R-epsilon-section5} extend analytically to $\e=0$. Moreover,
\[
R_\e=R_0+\cO(\e^2),\qquad
R_0=
\begin{pmatrix}
\frac12&\frac12&0\\
\frac12&-\frac12&0\\
0&0&1
\end{pmatrix},\qquad
\det R_0=-\frac12.
\]
Hence $R_\e$ is analytic and invertible for sufficiently small $|\e|$.

Write $\mathbf M_0:=\mathbf M(0,\e)$ and let $e_j$ denote the standard coordinate vectors. For $\e\ne0$, the definition of $R_\e$ gives
\[
R_\e e_2=\frac{\mathbf M_0e_3}{\im\omega\e},
\qquad R_\e e_3=e_3.
\]
Using $\mathbf M_0^2=0$, we obtain
\[
\mathbf M_0R_\e e_2=0,
\qquad
\mathbf M_0R_\e e_3=\im\omega\e R_\e e_2.
\]
The last row remains zero under this block diagonal change of coordinates. Therefore,
\[
\mathbf M^{(1)}(0,\e)=
\begin{pmatrix}
\mathbf M^{(1)}_{11}(0,\e)&0&0\\
\mathbf M^{(1)}_{21}(0,\e)&0&\im\omega\e\\
0&0&0
\end{pmatrix}.
\]
Its square is zero, so $\big(\mathbf M^{(1)}_{11}(0,\e)\big)^2=0$. Setting $a(\e):=\mathbf M^{(1)}_{21}(0,\e)$ proves \eqref{eq:M1-at-zero-section5}; $a$ is analytic, and the identity at $\e=0$ follows by continuity.

To prove \eqref{eq:M1-mu-derivatives-section5}, regard $\mathcal F_{\mu,\e}$ and $\mathcal G_{\mu,\e}$ as the right and left basis maps, and set
\[
u_k(\mu,\e):=\mathcal F_{\mu,\e}R_\e e_k,
\qquad
h_j(\mu,\e):=\mathcal G_{\mu,\e}(R_\e^{-1})^*e_j.
\]
These bases are biorthogonal, and
\[
\mathbf M^{(1)}_{jk}(\mu,\e)
=\big(\mathscr L_{\mu,\e}u_k(\mu,\e),
h_j(\mu,\e)\big)_{\mathcal X}.
\]
By \eqref{eq:M1-at-zero-section5} and invariance of the dual spectral subspace,
\begin{align} \label{u2 h1 h3}
\mathscr L_{0,\e}u_2(0,\e)=0,
\qquad
\mathscr L_{0,\e}^*h_j(0,\e)=0,
\quad j=1,3.
\end{align}
Differentiating the preceding matrix entry and using \eqref{eq:scriptL-family} thus gives
\begin{align*}
\partial_\mu\mathbf M^{(1)}_{j2}(0,\e)
&=\big(\mathscr A_\e u_2,h_j\big)_{\mathcal X}
+\big(\mathscr L_{0,\e}\partial_\mu u_2,h_j\big)_{\mathcal X}
+\big(\mathscr L_{0,\e}u_2,\partial_\mu h_j\big)_{\mathcal X}\\
&=\big(\mathscr A_\e u_2(0,\e),h_j(0,\e)\big)_{\mathcal X},
\qquad j=1,3,
\end{align*}
where all vectors in the first line are evaluated at $(0,\e)$.

We now establish the required parity for $\e\ne0$. Recalling Lemma \ref{lem:one-sided-Jordan-section5}, we have
\[
E_2(\e)=\vet{E_{21}(\e)}{E_{22}(\e)}=
-\frac{\im}{\e}
\left(
\mathcal P Z^{-1}
\begin{bmatrix}
\partial_x\eta_\e\\
\partial_x\psi_\e
\end{bmatrix}
-m_\e \vet{0}{1}
\right)\in \vet{\mathrm{odd}(x)}{\mathrm{even}_0(x)}.
\]
Its components are odd and even, respectively, and both have zero mean. The exact chain \eqref{exact-chain-section5} shows that
\[
\operatorname{Ran}
\big(\mathscr L_{0,\e}|_{\mathcal V_{0,\e}}\big)
=\operatorname{span}\{E_2(\e)\}.
\]
On the other hand, the intertwining identity and the definition of $R_\e$ give
\[
u_2(0,\e)
=\frac{\mathcal F_{0,\e}\mathbf M_0e_3}{\im\omega\e}
=\frac{\mathscr L_{0,\e}f_0(0,\e)}{\im\omega\e}\in \operatorname{span}\{E_2(\e)\}.
\]
Hence $u_2(0,\e)$ has an odd first component and an even, zero-mean second component.

To identify the parity of the adjoint kernel, define
\[
\widehat h_\e:=
\begin{bmatrix}
-E_{22}(\e)+\omega\partial_x^{-1}E_{21}(\e)\\
E_{21}(\e)
\end{bmatrix}.
\]
Its first component is even and its second component is odd. Since $\mathscr L_{0,\e}E_2(\e)=0$,
\begin{align*}
\partial_x\big((\cw+p_\e)E_{21}(\e)\big)
+|D| E_{22}(\e)&=0,\\
-(g+a_\e)E_{21}(\e)
+(\cw+p_\e)\partial_x E_{22}(\e)
+\omega\mathfrak H E_{22}(\e)&=0.
\end{align*}
Using $|D|=\partial_x\mathfrak H$, integration of the first equation gives
\[
(\cw+p_\e) E_{21}(\e)+\mathfrak H E_{22}(\e)=0,
\]
because its left-hand side is odd. Taking the adjoint in \eqref{eq:scriptL-zero-plus} and using $|D|\partial_x^{-1} E_{21}(\e)=\mathfrak H E_{21}(\e)$ yields
\[
\mathscr L_{0,\e}^*\widehat h_\e
=\begin{bmatrix}
-\omega\big((\cw+p_\e)E_{21}(\e)+\mathfrak H E_{22}(\e)\big)\\
-|D| E_{22}(\e)-\partial_x\big((\cw+p_\e)E_{21}(\e)\big)
\end{bmatrix}
=0.
\]
Also, $\mathscr L_{0,\e}^*(1,0)^\top=0$. For $\e\ne0$, the exact chain has rank one on the three-dimensional generalized kernel, and the operator is invertible on its complementary spectral subspace. Thus $\dim\ker\mathscr L_{0,\e}^*=2$. Since $E_{21}(\e)=\im\sin (x)+\im\cO(\e)$, the two adjoint kernel vectors are independent for sufficiently small nonzero $\e$, giving
\[
\ker\mathscr L_{0,\e}^*
=\operatorname{span}\left\{
\begin{bmatrix}1\\0\end{bmatrix},\widehat h_\e
\right\}.
\]
Recalling \eqref{u2 h1 h3}, we know $h_1(0,\e),\,h_3(0,\e)\in \ker\mathscr L_{0,\e}^*$. Consequently, $h_1(0,\e)$ and $h_3(0,\e)$ have even first components and odd second components.

Finally, $p_\e$ is even, and $\mathrm{sgn}(D)$ maps even functions to odd functions. The zero mean of the second component of $u_2(0,\e)$ removes the $\Pi_0$ term in \eqref{eq:A-eps}. Hence $\mathscr A_\e u_2(0,\e)$ has an odd first component and an even second component. It follows that
\[
\big(\mathscr A_\e u_2(0,\e),h_j(0,\e)\big)_{\mathcal X}=0,
\qquad j=1,3.
\]
This proves \eqref{eq:M1-mu-derivatives-section5} for $\e\ne0$, and the identities at $\e=0$ follow by analyticity.
\end{proof}

\section{Calculation of the bases and matrix at the first isola}
\label{sec:AppIsolaKato}

Throughout this appendix, we assume \eqref{isola:range} and $\bar\mu\ne0$, and all derivatives are evaluated at $(\delta,\e)=(0,0)$. We use the reference quantities in \eqref{isola:reference-parameters} and write
\[
\mathcal L_\delta:=
\begin{bmatrix}
\im\cw&\operatorname{sign}(D+\bar\mu)\\0&\im\cw
\end{bmatrix},\qquad
\mathcal L_n:=
\begin{bmatrix}
(\partial_x+\im\bar\mu)\circ p_n&0\\
-a_n&p_n(\partial_x+\im\bar\mu)
\end{bmatrix},\quad n=1,2.
\]
Thus $\mathcal L_\delta=\partial_\delta \mathcal L_{\bar\mu+\delta,\e}|_{(0,0)}$ and $\mathcal L_{\bar\mu,\e}=\mathcal L_{\bar\mu,0}+\e\mathcal L_1+ \e^2\mathcal L_2+\cO(\e^3)$. As in Appendix \ref{sec:AppKato}, a dot denotes the Bloch derivative and a prime denotes the amplitude derivative.

To compute the Kato derivatives, we first describe the action of $(\mathcal L_{\bar\mu,0}-\lambda)^{-1}$, following Lemma \ref{Lemma B2}. The other eigenvector on each resonant Fourier mode is
\begin{equation}\label{app:isola-complement-right}
\widetilde f_1=
\begin{bmatrix}1\\-\dfrac{2\im g}{\cw}\end{bmatrix}e^{\im k_1x},
\qquad
\widetilde f_2=
\begin{bmatrix}1\\\dfrac{2\im g}{3\cw}\end{bmatrix}e^{\im k_2x}.
\end{equation}
For $k\in\mathbb Z\setminus\{k_1,k_2\}$, set
\[
\mathcal W_k:=\{v e^{\im kx}:v\in\mathbb C^2\},\qquad
\mathcal W_{\mathcal X}:=
\overline{\bigoplus_{k\notin\{k_1,k_2\}}\mathcal W_k}^{\mathcal X},
\qquad
\mathcal W_{\mathcal Y}:=
\overline{\bigoplus_{k\notin\{k_1,k_2\}}\mathcal W_k}^{\mathcal Y}.
\]

\begin{lemma}\label{lem:app-isola-resolvent}
The periodic spaces decompose as
\begin{equation}\label{app:isola-space-decomposition}
\mathcal X=\mathcal V_0\oplus\mathcal U\oplus\mathcal W_{\mathcal X},
\qquad
\mathcal Y=\mathcal V_0\oplus\mathcal U\oplus\mathcal W_{\mathcal Y},
\end{equation}
where $\mathcal V_0=\operatorname{span}\{f_1,f_2\}$ and $\mathcal U=\operatorname{span}\{\widetilde f_1,\widetilde f_2\}$. The operator $\mathcal L_{\bar\mu,0}$ maps each summand of $\mathcal Y$ into the corresponding summand of $\mathcal X$. Its restrictions have the following resolvent formulas.
\begin{itemize}
\item[(i)] On $\mathcal V_0$, for $\lambda\ne\im s_*$,
\begin{equation}\label{app:isola-resolvent-V}
(\mathcal L_{\bar\mu,0}-\lambda)^{-1}f_a
=\frac1{\im s_*-\lambda}f_a,\qquad a=1,2.
\end{equation}
\item[(ii)] On $\mathcal U$, for $\lambda\notin\{\im(s_*-\cw-\omega),\im(s_*+3\cw-\omega)\}$,
\begin{equation}\label{app:isola-resolvent-U}
\begin{aligned}
(\mathcal L_{\bar\mu,0}-\lambda)^{-1}\widetilde f_1
=\frac{\widetilde f_1}{\im(s_*-\cw-\omega)-\lambda},\qquad (\mathcal L_{\bar\mu,0}-\lambda)^{-1}\widetilde f_2
=\frac{\widetilde f_2}{\im(s_*+3\cw-\omega)-\lambda}.
\end{aligned}
\end{equation}
\item[(iii)] For $|\lambda-\im s_*|$ sufficiently small, $\mathcal L_{\bar\mu,0}-\lambda:\mathcal W_{\mathcal Y}\to \mathcal W_{\mathcal X}$ is invertible. For every $w\in\mathcal W_{\mathcal X}$,
\begin{equation}\label{app:isola-resolvent-W}
\begin{aligned}
(\mathcal L_{\bar\mu,0}-\lambda)^{-1}w
={}&\Big[
\big(\cw(\partial_x+\im\bar\mu)-\im s_*\big)^2
+\cw^2|D+\bar\mu|-\im\omega s_*\mathfrak H_{\bar\mu}
\Big]^{-1}\\
&\quad\times
\begin{bmatrix}
\cw(\partial_x+\im\bar\mu)+\omega\mathfrak H_{\bar\mu}-\im s_*
&-|D+\bar\mu|\\
g&\cw(\partial_x+\im\bar\mu)-\im s_*
\end{bmatrix}w+(\lambda-\im s_*)\varphi_w(\lambda),
\end{aligned}
\end{equation}
where $\varphi_w$ is analytic with values in $\mathcal W_{\mathcal Y}$. The scalar inverse is taken componentwise on the Fourier modes $k\notin\{k_1,k_2\}$. The operator product on the right is understood by continuous extension from trigonometric polynomials to $\mathcal W_{\mathcal X}$.
\end{itemize}
\end{lemma}
\begin{proof}
The proof is similar to the proof in Lemma \ref{Lemma B2}.
\end{proof}

\begin{proof}[Proof of Lemma \ref{lem:isola-kato-expansions}: right basis]
By Lemma \ref{U mu epsilon P00},
\begin{equation}\label{isola:first-jets}
\dot f_a=\dot P_0f_a,\qquad f_a'=P'_0f_a.
\end{equation}
Here, for $\nu=\delta,\e$, with $\mathcal L_\e:=\mathcal L_1$,
\begin{equation}\label{app:isola-projector-jets}
P_\nu=\frac1{2\pi\im}\oint_\Gamma
(\mathcal L_{\bar\mu,0}-\lambda)^{-1}\mathcal L_\nu
(\mathcal L_{\bar\mu,0}-\lambda)^{-1}\,d\lambda.
\end{equation}
The positively oriented circle $\Gamma$ is centered at $\im s_*$ and is chosen small enough for Lemma \ref{lem:app-isola-resolvent} to apply.

Applying $\mathcal L_\delta$ and decomposing the resulting vector along $f_a,\widetilde f_a$ gives
\begin{equation}\label{app:isola-Ldelta-f}
\begin{aligned}
\mathcal L_\delta f_1
=\frac{2\im\cw^3}{2\cw^2-g}f_1
+\frac{\im\cw^3g}{(3\cw^2-2g)(2\cw^2-g)}\widetilde f_1,\qquad
\mathcal L_\delta f_2
=\frac{2\im\cw^3}{2\cw^2+g}f_2
-\frac{3\im\cw^3g}{(\cw^2+2g)(2\cw^2+g)}\widetilde f_2.
\end{aligned}
\end{equation}
Inserting \eqref{app:isola-resolvent-V}--\eqref{app:isola-resolvent-U} into \eqref{app:isola-projector-jets}, the terms along $f_1,f_2$ have double poles with zero residue. The remaining terms give
\begin{align*}
\dot f_1
&=\frac{\im\cw^3g}{(3\cw^2-2g)(2\cw^2-g)}
\frac1{2\pi\im}\oint_\Gamma
\frac{\widetilde f_1\,d\lambda}
{(\im s_*-\lambda)(\im(s_*-\cw-\omega)-\lambda)}=\frac{\cw^4g}{(3\cw^2-2g)(2\cw^2-g)^2}\widetilde f_1,\\
\dot f_2
&=-\frac{3\im\cw^3g}{(\cw^2+2g)(2\cw^2+g)}
\frac1{2\pi\im}\oint_\Gamma
\frac{\widetilde f_2\,d\lambda}
{(\im s_*-\lambda)(\im(s_*+3\cw-\omega)-\lambda)}=\frac{3\cw^4g}{(\cw^2+2g)(2\cw^2+g)^2}\widetilde f_2.
\end{align*}
This proves \eqref{app:isola-right-delta}.

Using $2\cos (x)=e^{\im x}+e^{-\im x}$ gives
\begin{equation}\label{app:isola-Lprime-f}
\begin{aligned}
\mathcal L_1f_1=-\frac{\cw^2+g}{2\cw}\sum_{\sigma=\pm1}
\begin{bmatrix}\im(\xi_1+\sigma)\\-3\cw/2\end{bmatrix}
e^{\im(k_1+\sigma)x},\qquad
\mathcal L_1f_2=-\frac{\cw^2+g}{2\cw}\sum_{\sigma=\pm1}
\begin{bmatrix}\im(\xi_2+\sigma)\\\cw/2\end{bmatrix}
e^{\im(k_2+\sigma)x}.
\end{aligned}
\end{equation}
Since $k_2-k_1=2$, both vectors belong to $\mathcal W_{\mathcal Y}$. Hence \eqref{app:isola-resolvent-V} and \eqref{app:isola-resolvent-W} give
\begin{align*}
f_a'
&=-\frac1{2\pi\im}\oint_\Gamma
\frac{(\mathcal L_{\bar\mu,0}-\lambda)^{-1}\mathcal L_1f_a}
{\lambda-\im s_*}\,d\lambda=-\frac1{2\pi\im}\oint_\Gamma
\left[
\frac{\left[(\mathcal L_{\bar\mu,0}-\im s_*I)
|_{\mathcal W_{\mathcal Y}}\right]^{-1}\mathcal L_1f_a}
{\lambda-\im s_*}
+\varphi_{\mathcal L_1f_a}(\lambda)
\right]d\lambda.
\end{align*}
The analytic term integrates to zero. Therefore,
\begin{equation}\label{app:isola-right-eps-residue}
f_a'=-\left[(\mathcal L_{\bar\mu,0}-\im s_*I)
|_{\mathcal W_{\mathcal Y}}\right]^{-1}\mathcal L_1f_a.
\end{equation}
A straightforward computation gives
\begin{equation*}
-\mathcal L_1f_2
=\frac{\cw^2+g}{2\cw}
\left(
\begin{bmatrix}
\im(\xi_2-1)\\[1mm]
\dfrac{\cw}{2}
\end{bmatrix}
e^{\im(j_*+1)x}
+
\begin{bmatrix}
\im(\xi_2+1)\\[1mm]
\dfrac{\cw}{2}
\end{bmatrix}
e^{\im(j_*+3)x}
\right),
\end{equation*}
and, hence, by \eqref{app:isola-resolvent-W}, we obtain $f'_2$ in \eqref{app:isola-right-eps}.

Similarly, we obtain
\begin{equation} \label{-L1 f1}
-\mathcal L_1f_1
=\frac{\cw^2+g}{2\cw}
\left(
\begin{bmatrix}
\im(\xi_1-1)\\[1mm]
-\dfrac{3\cw}{2}
\end{bmatrix}
e^{\im(j_*-1)x}
+
\begin{bmatrix}
\im(\xi_1+1)\\[1mm]
-\dfrac{3\cw}{2}
\end{bmatrix}
e^{\im(j_*+1)x}
\right).
\end{equation}
We notice that \eqref{-L1 f1} has two modes $j_*-1$ and $j_*+1$. Also,
\begin{align} \label{c large or small}
    j_*+1+\bar\mu=\frac{3\cw^2+2g}{4g}>0,\quad\text{but}\quad j_*-1+\bar\mu=\frac{3(\cw^2-2g)}{4g} \begin{cases}
\displaystyle
< 0,
&\displaystyle\frac{2g}{3}<\cw^2<2g,
\\ > 0,
&\cw^2>2g.
\end{cases}
\end{align}
We pause to remark that if $\cw^2=2g$, we have $\bar\mu=0$ (see \eqref{isola:anchor}). Combining \eqref{-L1 f1} and \eqref{c large or small} we obtain $f'_1$ in \eqref{app:isola-right-eps}.
\end{proof}

For the left basis, we proceed as in Lemma \ref{Lemma B2-star}. The adjoint reference operator and its first derivative with respect to $\delta$ are
\[
\mathcal L_{\bar\mu,0}^*=
\begin{bmatrix}
-\cw(\partial_x+\im\bar\mu)&-g\\
|D+\bar\mu|&-\cw(\partial_x+\im\bar\mu)-\omega\mathfrak H_{\bar\mu}
\end{bmatrix},
\qquad
\mathcal L_\delta^*=
\begin{bmatrix}
-\im\cw&0\\
\operatorname{sign}(D+\bar\mu)&-\im\cw
\end{bmatrix}.
\]
Define the complementary adjoint eigenvectors
\begin{equation}\label{app:isola-complement-left}
\begin{aligned}
\widetilde g_1&=
\begin{bmatrix}
\dfrac{\cw^2}{2(2\cw^2-g)}\\[1mm]
-\dfrac{\im\cw(3\cw^2-2g)}{4g(2\cw^2-g)}
\end{bmatrix}e^{\im k_1x},\qquad
\widetilde g_2=
\begin{bmatrix}
\dfrac{3\cw^2}{2(2\cw^2+g)}\\[1mm]
\dfrac{3\im\cw(\cw^2+2g)}{4g(2\cw^2+g)}
\end{bmatrix}e^{\im k_2x}.
\end{aligned}
\end{equation}

Similar to Lemma \ref{lem:app-isola-resolvent}, let us introduce the following useful lemma.
\begin{lemma}\label{lem:app-isola-adjoint-resolvent}
The periodic spaces decompose as
\[
\mathcal X=\mathcal V_0^*\oplus\mathcal U^*\oplus\mathcal W_{\mathcal X},
\qquad
\mathcal Y=\mathcal V_0^*\oplus\mathcal U^*\oplus\mathcal W_{\mathcal Y},
\]
where $\mathcal V_0^*=\operatorname{span}\{g_1,g_2\}$ and $\mathcal U^*=\operatorname{span}\{\widetilde g_1,\widetilde g_2\}$. The operator $\mathcal L_{\bar\mu,0}^*$ maps each summand of $\mathcal Y$ into the corresponding summand of $\mathcal X$, with the following resolvent formulas.
\begin{itemize}
\item[(i)] On $\mathcal V_0^*$, for $\lambda\ne-\im s_*$,
\begin{equation}\label{app:isola-adjoint-resolvent-V}
(\mathcal L_{\bar\mu,0}^*-\lambda)^{-1}g_a
=\frac{g_a}{-\im s_*-\lambda},\qquad a=1,2.
\end{equation}
\item[(ii)] On $\mathcal U^*$, for $\lambda\notin\{-\im(s_*-\cw-\omega),-\im(s_*+3\cw-\omega)\}$,
\begin{equation}\label{app:isola-adjoint-resolvent-U}
\begin{aligned}
(\mathcal L_{\bar\mu,0}^*-\lambda)^{-1}\widetilde g_1
=\frac{\widetilde g_1}{-\im(s_*-\cw-\omega)-\lambda},\qquad(\mathcal L_{\bar\mu,0}^*-\lambda)^{-1}\widetilde g_2
=\frac{\widetilde g_2}{-\im(s_*+3\cw-\omega)-\lambda}.
\end{aligned}
\end{equation}
\item[(iii)] For $|\lambda+\im s_*|$ sufficiently small, $\mathcal L_{\bar\mu,0}^*-\lambda:\mathcal W_{\mathcal Y}\to \mathcal W_{\mathcal X}$ is invertible. For every $w\in\mathcal W_{\mathcal X}$,
\begin{equation}\label{app:isola-adjoint-resolvent-W}
\begin{aligned}
(\mathcal L_{\bar\mu,0}^*-\lambda)^{-1}w
={}&\Big[
\big(\cw(\partial_x+\im\bar\mu)-\im s_*\big)^2
+\cw^2|D+\bar\mu|-\im\omega s_*\mathfrak H_{\bar\mu}
\Big]^{-1}\\
&\quad\times
\begin{bmatrix}
-\cw(\partial_x+\im\bar\mu)-\omega\mathfrak H_{\bar\mu}+\im s_*&g\\
-|D+\bar\mu|&-\cw(\partial_x+\im\bar\mu)+\im s_*
\end{bmatrix}w+(\lambda+\im s_*)\psi_w(\lambda),
\end{aligned}
\end{equation}
where $\psi_w$ is analytic with values in $\mathcal W_{\mathcal Y}$. The scalar inverse acts only on the modes $k\notin\{k_1,k_2\}$.
\end{itemize}
\end{lemma}

\begin{proof}[Proof of Lemma \ref{lem:isola-kato-expansions}: left basis]
By Lemma \ref{U inverse mu epsilon P00},
\[
\dot g_a=(\dot P_0)^*g_a.
\]
Let $\Gamma_*$ be the positively oriented circle of the same radius as $\Gamma$, centered at $-\im s_*$. Then
\begin{equation}\label{app:isola-adjoint-projector-jets}
P_\delta^*=\frac1{2\pi\im}\oint_{\Gamma_*}
(\mathcal L_{\bar\mu,0}^*-\lambda)^{-1}\mathcal L_\delta^*
(\mathcal L_{\bar\mu,0}^*-\lambda)^{-1}\,d\lambda.
\end{equation}
Since $k_a+\bar\mu=\xi_a>0$, the sign multiplier equals $1$ on both reference modes. Thus direct application of $\mathcal L_\delta^*$ gives
\begin{equation}\label{app:isola-Ldelta-star-g}
\begin{aligned}
\mathcal L_\delta^*g_1
=-\frac{2\im\cw^3}{2\cw^2-g}g_1
+\frac{\im g(3\cw^2-2g)}{\cw(2\cw^2-g)}\widetilde g_1,\qquad \mathcal L_\delta^*g_2
=-\frac{2\im\cw^3}{2\cw^2+g}g_2
-\frac{\im g(\cw^2+2g)}{3\cw(2\cw^2+g)}\widetilde g_2.
\end{aligned}
\end{equation}
Using items (i)--(ii) of Lemma \ref{lem:app-isola-adjoint-resolvent}, we obtain
\begin{align*}
\dot g_1=-\frac{g(3\cw^2-2g)}{(2\cw^2-g)^2}\widetilde g_1,\qquad
\dot g_2=-\frac{g(\cw^2+2g)}{3(2\cw^2+g)^2}\widetilde g_2.
\end{align*}
Here we used $\cw+\omega=(2\cw^2-g)/\cw$ and $3\cw-\omega=(2\cw^2+g)/\cw$. This proves \eqref{app:isola-left-delta}.
\end{proof}

\begin{proof}[Proof of Lemma \ref{lem:isola-coefficients}]
By \eqref{isola:flat-biorthogonality},
\[
\mathbf H(\bar\mu,0)
=\big((\mathcal L_{\bar\mu,0}f_b,g_a)_{\mathcal X}\big)_{a,b\in\{1,2\}}
=\im s_*I_2.
\]
Recall the fact $(f_b(\delta,\e),g_a(\delta,\e))_{\mathcal X}=\delta_{ab}$ and, hence, $(f_b',g_a)_{\mathcal X}+(f_b,g_a')_{\mathcal X}=0$. We have \begin{align*}
\partial_\e\mathbf H_{ab}(\bar\mu,0)
&=(\mathcal L_1f_b,g_a)_{\mathcal X}+(\mathcal L_{\bar\mu,0}f'_b,g_a)_{\mathcal X}+(\mathcal L_{\bar\mu,0}f_b,g'_a)_{\mathcal X}=(\mathcal L_1f_b,g_a)_{\mathcal X}+\im s_*\big[(f_b',g_a)_{\mathcal X}
+(f_b,g_a')_{\mathcal X}\big]=0.
\end{align*}
The last equality follows from Fourier orthogonality (see \eqref{app:isola-Lprime-f}).

Similarly, we have $(\dot f_b,g_a)_{\mathcal X}+(f_b,\dot g_a)_{\mathcal X}=0$. Thus
\begin{align*}
\partial_\delta\mathbf H_{ab}(\bar\mu,0)
&=(\mathcal L_\delta f_b,g_a)_{\mathcal X}
+\im s_*\big[(\dot f_b,g_a)_{\mathcal X}
+(f_b,\dot g_a)_{\mathcal X}\big]=(\mathcal L_\delta f_b,g_a)_{\mathcal X}.
\end{align*}
Recall \eqref{app:isola-Ldelta-f}. Fourier orthogonality gives $(\mathcal L_\delta f_b,g_a)_{\mathcal X}=0$ for $a\ne b$, while substitution of $g_1,g_2$ yields
\[
(\mathcal L_\delta f_1,g_1)_{\mathcal X}
=\frac{2\im\cw^3}{2\cw^2-g},
\qquad
(\mathcal L_\delta f_2,g_2)_{\mathcal X}
=\frac{2\im\cw^3}{2\cw^2+g}.
\]
We now compute the coefficient of $\e^2$. All basis derivatives below are evaluated at $(\delta,\e)=(0,0)$. Since $\mathcal L_{\bar\mu,\e}=\mathcal L_{\bar\mu,0} +\e\mathcal L_1+\e^2\mathcal L_2+\cO(\e^3)$, differentiating the matrix entries twice gives
\begin{align*}
\frac12\partial_\e^2\mathbf H_{ab}(\bar\mu,0)
={}&(\mathcal L_2f_b,g_a)_{\mathcal X}
+(\mathcal L_1f_b',g_a)_{\mathcal X}
+(\mathcal L_1f_b,g_a')_{\mathcal X}+\frac12(\mathcal L_{\bar\mu,0}f_b'',g_a)_{\mathcal X}
+(\mathcal L_{\bar\mu,0}f_b',g_a')_{\mathcal X}
+\frac12(\mathcal L_{\bar\mu,0}f_b,g_a'')_{\mathcal X}.
\end{align*}
The second derivative of $(f_b(\delta,\e),g_a(\delta,\e))_{\mathcal X}=\delta_{ab}$ yields
\[
(f_b'',g_a)_{\mathcal X}
+2(f_b',g_a')_{\mathcal X}
+(f_b,g_a'')_{\mathcal X}=0.
\]
Using the left and right eigenvector identities, we therefore obtain
\begin{align*}
&\frac12(\mathcal L_{\bar\mu,0}f_b'',g_a)_{\mathcal X}
+\frac12(\mathcal L_{\bar\mu,0}f_b,g_a'')_{\mathcal X}=\frac{\im s_*}{2}
\big[(f_b'',g_a)_{\mathcal X}+(f_b,g_a'')_{\mathcal X}\big]
=-\im s_*(f_b',g_a')_{\mathcal X}.
\end{align*}
Consequently,
\begin{align*}
\frac12\partial_\e^2\mathbf H_{ab}(\bar\mu,0)
={}&(\mathcal L_2f_b+\mathcal L_1f_b',g_a)_{\mathcal X}+\big(\mathcal L_1f_b
+(\mathcal L_{\bar\mu,0}-\im s_*I)f_b',g_a'\big)_{\mathcal X}=(\mathcal L_2f_b+\mathcal L_1f_b',g_a)_{\mathcal X}.
\end{align*}
The last equality holds because $(\mathcal L_{\bar\mu,0}-\im s_*I)f_a'=-\mathcal L_1f_a$ (see \eqref{app:isola-right-eps-residue}). We first compute $\mathcal L_2f_1$ and $\mathcal L_2f_2$. Recalling \eqref{isola:reference-parameters} and
\[
\xi_1=\frac{3\cw^2-2g}{4g},\qquad
\xi_2=\frac{3(\cw^2+2g)}{4g},
\]
direct multiplication gives
\begin{align*}
\mathcal L_2f_1=\begin{bmatrix}
\partial_xp_2+\im\xi_1p_2\\[1mm]
-a_2-\dfrac{\cw}{2}p_2
\end{bmatrix}e^{\im k_1x},\qquad
\mathcal L_2f_2=\begin{bmatrix}
\partial_xp_2+\im\xi_2p_2\\[1mm]
-a_2+\dfrac{3\cw}{2}p_2
\end{bmatrix}e^{\im k_2x}.
\end{align*}
Substituting $p_2=p_2^{[0]}+p_2^{[2]}\cos(2x)$ and $a_2=a_2^{[0]}+a_2^{[2]}\cos(2x)$ from \eqref{exp:p2}--\eqref{exp:a22}, we obtain
\begin{equation}\label{app:isola-L2-right}
\begin{aligned}
\mathcal L_2f_1
={}&\begin{bmatrix}
\im\xi_1p_2^{[0]}\\[1mm]
-a_2^{[0]}-\dfrac{\cw}{2}p_2^{[0]}
\end{bmatrix}e^{\im k_1x}+\frac12\sum_{\sigma=\pm1}
\begin{bmatrix}
\im(\xi_1+2\sigma)p_2^{[2]}\\[1mm]
-a_2^{[2]}-\dfrac{\cw}{2}p_2^{[2]}
\end{bmatrix}e^{\im(k_1+2\sigma)x},\\[2mm]
\mathcal L_2f_2
={}&\begin{bmatrix}
\im\xi_2p_2^{[0]}\\[1mm]
-a_2^{[0]}+\dfrac{3\cw}{2}p_2^{[0]}
\end{bmatrix}e^{\im k_2x}+\frac12\sum_{\sigma=\pm1}
\begin{bmatrix}
\im(\xi_2+2\sigma)p_2^{[2]}\\[1mm]
-a_2^{[2]}+\dfrac{3\cw}{2}p_2^{[2]}
\end{bmatrix}e^{\im(k_2+2\sigma)x}.
\end{aligned}
\end{equation}
Taking the inner products with the dual basis vectors in \eqref{isola:flat-left-basis} and using Fourier orthogonality, we obtain
\begin{equation} \label{L2 f g}
\begin{aligned}
(\mathcal L_2f_1,g_1)_{\mathcal X}
&=\frac{\im(\cw^2+g)(3\cw^2-2g)
(14\cw^6+g\cw^4-g^3)}
{32\cw^5g(2\cw^2-g)},\\(\mathcal L_2f_2,g_1)_{\mathcal X}
&=\frac{\im(\cw^2-g)^2(\cw^2+g)(3\cw^2-2g)(4\cw^2+g)}
{32\cw^5g(2\cw^2-g)},\\
(\mathcal L_2f_1,g_2)_{\mathcal X}
&=-\frac{3\im(\cw^2-g)^2(\cw^2+g)(\cw^2+2g)(4\cw^2+g)}
{32\cw^5g(2\cw^2+g)},\\
(\mathcal L_2f_2,g_2)_{\mathcal X}
&=\frac{3\im(\cw^2+g)(\cw^2+2g)
(-2\cw^6+7g\cw^4+4g^2\cw^2+g^3)}
{32\cw^5g(2\cw^2+g)}.
\end{aligned}
\end{equation}
Substituting \eqref{app:isola-right-eps} into the action of $\mathcal L_1$ and collecting equal Fourier modes gives the following formulas.

For $2g/3<\cw^2<2g$,
\begin{equation}\label{app:isola-L1-f1prime-low}
\begin{aligned}
\mathcal L_1f_1'
={}&\frac{(\cw^2+g)^2}{64\cw^5g^2}
\Bigg(
\begin{bmatrix}
-2\im(3\cw^2-10g)(\cw^2-g)(\cw^2-2g)\\[1mm]
-4\cw g(\cw^2-2g)(4\cw^2-g)
\end{bmatrix}
e^{\im(j_*-2)x}
\\[2mm]
&\qquad+
\begin{bmatrix}
-\im(3\cw^2-2g)(8\cw^4-5g\cw^2+2g^2)\\[1mm]
2\cw g(4\cw^2-g)(\cw^2+6g)
\end{bmatrix}
e^{\im j_*x}+
\begin{bmatrix}
-3\im(\cw^2+2g)(3\cw^2+2g)(2\cw^2-g)\\[1mm]
2\cw g(3\cw^2+2g)(4\cw^2-g)
\end{bmatrix}
e^{\im(j_*+2)x}
\Bigg).
\end{aligned}
\end{equation}

For $\cw^2>2g$,
\begin{equation}\label{app:isola-L1-f1prime-high}
\begin{aligned}
\mathcal L_1f_1'
={}&\frac{(\cw^2+g)^2}{64\cw^5g^2}
\Bigg(
\begin{bmatrix}
\im(3\cw^2-10g)(\cw^2-2g)(4\cw^2-g)\\[1mm]
-2\cw g(\cw^2-2g)(14\cw^2-5g)
\end{bmatrix}
e^{\im(j_*-2)x}
\\[2mm]
&\qquad+
\begin{bmatrix}
-2\im(3\cw^2-2g)(\cw^4+5g\cw^2-2g^2)\\[1mm]
-4\cw g(\cw^4-19g\cw^2+6g^2)
\end{bmatrix}
e^{\im j_*x}+
\begin{bmatrix}
-3\im(\cw^2+2g)(3\cw^2+2g)(2\cw^2-g)\\[1mm]
2\cw g(3\cw^2+2g)(4\cw^2-g)
\end{bmatrix}
e^{\im(j_*+2)x}
\Bigg).
\end{aligned}
\end{equation}

In both cases,
\begin{equation}\label{app:isola-L1-f2prime}
\begin{aligned}
\mathcal L_1f_2'
={}&
\frac{(\cw^2+g)^2(3\cw^2+2g)}{64\cw^5g}
\begin{bmatrix}
\im(3\cw^2-2g)\\[1mm]
2\cw(2\cw^2-g)
\end{bmatrix}e^{\im j_*x}
\\[2mm]
&+
\frac{(\cw^2+g)^2}{32\cw^5g^2}
\begin{bmatrix}
-\im(\cw^2+2g)^2(3\cw^2+g)\\[1mm]
2\cw g(3\cw^4-g\cw^2-6g^2)
\end{bmatrix}e^{\im(j_*+2)x}-
\frac{(\cw^2+g)^2(3\cw^2+10g)}{192\cw^5g^2}
\begin{bmatrix}
\im(2\cw^2+g)(3\cw^2+14g)\\[1mm]
6\cw g^2
\end{bmatrix}e^{\im(j_*+4)x}.
\end{aligned}
\end{equation}

Only the modes $j_*$ and $j_*+2$ contribute to the inner products with $g_1$ and $g_2$, respectively. Taking the inner products with the dual basis vectors in \eqref{isola:flat-left-basis} and using Fourier orthogonality, we obtain
\begin{equation}\label{app:isola-L1-inner-products}
\begin{aligned}
(\mathcal L_1f_1',g_1)_{\mathcal X}
&=
\begin{cases}
-\displaystyle
\frac{\im(\cw^2+g)^2(3\cw^2-2g)
(14\cw^6-4g\cw^4+5g^2\cw^2-2g^3)}
{64\cw^5g^2(2\cw^2-g)},
&\displaystyle \frac{2g}{3}<\cw^2<2g,
\\[3mm]
-\displaystyle
\frac{\im(\cw^2+g)^2(3\cw^2-2g)
(\cw^6+16g\cw^4-11g^2\cw^2+2g^3)}
{32\cw^5g^2(2\cw^2-g)},
&\cw^2>2g,
\end{cases}
\\[4mm]
(\mathcal L_1f_2',g_1)_{\mathcal X}
&=
-\frac{\im(\cw^2-g)^2(\cw^2+g)^2
(3\cw^2-2g)(3\cw^2+2g)}
{64\cw^5g^2(2\cw^2-g)},
\\[3mm]
(\mathcal L_1f_1',g_2)_{\mathcal X}
&=
\frac{3\im(\cw^2-g)^2(\cw^2+g)^2
(\cw^2+2g)(3\cw^2+2g)}
{64\cw^5g^2(2\cw^2+g)},
\\[3mm]
(\mathcal L_1f_2',g_2)_{\mathcal X}
&=
\frac{\im(\cw^2+g)^2(\cw^2+2g)
(3\cw^6-8g\cw^4-17g^2\cw^2-2g^3)}
{32\cw^5g^2(2\cw^2+g)}.
\end{aligned}
\end{equation}
Combining \eqref{L2 f g} and \eqref{app:isola-L1-inner-products}, we obtain \eqref{isola:C21}, \eqref{isola:C22}, \eqref{isola:C11}.
\end{proof}

We finally introduce a useful lemma.

\begin{lemma}\label{final lemma}
For every sufficiently small $\delta$, one has
\[
\mathbf H_{12}(\bar\mu+\delta,0)
=\mathbf H_{21}(\bar\mu+\delta,0)=0,
\qquad
\partial_\e\mathbf H(\bar\mu+\delta,0)=0.
\]
\end{lemma}

\begin{proof}
Fix $\delta$ sufficiently small. The flat operator $\mathcal L_{\bar\mu+\delta,0}$ and its resolvent preserve each Fourier index. Hence so do $P(\delta,0)$, $P_0$, and the Kato maps $U(\delta,0)$ and $(U(\delta,0)^{-1})^*$. Thus $f_b(\delta,0)$ and $g_a(\delta,0)$ have Fourier indices $k_b$ and $k_a$, respectively. Since $k_2-k_1=2$, Fourier orthogonality gives $\mathbf H_{ab}(\bar\mu+\delta,0)=(\mathcal{L}_{\bar\mu+\delta,0}f_b(\delta,0),g_a(\delta,0))_{\mathcal{X}}=0$ for $a\ne b$.

The operator $\left.\partial_\e\mathcal L_{\bar\mu+\delta,\e}\right|_{\e=0}$ shifts Fourier indices only by $\pm1$. The resolvent differentiation formula implies the same property for $\partial_\e P(\delta,0)$. Differentiating the convergent power series defining $U$ inserts exactly one factor $\partial_\e P$ into each term; all remaining factors preserve Fourier indices. Consequently, $\partial_\e U(\delta,0)$ and, by $\partial_\e U^{-1}=-U^{-1}(\partial_\e U)U^{-1}$, $\partial_\e U^{-1}(\delta,0)$ also shift indices only by $\pm1$. It follows that $\partial_\e f_b(\delta,0)$ and $\partial_\e g_a(\delta,0)$ have indices only in $\{k_b-1,k_b+1\}$ and $\{k_a-1,k_a+1\}$, respectively.

Differentiating the defining pairing for $\mathbf H_{ab}$ therefore gives three terms
\begin{align*}
\pa_\e \mathbf{H}_{ab}(\bar\mu+\delta,0)=(\pa_\e \mathcal{L}_{\bar\mu+\delta,0}f_b(\delta,0),g_a(\delta,0))_\mathcal{X}+( \mathcal{L}_{\bar\mu+\delta,0}\pa_\e f_b(\delta,0),g_a(\delta,0))_\mathcal{X}+( \mathcal{L}_{\bar\mu+\delta,0}f_b(\delta,0),\pa_\e g_a(\delta,0))_\mathcal{X}
\end{align*}
pairing indices $k_b\pm1$ with $k_a$, or $k_b$ with $k_a\pm1$. Since $k_a-k_b\in\{0,\pm2\}$, all three terms vanish.
\end{proof}

\footnotesize

\vspace{1em}
\noindent{\bf Acknowledgments:}
T.-Y. Hsiao is supported by the European Union  ERC CONSOLIDATOR GRANT 2023 GUnDHam, Project Number: 101124921. Views and opinions expressed are however those of the authors only and do not necessarily reflect those of the European Union or the European Research Council. Neither the European Union nor the granting authority can be held responsible for them.

\vspace{1em}
\noindent{\bf Data availability:}
Data will be made available on reasonable request.

\vspace{1em}
\noindent{\bf Conflict of interest:}
The authors declare that there are no conflicts of interest regarding the publication of this paper.

\bibliographystyle{plain}
\bibliography{references.bib}
\addcontentsline{toc}{section}{References}
\end{document}